\documentclass[11pt,twoside,reqno]{amsart}
\usepackage{amssymb, amsmath, graphicx, amsfonts, euscript, mathrsfs}
\usepackage{color}
\usepackage{tikz}
\usepackage{bm}
\usepackage{dsfont}
\usepackage{comment}

\DeclareMathOperator{\Height}{h}
\newtheorem{thm}{Theorem}[section]

\newtheorem{prop}[thm]{Proposition}

\newtheorem*{Hk}{Hypothesis k}
\newtheorem{lemma}[thm]{Lemma}
\theoremstyle{remark}
\newtheorem{remark}{Remark}
\newtheorem{question}{Question}
\newtheorem{example}{Example}

\newcommand{\I}{\mathbb{I}}
\newcommand{\e}{\epsilon}

\newcommand{\A}{\alpha}
\newcommand{\Ad}{\mathbb{A}}

\newcommand{\om}{\omega}

\newcommand{\la}{\lambda}
\newcommand{\La}{\Lambda}

\newcommand{\R}{\mathbb{R}}
\newcommand{\Z}{\mathbb{Z}}
\newcommand{\N}{\mathbb{N}}
\newcommand{\LL}{\rm {Log}}
\newcommand{\F}{\mathcal{F}}

\newcommand{\Q}{\mathbb{Q}}

\newcommand{\G}{\mathcal{G}}
\newcommand{\les}{\lesssim}

\newcommand{\T}{\mathbb{T}}

\newcommand{\E}{\mathbb{E}}

\newcommand{\ka}{\kappa}

\newcommand{\V}{\mathbf{V}}
\newcommand{\ind}[1]{\mathds{1}_{{#1}}}

\newcommand{\D}{\mathbb{D}}

\newcommand{\ma}{\mathbb{A}}
\newcommand{\Cramer}{\operatorname{Cram\acute{e}r}}
\newcommand{\HB}{\operatorname{HB}}
\newcommand{\Log}{\operatorname{Log}}

\newcommand{\Sample}{{\mathcal{S}}}

\newcommand{\Schwartz}{{\bm S}}
\renewcommand{\mod}{{\ \mathrm{mod}\ }}
\def\beq{\begin{equation}}
\def\eeq{\end{equation}}
\numberwithin{equation}{section}

\begin{document}
\title[Pointwise convergence]{Pointwise convergence of
 polynomial  multiple \\ ergodic averages along the primes}

\author{Mariusz Mirek}
\address[Mariusz Mirek]{
Department of Mathematics,
Rutgers University,
Piscataway, NJ 08854-8019, USA
\&
Instytut Matematyczny,
Uniwersytet Wroc{\l}awski,
Plac Grunwaldzki 2/4, 
50-384 Wro\-c{\l}aw,
Poland}
\email{mariusz.mirek@rutgers.edu}

\author{Renhui Wan}

\address[Renhui Wan]{Ministry of Education Key Laboratory of NSLSCS,  School of Mathematical Sciences, Nanjing Normal University, Nanjing 210023, People's Republic of China}

\email{wrh@njnu.edu.cn}

\author{James Wright}
\address[James Wright]{School of Mathematics and Maxwell Institute for Mathematical Sciences, 
James Clerk Maxwell Building,
The King's Buildings,
Peter Guthrie Tait Road,
City Edinburgh,
EH9 3FD, UK}
\email{J.R.Wright@ed.ac.uk}

\begin{abstract}
We establish pointwise almost everywhere convergence  for  the   polynomial multilinear ergodic averages 
$$\frac{1}{N} \sum_{n=1}^N \La(n) f_1(T^{P_1(n)} x)\cdots f_k(T^{P_k(n)} x)$$
as $N\to \infty$, where $\La$ is the von Mangoldt function, $T \colon X \to X$ is an invertible  measure-preserving transformation of a probability space $(X,\nu)$, $P_1,\ldots, P_k$ are polynomials with integer coefficients and distinct degrees, and $f_1,\ldots,f_k\in L^\infty(X)$.
This pointwise almost everywhere convergence result can be seen as a refinement of  the norm convergence result obtained  in Wooley--Ziegler (Amer. J. Math, 2012) in the case of  polynomials with distinct degrees. 

We develop a  multilinear circle method for 
  von Mangoldt-weighted (equivalently, prime-weighted) averages in the general $k$-linear setting. 
The advantage of our method, besides establishing Weyl-type inequalities for multilinear Cram{\'e}r-weighted averages and sharp $p$-adic $L^q$-improving multilinear estimates among other tools, is that for the first time it allows us to work with inverse theorems having subpolynomial bounds in the general multilinear setting. This, in turn, yields sharp $r$-variational estimates $r > 2$ for our weighted polynomial multilinear ergodic average and, more importantly, offers prospects for addressing other multilinear problems involving inverse theorems lacking polynomial bounds.
\end{abstract}

\date{\today}

\maketitle

    

\section{Introduction}
\label{intrs1}
\subsection{Motivation and statement of the main results}
\label{Mr}
Multilinear polynomial ergodic averages lie at the heart of modern ergodic theory, and understanding their norm and almost everywhere behavior is a fundamental problem in the field. This began in the early 1930s with von Neumann's mean ergodic theorem \cite{vN} and Birkhoff's pointwise ergodic theorem \cite{BI}. In  Furstenberg's fundamental work \cite{Fur0}, multilinear ergodic averages for linear polynomials and a single transformation emerged as a natural tool for identifying recurrent points and consequently, arithmetic progressions within subsets of integers with positive upper density. This approach led to an ergodic theoretic proof of 
Szemer{\'e}di's famous theorem \cite{Sem1} and sparked the development of a new field, ergodic Ramsey theory. Over the past century, significant progress has been made in this area of research, which we now briefly summarize.

Let $X=(X,\nu)$ be a probability space and $T:X\to X$ is an invertible measure preserving map; that is, $\nu(T^{-1}(E))=\nu(E)$ for all measurable sets $E\subset X$.
The triple $X=(X,\nu,T)$ is referred to as a measure-preserving system. Given complex-valued  functions $f_1,\ldots,f_k \in L^\infty(X)$ and  polynomials $P_1,\ldots, P_k$ with integer coefficients, a scale $N \geq 1$, and a weight function $w \colon \N \to \mathbb C$, we define the weighted polynomial multilinear ergodic averages  by 
\begin{align}\label{eq:weightedweighted}
 A_{N, w; X}^{P_1,\ldots,P_k}(f_1,\ldots,f_k)(x):=\,\E_{n \in [N]} w(n) f_1(T^{P_1(n)} x)\cdots f_k(T^{P_k(n)} x),\qquad x\in X,
\end{align}
where $\E_{n\in S} f(n) = \frac{1}{\# S} \sum_{n\in S} f(n)$
 (see Subsection \ref{subsectionbasicnotation} for notation and basic definitions).  In this paper, we consider  the problem of  pointwise almost everywhere convergence of 
$A_{N, \La; X}^{P_1,\ldots,P_k}$, where  $\La$ is  the von Mangoldt function given  by 
\begin{align*}
\Lambda(n)=\begin{cases}
\log p,\quad n \textnormal{ is a power of a prime } p,\\
0,\quad \quad \,\,\,\textnormal{otherwise}.\end{cases}    
\end{align*}

\subsubsection{Unweighted ergodic average}
\label{subsubsection:Unweighted ergodic averages}
By setting  $w=1$ in (\ref{eq:weightedweighted}), we obtain the 
 unweighted ergodic averages  
\begin{align}\label{eq:unweighted}
A_{N, 1; X}^{P_1,\ldots,P_k}(f_1,\ldots,f_k)(x):=\E_{n \in [N]}  f_1(T^{P_1(n)} x)\cdots f_k(T^{P_k(n)} x).
\end{align}
Both norm and pointwise almost everywhere convergence of (\ref{eq:unweighted}) as $N\to \infty$ has been studied extensively. 

The $L^2(X)$ norm
 convergence of the averages   \eqref{eq:unweighted} and more general cases is well understood, largely due to the definitive work of Walsh \cite{W}.
 Before \cite{W}, Host and Kra \cite{HK05}, and independently Ziegler \cite{Z1}, established
 $L^2(X)$
convergence of \eqref{eq:unweighted} for linear polynomials. 
We refer the reader to  Leibman's appendix in \cite{Ber2}, which provides a comparative analysis of the characteristic factor constructions introduced by Host--Kra \cite{HK05} and Ziegler \cite{Z1}, showing their equivalence.  The case of general polynomials was addressed  by Leibman \cite{Le05}, 
Frantzikinakis--Kra \cite{FraKra}, and Host--Kra \cite{HK2}. In addition,
 the limiting function for \eqref{eq:unweighted} can be identified through the theory of Host--Kra factors \cite{HK05} and  equidistribution on nilmanifolds \cite{Le05}.
For the more general  averages  
\beq\label{commaveragesBerglesonopen}
\E_{n\in [N]} f_1(T_1^{P_1(n)}x)\cdots f_k(T_k^{P_k(n)}x),\qquad x\in X,
\eeq
where $T_1,\ldots,T_k:X\to X$ are    commuting invertible  measure-preserving transformations,  Chu--Frantzikinakis--Host  \cite{CFH} established $L^2(X)$ norm
convergence for \eqref{commaveragesBerglesonopen} under the assumption that the
polynomials $P_1,\ldots,P_k$ have distinct degrees;     we  refer to \cite{A2,H,Tao} for the case of linear polynomials.  Walsh \cite{W} finally established
norm convergence of~\eqref{commaveragesBerglesonopen} in general cases, including the case of noncommutative transformations $T_1,\ldots,T_k$ which generate  a nilpotent group. 
On the other hand, identifying the limit for general polynomial ergodic averages \eqref{commaveragesBerglesonopen} remains an open problem in ergodic theory.  Recently, significant progress was made by Frantzikinakis and Kuca \cite{FK25}, who identified the  
$L^2(X)$ limit for \eqref{commaveragesBerglesonopen}  in the case of commuting transformations and linearly independent polynomials.

Our current understanding of pointwise a.e. convergence for multilinear polynomial ergodic averages is less developed. Although pointwise a.e.  convergence is the most intuitive form of convergence, it is quite subtle and can differ markedly from norm convergence, which makes investigating pointwise convergence problems particularly challenging. By the dominated convergence theorem, pointwise a.e.
convergence of the averages \eqref{eq:unweighted} 
implies norm convergence for all
$f_1, \ldots, f_k\in L^{\infty}(X)$ on a probability space
$(X,  \nu)$. Furthermore, a maximal function bound for the averages \eqref{eq:unweighted} shows that pointwise a.e. convergence implies norm convergence for general Lebesgue space functions. Thus, pointwise a.e.
convergence can be viewed as a refinement or upgrade of norm convergence. 
The question of pointwise a.e. convergence is the celebrated Furstenberg--Bergelson--Leibman conjecture~\cite[Section 5.5]{BL02}.  In two  notable instances, \cite{Bour89} and \cite{Bou90}, Bourgain established the pointwise a.e. convergence of \eqref{eq:unweighted} in the linear ($k=1$) case, and in the blinear case ($k=2$) when $P_1,P_2$ are linear, respectively. 
Krause, the first author  and Tao   \cite{KMT22}, proved   pointwise a.e. convergence of \eqref{eq:unweighted} in 
the bilinear case with polynomials of the form $(n,P(n))$ where $\deg P\ge 2$.  
Very recently, Kosz,  Peluse and  the authors  \cite{KMPWW24} established pointwise a.e.  convergence of  \eqref{commaveragesBerglesonopen}
for arbitrary  $k\in \Z_+$ and   polynomials $P_1,\ldots, P_k$ with  distinct degrees. This result 
gives an affirmative  answer to a  question posed  by  Bergelson \cite[Question~9]{Ber1} in the case of polynomials with distinct degrees. 

\subsubsection{Prime-weighted ergodic averages}
\label{subsubsect:weightea}
The objective of this paper is to prove the pointwise a.e. convergence for the weighted polynomial multilinear ergodic averages $A_{N, \La; X}^{P_1,\ldots,P_k}$.  
The pointwise a.e. convergence for these weighted averages is equivalent to those for
the prime-weighted averages given by\footnote {All sums and products over the symbol $p$ will be understood to be over primes; other sums will be understood to be over positive integers unless otherwise specified.}
\beq\label{prime-weightednew}
 \frac{1}{N/\log N} \sum_{p \leq N} f_1(T^{P_1(p)}x)\cdots f_k(T^{P_k(p)}x),\qquad x\in X.
\eeq
Throughout this paper we will not distinguish between these two weighted averages.
The norm convergence of polynomial  
ergodic  prime-weighted averages \eqref{prime-weightednew} is known for general $k$ and arbitrary polynomials, thanks to the work of Wooley--Ziegler~\cite{wooley-ziegler} (see Frantzikinakis--Host--Kra~\cite{fhk} for the case of linear polynomials).

However, similar to the unweighted averages, the problem of pointwise
a.e. convergence of ergodic averages along the primes is more
involved.  The case of $k=1$ with linear polynomials was established
by Bourgain~\cite{bourgain-arithmetic} and Wierdl~\cite{wierdl} (with
the latter work allowing $L^q$ functions for any $q>1$), while the
case of an arbitrary (single) polynomial was handled by
Nair~\cite{nair1,nair2}.  We refer to \cite{Tro19, MT17, MS24} for
related work on quantitative forms of pointwise convergence for ergodic averages along the primes, including
variational, oscillation, and jump inequalities.  For the bilinear
case, by employing the method from Krause-Mirek-Tao \cite{KMT22} and
the generalized von Neumann theorem established by Ter\"{a}v\"{a}inen
\cite{Te24} (which handles a broad class of weight functions),
Krause--Mousavi--Tao--Ter\"{a}v\"{a}inen \cite{KMTT24} established the
pointwise a.e. convergence for polynomials of the form $(n,P(n))$ with
$\deg P\ge 2$. Recently, Fornal and Krause \cite{FK} proved pointwise a.e. convergence for $A_{N, \La; X}^{P_1,P_2}$ where   $P_1,P_2$ are liner polynomials with vanishing constant terms. Beyond these specific bilinear results, to the best of
our knowledge, no other pointwise a.e. convergence results exist for
the multilinear averages $A_{N, \La; X}^{P_1,\ldots,P_k}$ or
\eqref{prime-weightednew}.  We also mention that the problem of
pointwise a.e.  convergence of bilinear weighted ergodic averages was
discussed by Frantzikinakis \cite[Problem 12]{Fra16}.
 \subsubsection{Statement of our main result}
Motivated by previous work on 
\begin{itemize}
    \item[(a)] 
   the pointwise a.e. convergence of  multilinear unweighted ergodic averages \cite{KMT22,KMPWW24},

    \item[(b)] 
  the pointwise  a.e. convergence of  bilinear prime-weighted ergodic averages \cite{KMTT24},

  \item[(c)] 
  and the norm convergence of  multilinear prime-weighted ergodic averages \cite{wooley-ziegler},
    
\end{itemize}
we address the following question:
\begin{question}
\label{QQQ21}
Let $k\in \Z_+$ with $k\ge2$, and let $(X,\nu,T)$ be a measure-preserving system.\footnote{
 Throughout this paper, all measure-preserving systems are assumed to have finite measure, while all $\sigma$-finite measure-preserving systems will be explicitly stated as such.}  Let
$P_1,\ldots, P_k$ be polynomials with integer
coefficients. Is it true that for any functions $f_1, \ldots, f_k\in L^{\infty}(X)$ the
multilinear polynomial ergodic averages $A_{N,\La;X}^{P_1,\ldots, P_k}$ 
converge pointwise almost everywhere  on $X$ as $N\to \infty$?
\end{question}
This question extends \cite[Question~9]{Ber1} (in the single transformation case) posed by Bergelson to the prime-weighted setting, while  generalizing a question from Frantzikinakis's survey \cite[Problem 12]{Fra16} to general polynomial cases.
We now state the  main result of this paper.
\begin{thm}[Main result]\label{t1}
Let $k\in \Z_+$ and let $(X,\nu,T)$ be a $\sigma$-finite measure-preserving system.
 Suppose that ${\mathcal P} = \{P_1,\ldots, P_k\}$ is a family
of polynomials  with integer
coefficients and   distinct degrees. 
Let  $1<q_1,\ldots,q_k<\infty$ with 
$\frac{1}{q_1}+\cdots+\frac{1}{q_k}=\frac{1}{q}\le 1$. 
Then, for any   $f_1\in L^{q_1}(X),\ldots,f_k \in L^{q_k}(X)$, the following results hold.
\begin{itemize}
\item[\rm (i)] \textit{(Mean ergodic theorem)} The averages
$A_{N,\La;X}^{P_1,\ldots, P_k}(f_1,\ldots, f_k)$ converge in $L^q(X)$ as $N \to \infty$.

\smallskip

\item[\rm (ii)] \textit{(Pointwise ergodic theorem)} The  averages
$A_{N,\La;X}^{P_1,\ldots, P_k}(f_1,\ldots, f_k)$ converge pointwise a.e. as $N \to \infty$.

\smallskip

\item[\rm (iii)] \textit{(Maximal ergodic theorem)}
One has
\begin{align*}
\big\|\sup_{N\in\Z_+}\big|A_{N,\La;X}^{P_1,\ldots, P_k}(f_1,\ldots, f_k)\big|\big\|_{L^q(X)}\lesssim_{{\mathcal P},q_1,\ldots,q_k} \|f_1\|_{L^{q_1}(X)}\cdots\|f_k\|_{L^{q_k}(X)}.
\end{align*}

\smallskip

\item[\rm (iv)] \textit{(Variational ergodic theorem)}
If  $r>2$ and $\la>1$, then one has
\begin{align}
\label{eq:101}
\big\|\big(A_{N,\La;X}^{P_1,\ldots, P_k}(f_1,\ldots, f_k)\big)_{ N\in\D} \big\|_{L^{q}(X;\V^r)}\lesssim_{{\mathcal P},\la,r,q_1,\ldots,q_k} \|f_1\|_{L^{q_1}(X)}\cdots\|f_k\|_{L^{q_k}(X)} ,
\end{align}
where $\mathbb D = \{\lambda_n\in\N:n\in\N\} \subset [1,+\infty)$ is $\lambda$-lacunary, i.e., $\inf_{n\in\N}\frac{\lambda_{n+1}}{\lambda_n}\ge\lambda$.
\end{itemize}
\end{thm}

 We now make a few remarks on Theorems \ref{t1}.

\begin{enumerate}
  
\item 
The conclusion from part (ii)  provides an affirmative answer to {\it Question} \ref{QQQ21} when the polynomials have distinct degrees and establishes the pointwise a.e. convergence of 
 the prime-weighted averages  \eqref{prime-weightednew}. Additionally, with necessary adjustments, our method remains applicable to the unweighted averages \eqref{eq:unweighted}.

\smallskip

\item 
The conclusion from part (iv) in  Theorem \ref{t1}  serves as the central result of the theorem. From this, parts (ii) and (iii) follow directly, while part (i),
which was originally established by Wooley--Ziegler \cite{wooley-ziegler}, 
is immediately deduced via combining parts (ii)-(iii) and
the dominated convergence theorem. Furthermore, 
the variational  inequality \eqref{eq:101} removes the  restrictions that $k=2$ and that the  
 polynomials have the special form   $(n,P(n))$ with  $\deg P\ge 2$ required in Krause--Mousavi--Tao--Ter\"{a}v\"{a}inen \cite{KMTT24}.

\smallskip

\item    
 The range 
$r>2$ is  sharp,
as variational inequality \eqref{eq:101} fails for $r\le 2$, except for
certain specific operators (see, for example, \cite{DMT}). Moreover,
following the arguments in \cite[Section 11]{KMT22}, Theorem \ref{t1}
also holds for a certain range of exponents $1<q_1,\ldots, q_k<\infty$
such that $\frac{1}{q_1}+\cdots+\frac{1}{q_k}=\frac{1}{q}>1$, at the
cost of slightly increasing $r$ in \eqref{eq:101}.

\smallskip

 \item

 As noted earlier, the linear case ($k=1$) has been established in prior works \cite{Tro19, MT17, MS24}.
 
\end{enumerate}

\subsection{Reduction to the integer shift system}
\label{subsubsection:Trantoshift}
In the study of pointwise convergence problems for ergodic averages with polynomial iterates or other arithmetic features, the integer shift system emerges as the most significant dynamical system, as exemplified below.
\begin{example}[integer shift system]\label{ex:1}
    The integer shift system $(\Z,\nu_\Z,T_\Z)$ is the set of integers $\Z$ equipped with counting measure $\nu_\Z$ and the shift $T_\Z(x):= x-1$. 
  The averages
  $A_{N,\La; X}^{P_1, \ldots, P_k}$ with
  $T=T_\Z$ and $X=\Z$ are
\beq\label{moz1}
A_{N,\La;\Z}^{P_1,\ldots,P_k}(f_1,\ldots,f_k)(x):=\E_{n\in [N]}\La(n)
f_1(x-P_1(n))\cdots f_k(x-P_k(n)),\qquad x\in\Z.
\eeq
The associated   truncated version are given by 
\beq\label{moz2}
\tilde{A}_{N,\La;\Z}^{P_1,\ldots,P_k}(f_1,\ldots,f_k)(x):=\frac{1}{\lfloor N\rfloor}\sum_{n\in J_N}\La(n)  f_1(x-P_1(n))\cdots f_k(x-P_k(n)),\quad x\in\Z,
\eeq  
where $J_N:=[N]\setminus[N/2]$.  We will often abbreviate
$$A_{N,\La}^{{\mathcal P}}:=A_{N,\La; \Z}^{P_1, \ldots, P_k}\qquad {\rm and}\qquad \tilde{A}_{N,\La}^{{\mathcal P}}:=\tilde{A}_{N,\La; \Z}^{P_1, \ldots, P_k}.$$
\end{example}

In view of the Calder{\'o}n transference principle \cite{C1} (or, more
precisely, the arguments from \cite[Proposition~3.2 (ii)]{KMT22} or
\cite[Theorem~1.6]{Kosz}), to prove part (iv) in Theorem \ref{t1}, it
suffices to work with the integer shift system. In fact, the
Calder{\'o}n transference principle allows us to transfer the
quantitative estimates \eqref{eq:101} from the integer shift system to
corresponding estimates for $A_{N,\La; X}^{P_1,\ldots,P_k}$ in the
abstract ($\sigma$-finite) measure-preserving system $(X,\nu,T)$. It is important to
emphasize that this transference is only possible because the
parameters \( 1 < q_1, \ldots, q_k < \infty \) in Theorem \ref{t1}
satisfy the H{\"o}lder condition
\(\frac{1}{q_1} + \cdots + \frac{1}{q_k} = \frac{1}{q} \le 1\). Otherwise,
the transference does not hold. By working over the
integers,  we can employ Fourier methods on $\Z$ and utilize the
algebraic structure of $\Z$, which are generally not available in
abstract measure-preserving systems. It is also worth noting that the
Calder{\'o}n transference principle \cite{C1} only transfers
quantitative bounds that imply pointwise a.e. convergence, but does
not transfer pointwise a.e. convergence itself. Therefore, we will
focus on proving quantitative bounds for $A_{N,\La}^{\mathcal P}$ or
$\tilde A_{N,\La}^{\mathcal P}$ in the integer shift system, rather
than pointwise convergence on $\Z$. For technical reasons, we will
restrict our attention to the truncated averages
$\tilde A_{N,\La}^{\mathcal P}$.

After these reductions, Theorem \ref{t1} follows from the next theorem.
\begin{thm}\label{t3shifted}
Let $k\in \Z_+$ with $k\ge2$, and  let
${\mathcal P} = \{P_1,\ldots, P_k\}$ be a family
of polynomials  with integer
coefficients and   distinct degrees. 
Let  $1<q_1,\ldots,q_k<\infty$ with 
$\frac{1}{q_1}+\cdots+\frac{1}{q_k}=\frac{1}{q}\le 1$. 
Then, for any   $f_1\in \ell^{q_1}(\Z),\ldots,f_k \in \ell^{q_k}(\Z)$, 
we have
\beq\label{shift}
\begin{aligned}
\big\|\big(\tilde A_{N,\La}^{\mathcal P}(f_1,\ldots, f_k)\big)_{ N\in\D} \big\|_{\ell^{q}(\Z;\V^r)}\lesssim_{\mathcal P,\la,r,q_{1},\ldots,q_{k}} 
\|f_1\|_{\ell^{q_1}(\Z)}\cdots \|f_k\|_{\ell^{q_k}(\Z)}
\end{aligned}
\eeq
where $\tilde A_{N,\La}^{\mathcal P}$ is given by \eqref{moz2}, and $\mathbb D\subset [1,+\infty)$ is
$\lambda$-lacunary.
\end{thm}
Following the   arguments yielding \cite[Proposition~3.2 (ii) and (iii)]{KMT22}, we can obtain Theorem \ref{t1} from Theorem \ref{t3shifted}. Thus,
it remains to prove Theorem \ref{t3shifted}, which is the main objective of the remainder of the paper.

\subsection{Overview of the proof}
\label{overview proof}
Since the averages considered in this paper are multilinear, the
classical circle method for the linear case cannot be applied, as
discussed in \cite{KMT22,KMPWW24}.  In fact, our proof primarily
relies on the framework of adelic harmonic analysis developed in
\cite{KMT22}, while also drawing inspiration from \cite{KMPWW24} to
overcome the bilinear barrier discussed in \cite{KMTT24}. However,
when we consider the pointwise convergence problem of multilinear
averages with the von Mangoldt weight, we encounter new difficulties
that we now highlight.
\subsubsection{Multilinear circle method for unweighted averages}
\label{subsubsub MCMFUA}
We now start by discussing some limitations of the multilinear circle
methods presented in \cite{KMT22} and \cite{KMPWW24}, which deal with
the unweighted averages $A_{N, 1; X}^{P_1,\ldots,P_k}$.  For specific
details about their methods, we refer to \cite[Subsection 1.4]{KMT22}
and \cite[Subsection 1.6]{KMPWW24}.

\begin{enumerate}

\item 

The method from \cite{KMT22} is restricted to the bilinear operators,
the single transformation case, and polynomials of the form $(n,P(n))$
where $\deg P\ge2$.  These limitations are necessary for various parts
of their proof, such as deriving an appropriate version of Peluse's
inverse theorems, establishing the bilinear Rademacher-Menshov
inequality, transitioning variational inequalities from integers to
adelic integers, and proving a desired bilinear arithmetic estimate
(see \cite[Theorem 9.9]{KMT22}) based on the trivial (yet significant)
inequality \beq\label{linear necessary}
\|A_{\Z_p}^{n,P(n)}\|_{L^1(\Z_p)\times L^\infty(\Z_p)\to
L^\infty(\Z_p)}\le1.  \eeq See Subsection \ref{abstraction HA} below
for the definition of the $p$-adic integers $\Z_p$ and other related
concepts.  It is noteworthy that the linear polynomial in
\eqref{linear necessary} is indispensable for this argument.  Although
the method from \cite{KMT22} has its limitations, its main advantage
is that the decay $2^{-cl}+\langle \Log N \rangle^{-C_1}$ (where $c$
is small and $C_1$ is sufficiently large) for the minor arcs estimates
(see the bilinear Weyl inequality in \cite[Theorem 5.12]{KMT22}) is
sufficient and one does not require stronger polynomially decaying bounds as in
\cite{KMPWW24}.  \smallskip

  \item
The approach introduced in \cite{KMPWW24} overcomes the limitations previously identified in \cite{KMT22}. Crucially, this approach includes  minor arc estimates (see \cite[Theorem 6.1]{KMPWW24}), which exhibit a polynomial decay of the form 
 $2^{-cl}+N^{-c}$ (with small $c$).   These estimates enable the effective application of a multiparamater norm interchanging inequality, as detailed in \cite[Subsection 7.3.6 and Remark 7.99]{KMPWW24}.
However, due to the intricate nature of prime numbers, the weighted averaging operator considered in this paper does not exhibit such favorable decay (precisely, we obtain bounds of the form $2^{-cl}+\exp(-c\Log^{1/C_0}N)$  with small $c$ and large $C_0$, see Theorem \ref{Thm:uweyl1}),  rendering the approach in \cite{KMPWW24} directly inapplicable. 

\end{enumerate}
 \subsubsection{Multilinear circle method for  averages with von Mangoldt weight}
 To our knowledge, the general multilinear circle method, which
 addresses the weighted averages $A_{N,\La; X}^{P_1,\ldots,P_k}$,
 remains unexplored, aside from the bilinear cases in \cite{KMTT24, FK}.
 In particular, the limitations of the method in \cite{KMT22} also
 apply to the method in \cite{KMTT24}, as the latter is heavily
 reliant on the former.  Based on the discussion in Subsection
 \ref{subsubsub MCMFUA}, we cannot directly utilize the two methods
 mentioned. Instead, we will develop a multilinear circle method for
 the weighted averages $A_{N,\La; X}^{P_1,\ldots,P_k}$, building on
 the approaches in \cite{KMT22, KMPWW24, KMTT24} while introducing
 some striking novelties.  Precisely, the tools used to prove Theorem
 \ref{t3shifted} make up our multilinear circle method.
 \smallskip
 
 We first highlight these new features.
\begin{enumerate}

\item We establish an inverse theorem (see Theorem \ref{inverse2}) and
a version of Weyl's inequality (we call the multilinear Weyl
inequality) for multilinear Cram\'{e}r-weighted averages (see Section
\ref{section:MWIII}) involving general polynomials with distinct
degrees (see Theorem \ref{Thm:uweyl1}).  It is important to emphasize
that the inverse theorem (see Theorem \ref{inverse2}), due to its
arithmetic nature arising from the prime numbers, involves
subpolynomial bounds, as stated in condition
\eqref{lcondition1}. This, in turn, implies that the bound in our
multilinear Weyl inequality behaves like
\( 2^{-cl} + \exp(-c \log^{1/C_0} N) \) with small \( c \) and large
\( C_0 \), as seen in \eqref{uweyl11}. 
This is a
disadvantage that prevents us from applying the methods from \cite{P2,
KMPWW24}, which are aligned with the inverse theorems involving 
polynomial bounds.
Nevertheless, the latter bound
provides the desired control on the minor arcs.

\smallskip

\item Additionally, in the specific bilinear case involving the
polynomial pair \( (P_1, P_2) \), where \( P_1 \) is linear and
\( P_2 \) has a degree of at least \( 2 \), the decay
\( 2^{-cl} + \exp(-c \log^{1/C_0} N) \) with small \( c \) and large
\( C_0 \) achieved in our multilinear Weyl inequality improves upon
the decay \( 2^{-cl} + \langle \log N \rangle^{-C_1} \) with small
\( c \) and large \( C_1 \), obtained in \cite[Theorem 3.2]{KMTT24}.
This improvement is possible by applying the new Ionescu-Wainger
multiplier theorem for the set of canonical fractions established in
\cite{KMPWW24}.

\smallskip

\item The previous remark is of significant importance in other
prospective multilinear problems involving inverse theorems lacking
polynomial bounds. At this moment, it is believed that addressing
pointwise convergence problems for multilinear averages involving
general polynomial iterates, in the best possible scenario, will rely
on inverse theorems involving subpolynomial bounds.  From this
perspective, the weighted multilinear polynomial averages in this paper serve as toy
models to identify these difficulties and develop new robust tools
beyond additive combinatorics to understand the new challenges that
arise with inverse theorems involving subpolynomial bounds.

\smallskip
\item In order to break the bilinear barrier from \cite{KMTT24} and to
compensate for the lack of polynomial bounds in our multilinear Weyl
inequality, we adopt the approach of low/high frequencies from
\cite{KMPWW24}. The high-frequencies are handled in much the same way
as in \cite{KMPWW24}.  For the low-frequency regime, we proceed
analogously to the approach in \cite{KMT22, KMTT24}, developing a
multilinear Rademacher-Menshov inequality (see Lemma \ref{mil})
together with an arithmetic multilinear estimate (see Theorem
\ref{lend1}). The former is employed to prove the major arc estimates
for the low-frequency case at a small scale, while the latter
addresses the major arc estimates for the low-frequency case at a
large scale.

\smallskip
\item The key mechanism that permits handling the low-frequency case
at a large scale and compensates for the lack of polynomial bounds in
our multilinear Weyl inequality is an
$L^{\mathrm{p}}$-improving\footnote{In order to avoid confusion with
the prime $p$, we use the notation ${\rm p}$ throughout this paper.}
inequality for multilinear unweighted averages on the $p$-adic fields
$\mathbb{Q}_{p}$ for all primes $p$ (see Section \ref{section:local
field}).  An important feature of our $L^{\mathrm{p}}$-improving
multilinear inequality (see inequality \eqref{bound} in Theorem
\ref{thm:padiccccc}) is that it is proved in the sharp range of
exponent $q$; see Theorem
\ref{thm:padiccccc} and Theorem \ref{thm:padic}. This plays a key role
in establishing the  arithmetic multilinear estimate
in Theorem \ref{lend1}.

\smallskip
\item We finally emphasize that Theorem
\ref{thm:padiccccc} is a new and important tool that makes it
realistic to expect that inverse theorems involving subpolynomial
bounds may be sufficient for handling pointwise convergence problems
for multilinear averages with polynomial iterates. We plan to
investigate this direction further in the near future.
 \end{enumerate}

 \subsubsection{Outline the proof of Theorem \ref{t3shifted}}
 Utilizing the little Gowers norm estimate (refer to
 \eqref{eqn:La-LaN} below) and the generalized von Neumann theorem
 (see Lemma \ref{VNT}), we reduce the proof of Theorem \ref{t3shifted}
 to demonstrating Theorem \ref{Thm:uweyl1} and Theorem
 \ref{lemma:majorarcsest}, which provide the minor arc estimates and
 major arc estimates for the Cram\'{e}r-weighted averages,
 respectively.

\smallskip

\textit{Sketch of the proof of Theorem \ref{Thm:uweyl1}:}
We prove the
Cram\'{e}r-weighted inverse theorem (see Theorem \ref{inverse2}) and
subsequently prove Theorem \ref{Thm:uweyl1} by combining Theorem
\ref{inverse2}, the weighted $L^{\rm p}$-improving estimates (see
Lemma \ref{lp-improvvv}), and the arguments that address unweighted
averages, such as the Hahn-Banach theorem (see \cite[Lemma
6.9]{KMT22}), dual arguments, and the Ionescu-Wainger multiplier
theorem for the set of canonical fractions (see Theorem
\ref{thmIW}). The latter simplifies notation and allows us to work
with the naive height function (see Subsection \ref{subbb:Terminor}) in place of the more complicated
Ionescu--Wainger height function \cite{KMT22, KMTT24}.

\smallskip To obtain this Cram\'{e}r-weighted inverse theorem (see
Subsection \ref{winver2}), we combine a reduction argument with an
unweighted inverse theorem. The proof needs several key techniques,
including the alternate inverse theorem in \cite{KMT22}, the Peluse's
inverse theorem \cite[Theorem 3.3]{P2}, and the little Gowers norm
estimates related to the Cram\'{e}r and Heath-Brown approximants (see
Subsection \ref{subsection:approx}).  Notably, the subpolynomial lower bound of
$\delta$ (see condition \eqref{lcondition1}) in the Cramér-weighted inverse
theorem will determine the decay in the multilinear Weyl inequality
\eqref{uweyl11}, but the polynomial decay required in \cite{KMPWW24}
is not achieved under our current framework even though we rely on the Ionescu-Wainger multiplier
theorem for the set of canonical fractions.

\smallskip
\textit{Sketch of the proof of Theorem  \ref{lemma:majorarcsest}:}
We employ both arithmetic and continuous dyadic decompositions to reduce the proof of Theorem \ref{lemma:majorarcsest} to demonstrating Propositions \ref{prop:highfre}, \ref{t4.1}, and \ref{t4.2}, which correspond to the high-frequency case, the low-frequency case at a small scale, and the low-frequency case at a large scale, respectively.
\smallskip

Proposition \ref{prop:highfre} (the high-frequency case) will be established by proving a multilinear Weyl inequality in the continuous setting 
given by 
Theorem \ref{appen11}. This is achieved by using the methods from \cite[Section 7]{KMPWfield}. 
In the proof of the low-frequency case, the multiparamater norm interchanging inequality from \cite{KMPWW24} fails to apply because the multilinear Weyl inequality in this setting lacks the required polynomial decay. To overcome this, we develop a nontrivial variant of the method introduced in \cite{KMT22}.
\smallskip

  The proof of Proposition \ref{t4.1} relies on establishing a multilinear Rademacher–Menshov inequality (see Lemma \ref{mil}), which we derive by combining an inductive argument with the bilinear version obtained in \cite{KMT22}.
  For the proof of Proposition \ref{t4.2}, in addition to applying the quantitative Shannon sampling theorem (Theorem  \ref{Sampling}) and the norm interchange technique (see \eqref{interchange}), we establish two new estimates: an arithmetic multilinear estimate (Theorem \ref{lend1}) and a $p$-adic multilinear estimate (Theorem \ref{thm:padic}).

\subsection{Organization} In Section \ref{sect.preparation}, we
summarize the notation we will use in the paper and compile several
key theorems and lemmas from the literature.  Section
\ref{section:MWIII} presents a multilinear Weyl inequality for
Cram\'er-weighted averages which yields the desired minor arc
estimates. In Section \ref{section:mainreduction}, we reduce the proof
of Theorem \ref{t3shifted} to showing Propositions \ref{prop:highfre},
\ref{t4.1} and \ref{t4.2}.  These propositions correspond to the major
arc estimates for the high-frequency case, the low-frequency case at
a small scale, and the low-frequency case at a large scale
respectively. We then proceed to prove Proposition \ref{prop:highfre}.
In Sections \ref{section:small scale} and \ref{section:large scale
Adelic group}, we establish Propositions \ref{t4.1} and \ref{t4.2}.
In Section \ref{section:local field}, we prove Theorem
\ref{thm:padiccccc} giving the $L^{\rm p}$-improving inequality for
multilinear unweighted averages on $p$-adic fields.  Appendix \ref{section:Appen1}
provides a multilinear Weyl inequality in the continuous setting,
which is utilized in the proof of Proposition
\ref{prop:highfre}. Appendix \ref{appenmap} includes the definition of
the sampling map $\Sample$ and the quantitative Shannon sampling
theorem, both of which are instrumental in proving Proposition
\ref{t4.2}. In Appendix \ref{Appendix:prop-1}, we prove Proposition
\ref{prop-1} which is essential to establish Theorem
\ref{thm:padiccccc}. Proposition \ref{prop-1} is a key structural
result for polynomials over the $p$-adic fields, identifying regions
on which these polynomials are injective.

\section*{Acknowledgments} 
Mariusz Mirek was partially supported  by NSF CAREER grant DMS-2236493. James Wright was partially supported by a Leverhulme Research Fellowship RF-2023-709$\backslash$9.

\section{Notation and preliminaries}
\label{sect.preparation}
\subsection{Basic notation}
\label{subsectionbasicnotation}
For any two quantities $A$ and $B$, we will write $A\lesssim B$  to denote   $A\le C B$ for some absolute constant $C$. The notation $A=B+{O}(X)$ means
   $|A-B|\les X$.
   If both $A\lesssim B$ and $B\les A$ hold, we use $A\sim B$.
   If the implied constant $C$
depends on additional parameters, we will denote this by subscripts. 
  To abbreviate the notation, we will sometimes suppress
  the dependence on certain parameters in the implied constant when the issue of uniformity with respect to such parameters is not important. 
 
 We denote the prime numbers by $\mathbb P:=\{2,3,5,\ldots\}$, 
 the positive integers by $\Z_+:=\{1,2,\dots\}$ and the natural numbers by  $\N:=\Z_+\cup\{0\}$.  
 For any $x\in \R$, $\lfloor x\rfloor$ denotes the largest integer not larger than $x$.
 For any $N > 0$, we use $[N]$  to denote the discrete interval $\{ n \in \Z_+: n \leq N \}$,
 let $\N_{\le N}:=[N]\cup \{0\}$  
  and $[\pm N]:=[-N, N]\cap \Z$.   For  $q\in\Z_+$ and $a = (a_1,\ldots, a_k) \in \Z^k$ with $k\in \Z_+$, we denote by $(a,q)$ the greatest
common divisor of $a$ and $q$; that is, the largest  $d\in\Z_+$ that divides $q$ and all the components
$a_1, \ldots, a_k$. 
Clearly, any vector in $\Q^k$ (the direct product of $k$ copies of the rational numbers $\Q$) has a unique representation as $a/q$ with $q\in \Z_{+}$,
$a \in \Z^k$ and $(a,q)=1$.  We let $[q]^\times:= \{ b \in [q]: (b,q) = 1 \}$ denote the elements of $[q]$ that are coprime to $q$. Additionally, 
if $\rm R$ is a commutative ring, we use $\rm R^\times$ to denote the multiplicative group of invertible elements of $\rm R$.

We use the averaging notation
$
 \E_{n \in S} f(n):= \frac{1}{\# S} \sum_{n \in S} f(n)
$
for any finite non-empty set $S$, where $\# S$ denotes the cardinality of $S$.  
We also need  
the Euler totient function by  $\varphi(q):= \# ( [q]^\times)$ for any $q\in \Z_+$.
  Moreover, $\ind{E}$ denotes the indicator function of a set $E$. If $P$ is a statement, we use $\ind{P}$ to denote the indicator function, equal to 1 if $P$ is true and 0 if $P$ is false; in particular, $\ind{E}(x) = \ind{x \in E}$.  We use the Japanese bracket notation $ \langle x \rangle := (1 + |x|^2)^{1/2}$ for any real or complex $x$.
  All logarithms in this paper are taken base $2$, and for any $N \geq 1$, we define the \emph{logarithmic scale} $\Log N$ of $N$ by the formula
\begin{equation*}
 \Log N:= \lfloor \log N  \rfloor,
 \end{equation*}
thus $\Log N$ is the unique natural number such that $2^{\Log N} \leq N < 2^{\Log N+1}$. Throughout this paper we fix a cutoff function $\eta \colon \R \to [0,1]$ that is a smooth even function supported on $[-1,1]$ that equals one on $[-1/2,1/2]$.    For any $\mathrm k \in \Z$, we define  $\eta_{\leq \mathrm k } \colon \R \to [0,1]$ by
\begin{equation}\label{eta-resc}
\eta_{\leq \mathrm k }(\xi):= \eta(\xi/2^{\mathrm k} ).
\end{equation}
\subsection{Functional spaces}
All vector spaces considered in this paper are defined over the complex field  $\mathbb C$. 
Given a measure space  $(X,\nu)$, we define $L^0(X)$ as the space of all $\nu$-measurable
complex-valued functions on $X$, identifying functions that are equal almost everywhere with respect to $\nu$. For
$q\in(0, \infty)$, the subspace of functions in $L^0(X)$ whose absolute value raised to the
$q$-th power is integrable is denoted by  $L^q(X)$;  the space $L^{\infty}(X)$ consists of all essentially bounded functions in $L^0(X)$. For any exponent $q\in [1,\infty]$, 
its H\"older conjugate $q'\in [1,\infty]$  is defined by the relation  $1/q + 1/q' = 1$.  When the measure on $X$ is counting measure, we  abbreviate $L^q(X)$ as $\ell^q(X)$, or simply $\ell^q$.

We can extend these notation  to functions that take values in a finite-dimensional normed vector space $V = (V, \|\cdot\|_V)$. For example,  $L^0(X;V)$ refers to the space of measurable functions from $X$ to $V$ (up to almost everywhere equivalence), as well as 
\begin{align}
\label{eq:11}
L^{p}(X;V)
:=\big\{F\in L^0(X;V):\|F\|_{L^{q}(X;V)}:= \big\|\|F\|_V\big\|_{L^{q}(X)}<\infty\big\}.
\end{align}
These notions can also be extended to functions taking values in an infinite-dimensional normed vector space $V$, at least when  $V$
 is separable. However, in most cases we will be able to work within finite-dimensional settings, or reduce to them via a standard approximation argument.

For any finite dimensional normed vector space $(B,\|\cdot\|_B)$,  any sequence
 $(\mathfrak a_t)_{t\in\I}$ of elements of $B$ indexed by a totally
 ordered set $\I$, and any exponent $1 \leq r < \infty$, the
 $r$-variation seminorm  is defined by the following formula:
 \begin{equation}\label{var-seminorm}
  \| (\mathfrak a_t)_{t \in \I} \|_{V^r(\I; B)}:=
 \sup_{J\in\Z_+} \sup_{\substack{t_{0} < \dotsb < t_{J}}}
\big(\sum_{j=0}^{J-1}  \|\mathfrak a_{t_{j+1}}-\mathfrak a_{t_{j}}\|_B^{r} \big)^{1/r},
 \end{equation}
where the  supremum is taken over all finite increasing sequences $\{t_{j}\}_{j\in [J]}\subset\I$, and it is set to zero by convention if $\I$ is empty.  
When $r=\infty$,  we define 
$
 \| (\mathfrak a_t)_{t \in \I} \|_{V^\infty(\I;B)}:= \sup_{t \leq t'\in \I} \|\mathfrak a(t') - \mathfrak a(t)\|_B.
$

The $r$-variation norm for $1 \leq r \leq \infty$ is defined by 
\begin{equation}\label{vardef}
  \| (\mathfrak a_t)_{t \in \I} \|_{\V^r(\I; B)} 
:= \sup_{t\in\I}\|\mathfrak a_t\|_B+
\| (\mathfrak a_t)_{t \in \I} \|_{V^r(\I;B)}.
\end{equation}
This clearly defines a norm on the space of functions mapping from $\I$ to $B$.
When $B=\mathbb C$,  we will denote $V^r(\I;B)$ by $V^r(\I)$ or simply $V^r$, and $\V^r(\I;B)$ by $ \V^r(\I)$ or $\V^r$. If $(X,\mu)$ is a measure space, by \eqref{vardef} and  \eqref{eq:11}, we write
\[
L^q(X;\V^r)=\big\{F\in L^0(X;\V^r):\|F\|_{L^{q}(X;\V^r)} :=\left\|\|F\|_{\V^r}\right\|_{L^{q}(X)}<\infty\big\}.
\]
Note that the $\V^r$ norm is non-increasing in $r$, and comparable to the $\ell^\infty$ norm when $r=\infty$. We also observe the simple triangle inequality
\begin{equation}\label{simple}
\| (\mathfrak a_t)_{t \in \I} \|_{\V^r(\I;X)} \lesssim \| (\mathfrak a_t)_{t \in \I_1} \|_{\V^r(\I_1;X)} + \| (\mathfrak a_t)_{t \in \I_2} \|_{\V^r(\I_2;X)} 
\end{equation}
whenever $\I = \I_1 \uplus \I_2$ is an ordered partition of $\I$.  In a similar spirit,  we have the bound
\begin{equation}\label{varsum}
 \| (\mathfrak a_t)_{t \in \I} \|_{\V^r(\I;X)} \lesssim \| (\mathfrak a_t)_{t \in \I} \|_{\ell^r(\I;X)}. 
 \end{equation}

\subsection{Abstract harmonic analysis}\label{abstraction HA}
Throughout this paper, we use Fourier analysis on multiple groups, with special emphasis on the intricate connection between major arc Fourier analysis on the integer group
$\Z$ and low frequency Fourier analysis on the adelic integers $\Ad_\Z$. To carry out this analysis systematically, we establish some abstract harmonic analysis notation. We define the unit circle as $\T:= \R/\Z$ and introduce the standard character $e: \T \to \mathbb{C}$ given by $e(\theta) := e^{-2\pi {\rm i} \theta}$, where ${\rm i}^2=-1$.

\subsubsection{Fourier transform}
Let $(\G, +)$ be a locally compact abelian (LCA) group equipped with a Haar measure $\nu_{\G}$. It is well known (see for instance \cite{F16,Rudin}) that every LCA group $\G$ has a Pontryagin dual ${\G}^* = ({\G}^*,+)$, 
a LCA group  with a Haar measure $\nu_{\G^*}$ and a pairing, i.e.,~a continuous bihomomorphism $\G \times \G^*\ni (x,\xi) \mapsto x\cdot \xi\in\T$, such that the Fourier transform $\F_{\G} : L^1(\G) \to C(\G^*)$ given by
\[
 \F_{\G} f(\xi):= \int_{\G} f(x) e(x\cdot \xi)d\nu_{\G}(x), \qquad \xi\in\G^*,
\]
extends to a unitary map from $L^2(\G)$ to $L^2(\G^*)$; in particular we have Plancherel's identity
$
\|\F_\G f\|_{L^2(\G^*)}=\|f\|_{L^2(\G)}$ for $f\in L^2(\G)$.
Moreover, the inverse Fourier transform $\F_\G^{-1} : L^2(\G^*) \to L^2(\G)$ is given by the formula
\[
 \F_\G^{-1} f(x) = \int_{\G^*} f(\xi) e(-x\cdot \xi)d\nu_{\G^*}(\xi), \qquad {\rm where}\quad f\in L^1(\G^*) \cap L^2(\G^*).
\]
We will utilize the following concrete pairs $(\G,\G^*)$ of Pontryagin dual LCA groups:
\begin{itemize}
    \item   $(\G,\G^*)=(\R,\R)$.
    \item   $(\G,\G^*)=(\Z,\T)$.
    \item   $(\G,\G^*)=(\Z/Q\Z,\frac{1}{Q} \Z/\Z)$ for some $Q \in \Z_+$.
     \item   $(\G,\G^*)=(\Z_p:= \varprojlim_j \Z/p^j\Z,\ \ \Z_p^*:= \Z[\frac{1}{p}]/\Z )$ for some prime $p \in \mathbb{P}$.
    \item  $(\G,\G^*)=(\hat \Z:= \prod_{p \in \mathbb{P}} \Z_p,\ \ \Q/\Z )$.
    \item   $(\G,\G^*)=(\Ad_\Z:= \R \times \hat \Z,\ \ \Ad_\Z^*:= \R \times \Q/\Z)$.
    More generally, if $\G_1,\ldots, \G_k$ are LCA groups with Pontryagin duals $\G^*_1,\ldots, \G^*_k$, then the product $\G=\G_1 \times\cdots \times \G_k$ with product Haar measure $\nu_{\G_1}\times\cdots\times\nu_{\G_k}$ is a LCA group with Pontryagin dual ${\G}^*=\G^*_1 \times\cdots\times \G^*_k$ with product Haar measure $\nu_{\G^*}=\nu_{\G^*_1}\times\cdots\times\nu_{\G^*_k}$.
    \end{itemize}
For further details, we refer to \cite[Section 4]{KMT22}.  We consider the space $\Schwartz(\G) \subset \big(L^1(\G) \cap L^\infty(\G)\big)$ of Schwartz--Bruhat functions $f \colon G \to \mathbb{C}$, which generalizes the classical Schwartz functions on  $\R$.  These functions are dense in $L^q(\G)$ for every $1 \leq q < \infty$ 
 and are well-suited for Fourier analysis. A detailed definition of the space $\Schwartz(\G)$ 
 can be found in \cite[Section 4]{KMT22}. For the general case,  see also \cite{bruhat,osb}. 
  For  $k\in\Z_+$, we denote by   $(\Schwartz(\G))^k$   the direct product of $k$  copies of
space $\Schwartz(\G)$.

\subsubsection{Fourier multiplier operators}
We now define Fourier multiplier operators.  A continuous function $\mathfrak m \colon \G^* \to \mathbb C$ is said to be \emph{smooth tempered} if $\mathfrak m F \in \Schwartz(\G^*)$ whenever $F \in \Schwartz(\G^*)$.  

If $\mathfrak m\colon \G^* \to \mathbb C$ is a smooth tempered function,
 we define the linear Fourier multiplier operator $T_\G[\mathfrak m]\colon \Schwartz(\G)  \to \Schwartz(\G)$  by 
\begin{align}
\label{eq:100}
T_{\G}[\mathfrak m]f(x) := \int_{\G^*}e(-x\cdot \xi)\mathfrak m(\xi)\mathcal F_{\G}f(\xi)d\nu_{\G^*}(\xi), \qquad  x\in\G.
\end{align}
We refer to $\mathfrak m$ as the \emph{symbol} of $T_{\G}[\mathfrak m]$.  
In addition, for any finite subset  $\mathfrak{Q}\subseteq \G^*$ we define
\begin{align}
\label{eq:106}
T_{\G}^{\mathfrak{Q}}[\mathfrak m]f(x) := T_{\G}\big[\sum_{\theta\in \mathfrak{Q}}\tau_{\theta}\mathfrak m\big]f(x), \qquad  x\in\G,
\end{align}
where $\tau_{\theta}\mathfrak m(\xi) := \mathfrak m(\xi-\theta)$ for $\xi\in\G^*$. In particular, when $\G = {\mathbb R}$ and ${\mathfrak m} \colon {\mathbb R} \to {\mathbb C}$ is a Fourier multiplier on $\R$ supported on $[-1/2,1/2)$, we shall occasionally abuse notation by identifying
 ${\mathfrak m}$ with its periodization 
	$$
	{\mathfrak m}_{\rm per}(\xi):=\sum_{n \in {\mathbb Z}} {\mathfrak m}(\xi + n),
	$$
	and consequently view ${\mathfrak m}$ as a Fourier multiplier on ${\mathbb T}$. For  a detailed comparison of $T_{{\mathbb R}}[{\mathfrak m}]$ and $T_{{\mathbb Z}}[{\mathfrak m}]$, we refer to \cite[Section 2]{MSW}.

For every $k\in\Z_+$, we denote by $(\G^*)^k$  the  direct product of $k$  copies of the group $\G^*$, and  denote by $m$  a smooth tempered function from $(\G^*)^k$ to $\mathbb{C}$.  We define the multilinear Fourier multiplier operator $B_\G[m]:(\Schwartz(\G))^k\to \Schwartz(\G)$
by the formula
    \beq\label{multioperator}
    B_\G[m](f_1,\ldots,f_k)(x):= \int_{(\G^*)^k} m(\xi)
    \big(\prod_{i\in [k]}\F_{\G} f_i(\xi_i)\big) e(-x \cdot (\xi_1+\cdots+\xi_k))\ d\nu_{(\G^*)^k}(\xi),
    \eeq
    where   $x\in \G$,   $\xi=(\xi_1,\ldots,\xi_k)$ and $\nu_{(\G^*)^k}(\xi)=\prod_{i\in [k]}\nu_{\G^*}(\xi_i)$. 
    We refer to $m$ as the \emph{symbol} of $B_\G[m]$. 
    \subsection{Terminology from the classical circle method}\label{subbb:Terminor}
We  fix some notation and  terminology from the classical circle method. 
Define the \emph{height} $\Height (\alpha) \in 2^\N$ of an arithmetic frequency $\alpha = \frac{b}{q} \mod 1 \in \Q/\Z$ by the formula
$$
\Height \big( \frac{b}{q} ~{\rm mod} ~1 \big):= \inf \{ 2^l:  2^l \ge q, l \in \N \}  \sim q
$$
whenever $q \in \Z_+$ and $b \in [q]^\times$.  For any $l \in \N$, $\mathrm k \in \Z$, we   define the arithmetic frequency sets
\beq\label{xinn0001}
(\Q/\Z)_{\leq l}:= \{\alpha \in \Q/\Z: \Height(\alpha) \leq 2^l\}, \qquad (\Q/\Z)_{ l}:=(\Q/\Z)_{\leq l}\setminus (\Q/\Z)_{\leq l-1},
\eeq
(with the convention that $(\Q/\Z)_{\leq -1}=\emptyset$) 
and the continuous frequency sets
$$ \R_{\leq \mathrm k }:= [-2^{\mathrm k},2^{\mathrm k}].$$
We define the major arcs
\beq\label{majorarcs22}
{\mathfrak M}_{\leq l, \leq \mathrm k}:=\pi( \R_{\leq \mathrm k} \times (\Q/\Z)_{\leq l} )
\eeq
where the addition homomorphism $\pi: \R \times \Q/\Z \to \T$ is given by
 $\pi(\theta,\alpha):= \alpha + \theta$; see
(\ref{pimap}) in Appendix \ref{appenmap}.
Thus, ${\mathfrak M}_{\leq l, \leq \mathrm k}$ consists of all elements of $\T$ of the form $\frac{b}{q} + \theta ~{\rm mod} ~ 1$ for some $q \in [2^l]$, $b \in [q]^\times$, and $\theta \in [-2^{\mathrm k},2^{\mathrm k}]$.  
These major arcs are the natural focus for our Fourier-analytic manipulations. Note that the height in this paper is the naive height from \cite{KMT22}. Since we will apply the  Ionescu–Wainger multiplier theorem for canonical fractions (see Section \ref{section:IW}), there is no need to introduce the more complicated height from \cite{KMT22}.

Let $(l,\mathrm k)\in \N\times \Z$, and recall  definitions  \eqref{eq:106},   (\ref{eta-resc}). We introduce the Ionescu--Wainger projections by setting
\begin{align}
\label{eq:132}
\Pi_{\le l, \le \mathrm k}f(x):= T_{\Z}^{\mathcal R_{\le 2^l}}[\eta_{\le \mathrm k}]f(x)
\quad \text{for} \quad x\in\Z,
\end{align}
where  $\mathcal R_{\le N}$ ($N\ge 1$) are
 1-periodic sets of the so-called canonical fractions 
defined by 
\begin{align}
\label{eq:129}
\mathcal R_{\le N} := \big \{ {b}/{q}\in\mathbb Q/\Z: q\in[N] \text{ and } (b, q)=1 \big \}.
\end{align}
Precisely, the elements of \eqref{eq:129} are equivalence classes $\frac{b}{q} + {\Z}$ in $\T = \R/\Z$. We often identify $\frac{b}{q} \in \Q$ with the equivalence class $\frac{b}{q} + \Z$.
\begin{remark}\label{rtwonotation}
One can easily see that 
the set  $\mathcal R_{\le 2^l}$ plays the same role as  the above $(\Q/\Z)_{\le l}$. In fact,
the notation \eqref{xinn0001} is employed to streamline the analysis of adelic groups and related abstract groups in Section \ref{section:large scale Adelic group}, whereas \eqref{eq:129} is used mainly for convenience in the notation of Fourier multiplier operators
like \eqref{eq:132}.
\end{remark}
\subsection{Ionescu--Wainger multiplier theorem for canonical fractions}
\label{section:IW}
Invoking   the  definitions of $T_{\G}^{\mathfrak{Q}}[\mathfrak m]$ from (\ref{eq:106}) and $T_{\G}[\mathfrak m]$ from \eqref{eq:100}, we introduce below the   Ionescu--Wainger multiplier theorem for canonical fractions; see  \cite[Theorem 3.3]{KMPWW24}.    
\begin{thm} \label{thmIW}
	Let $q \in [p_0',p_0]$ for some $p_0 \in 2\Z_+$. Then there exists a constant ${C}_{p_0} \in \R_+$ such that for  every   integer  $N\ge 100$ the following is true. Assume that 
	\[
	0 < \vartheta \leq (2 p_0 N^{p_0})^{-1}
	\] 
	and let  
	$\mathfrak m : \R \to L(H_1,H_2)$ be a measurable function supported on $[-\vartheta, \vartheta]$, whose values are bounded linear operators between two separable Hilbert spaces $H_{1}$ and $H_{2}$.
	Then
	\begin{align}
	\label{Iweq:376}
	\big\|T_{\Z}^{\mathcal R_{\le N}}[\mathfrak m]\big\|_{\ell^q(\Z;H_1)\to \ell^q(\Z;H_2)}
	\lesssim_{p_0} 2^{{\bf C}_{p_0}(N)}\big\|T_{\R}[\mathfrak m]\big\|_{L^{p_0}(\R;H_1)\to L^{p_0}(\R;H_2)},
	\end{align}
    where the power ${\bf C}_{p_0}(N)$ is given by   
    \beq\label{constantspecia}
    {\bf C}_{p_0}(N):={{ C}_{p_0} \log N \frac{\log \log \log N}{\log \log N}}.
    \eeq
\end{thm}
In \cite{KMPWW24}, a general dimensional Ionescu–Wainger multiplier theorem was obtained for canonical fractions. Here we focus on the one-dimensional result presented  in Theorem \ref{thmIW}, as it suffices for our purposes. In particular, we will  use (\ref{Iweq:376})  with  $H_1=H_2=\mathbb{C}$, and  
	$H_1=H_2=\ell^2$.  The lower bound for $N$ is imposed to ensure the logarithm in \eqref{constantspecia} remains well-defined, though it is not essential in our proof.
\subsection{Approximant to the von Mangoldt function}
\label{subsection:approx}
Let ${ \mu}$ be the M\"obius function, defined by $\mu(n)=(-1)^k$ if $n$ is the product of $k$ distinct primes and by $\mu(n)=0$ if $n$ is divisible by the square of a prime; see \cite[Chapter 1]{Iwaniec04}.
For $\om\ge 2$, we introduce the Cram\'er approximant $\La_{\Cramer,\om}$ given by 
\beq\label{defn:Crame}
\La_{\Cramer,\om}(n)=\frac{W}{\varphi(W)}{\ind {(n,W)=1}},\quad {\rm where} \quad W=\prod_{p\le \om}p,
\eeq
 and   the  Heath-Brown approximant $\La_{\HB,\om}$  given by 
\beq\label{HBE1}
\La_{\HB,\om}(n)=\sum_{q<\om}\frac{\mu(q)}{\varphi(q)}c_q(n),\quad {\rm where}\quad  c_q(n)=\sum_{r\in [q]^\times}e(rn/q).
\eeq
By convention, $\La_{\Cramer,\om}=\La_{\HB,\om}=1$ if $1\le \om<2$. 
For any $k\in\Z_+$ and $N\ge 1$,  the Cram\'er approximant $\La_{\Cramer,\om}$ and the von Mangoldt function $\La$ satisfy 
\begin{equation}\label{up1}
    0\le \La_{\Cramer,\om}\les \langle \Log\ \om\rangle \qquad {\rm and}\qquad
{\E}_{n\in [N]} \La(n)\les 1,
\end{equation}
while 
 the  Heath-Brown approximant $\Lambda_{\HB,\om}$  obeys 
the moment bounds 
\begin{equation}\label{1moment} 
    \E_{n \in [N]} |\Lambda_{\HB,\om}(n)|^k \lesssim_k \langle \Log \om\rangle^{2^k+k}.
\end{equation}
Indeed, 
 (\ref{up1}) is a direct result of $\prod_{p\le \om}p/(p-1)\les \langle \LL\  \om\rangle$ and the definition of $\La$; for the proof of  (\ref{1moment}), 
  we refer  to \cite[Lemma 4.6]{KMTT24}.

In the proof of our main result,  more precise  estimates related to  the ``little"  Gowers uniformity norm for  the above approximants are needed. We first introduce the definition of this norm.  
For $s\in \N$, we  define the little  Gowers norms $u^{s+1}[N]$  of a function $f: \mathbb{Z}\to \mathbb{C}$ by
\begin{align}\label{little}
\|f\|_{u^{s+1}[N]}:=\sup_{\deg P\leq s}\left|\mathbb{E}_{n\in [N]}f(n)e(P(n))\right|,
\end{align}
where $P$ ranges over all real-coefficient polynomials of degree at most $s$.
For more details on the little Gowers norms and the related   Gowers norms, we refer to \cite{Te24, Gowers01, GTZ12}. 
We now list some useful  estimates involving  the above approximants.
\begin{lemma} \label{stable}  
{\rm (i)}~ Let $N \geq 100$  and $1 \leq z \leq \exp({\Log^{1/10} }N)$. 
Then there exists a positive constant  $c$ independent of $N$ such that 
\beq\label{UP1}
\E_{n\in [N]} \La_{\Cramer,z}(n)=1+O(\exp(-c\Log^{4/5}N));
\eeq
moreover, if  $1 \leq q \leq z$, $b\in [q]$ and $I$ is an  interval  of length $N$,
\beq\label{approx;resides}
\E_{n \in I} \Lambda_{\Cramer,z}(n) \ind{n\equiv b\  ({\rm mod}\  q)} = \frac{\ind{(b,q)=1}}{\varphi(q)} + O(\exp(-c \LL^{4/5} N)).
\eeq
{\rm (ii)}~ Let $N \geq 100$,   $d \in \Z_+$ and $1 \leq z,w \leq \exp(\Log^{1/10} N)$. Then there exists a positive constant  $c_d$, depending only on $d$, such that 
\beq\label{stable:cramel}
\| \Lambda_{\Cramer, w} - \Lambda_{\Cramer, z} \|_{u^{d+1}(I)} \lesssim_{d} w^{-c_d} + z^{-c_d}
\eeq
for any interval $I$ of length $N$.\\
{\rm (iii)}~   Let $N \geq 1$, $d\in \Z_+ $ and $1 \leq w, Q, Q_1,Q_2 \leq \exp(\Log^{1/20} N)$.  Then 
there exists  a positive constant  $c_d'$, depending only on $d$,  such that 
\begin{align}
\| \Lambda_{\Cramer,w} - \Lambda_{\HB,Q} \|_{u^{d+1}(I)} \lesssim_d&\  w^{-c_d'} + Q^{-c_d'}\quad{\rm and}\label{stable:cramelHB}\\
\| \Lambda_{\HB,Q_1} - \Lambda_{\HB,Q_2} \|_{u^{d+1}(I)} \lesssim_d&\  Q_1^{-c_d'} + Q_2^{-c_d'}\label{stable:HBHBHB}
\end{align}
 for any interval $I$ of length $N$. 
\end{lemma}
For the proofs of \eqref{UP1}-\eqref{approx;resides}, \eqref{stable:cramel}, and \eqref{stable:cramelHB}-\eqref{stable:HBHBHB}, we refer to \cite{KMTT24}, 
Corollary 4.4, Lemma 4.5, and Proposition 4.7, respectively. In particular,  the fundamental lemma of sieve theory (see \cite{Iwaniec04}) and  some important arguments in \cite{TT25,MS21} were used 
in these proofs. 
Here,  the restriction 
$N\ge 100$ in both (i) and (ii) does not   affect our subsequent proof, as we will only focus on the case where $N\ge C$ with  sufficiently large $C>0$.  
\subsection{Generalized von Neumann  theorem by Ter\"av\"ainen}
Due to  the potential   presence of   Siegel zeros  (see \cite[Subsections 1.4 and 6.2]{KMTT24} for the details)  when analyzing  the symbol of  the averages weighted by the von Mangoldt function 
$\La$, we employ the following generalized von Neumann  theorem   to reduce the matter to  bounding  the  Cram\'{e}r-weighted averages.
\begin{lemma}
\label{VNT} 
Let $d,k\in \Z_+$, $C\ge 1$ and  let  $\mathcal{P}=\{P_1,\ldots, P_k\}$ be a family of polynomials with integer coefficients and distinct degrees, where $\max\{\deg P_i: i\in [k]\}=d$.  Let $N\geq 1$,  $f_0,\ldots,f_k: \mathbb{Z}\to \mathbb{C}$ be 1-bounded functions supported on $[\pm C N^d]$, and let $\theta:[N]\to \mathbb{C}$ be a 1-bounded function. Then there exists some  $1\leq K\les_{d}1$ such that 
\begin{align*}
\big|\frac{1}{N^{d+1}}\sum_{m\in \mathbb{Z}}\sum_{n\in [N]}\theta(n)f_0(m)f_1(m+P_1(n))\cdots f_k(m+P_k(n))\big|\les_{C,\mathcal{P}} (N^{-1}+\|\theta\|_{u^{d+1}[N]})^{1/K},   
\end{align*}
where $u^{d+1}[N]$ is  defined by (\ref{little}) with $s=d$.
\end{lemma} 
Ter\"av\"ainen \cite{Te24} proved this lemma using Peluse's inverse theorem \cite[Theorem 3.3]{P2} and some key reduction arguments. Additionally, the condition $|\theta|\le 1$  
 explains why, in the proof of Proposition \ref{Prop:Reduction linear} below, we cannot directly reduce the matter to estimating the averages weighted by  the  function $\La_{\HB, \om}$ (which does not have a good upper bound). 
\section{Multilinear Weyl inequality with Cram\'{e}r weight}
\label{section:MWIII}
In this section, we  prove a multilinear Weyl inequality\footnote{We   use 
the terminology  ``multilinear Weyl inequality" throughout this paper
since this inequality  can be regarded as  a variant of the linear Weyl inequality.
Moreover, this terminology  has also been employed  in \cite{KMPWW24} to investigate 
the pointwise convergence of the general unweighted averages \eqref{commaveragesBerglesonopen}.} which will give us minor arc estimates for the Cram\'{e}r-weighted averages. More precisely, this inequality  
 asserts that the Cram\'{e}r-weighted averages
are
negligible when the  Fourier transform of $f_j$, for at least
one $j\in[k]$, vanishes on appropriate major arcs.   
Consider a polynomial mapping
\begin{align}
\label{eq:42}
{\mathcal P}:=(P_1,\ldots, P_k): \mathbb \Z\to \mathbb \Z^k,
\end{align}
where 
$P_1, \ldots ,P_k$ are polynomials with integer coefficients and distinct degrees such that
\begin{align}
\label{eq:41}
1\le d_1:= \deg P_1<\cdots < d_k:= \deg P_k.
\end{align}

Throughout this paper, we fix a  large constant  $C_0$ (say, $C_0=100$), 
and  consider the Cram\'{e}r approximant  (see \eqref{defn:Crame} for  the definition)
\beq\label{defn:choicespeciao}
\La_N(n):= \La_{\Cramer,\exp(\Log^{1/C_0}N)}(n).
\eeq
We  now present  the   multilinear Weyl inequality for the Cram\'{e}r-weighted averages, expressed by \eqref{moz1} and \eqref{moz2} with $\La_N$ in place of 
the von Mangoldt function  $\La$.
\begin{thm}\label{Thm:uweyl1}
  Let  $N \geq 1$, $l\in \N$, $k\in\Z_+$,    and let $\mathcal P$ be a polynomial mapping satisfying (\ref{eq:42}) and  (\ref{eq:41}).  Let
$1<q_1,\ldots, q_k<\infty$ be exponents such that
$\frac{1}{q_1}+\cdots+\frac{1}{q_k}=\frac{1}{q}\le 1$. 
There exist two positive constants  $c$ and  $C$,   depending on
$\mathcal P, q_1, \dots, q_k$, such that the following holds. Let 
$f_i \in \ell^{q_i}(\Z)$ for each  $i\in [k]$, and fix $j\in [k]$.
If $f_j \in \ell^{q_j}(\Z) \cap \ell^{2}(\Z)$,  and
 $\F_{\Z} f_j$ vanishes on the major arcs
${\mathfrak M}_{\le l,\le -d_j \Log N+d_j l}$ given by (\ref{majorarcs22}),
\begin{align}
\label{uweyl11}
\big\| \tilde A^{\mathcal{P}}_{N, \La_N}(f_1,\ldots, f_k) \big\|_{\ell^q(\Z)} 
\le\ C\big(2^{-cl}+\exp(-c\Log^{1/C_0}N)\big) \|f_1\|_{\ell^{q_1}(\Z)}\cdots \|f_k\|_{\ell^{q_k}(\Z)}.
\end{align}
The same conclusion holds for $A^{\mathcal{P}}_{N, \La_N}$ in place of
 $\tilde A^{\mathcal{P}}_{N, \La_N}$.
\end{thm}
The proof of Theorem \ref{Thm:uweyl1} will be carried out in several steps. We first prove an inverse theorem for weighted averages $\tilde A^{\mathcal{P}}_{N, \La_N}$ which we then use to give us structural information on the adjoints of these averages. 

\begin{remark}\label{Remark:weylmulti}
(1)  This bound $\exp(-c\Log^{1/C_0}N)$ comes from the lower bound in  \eqref{lcondition1} of 
Theorem \ref{inverse2}. Furthermore, 
 the upper bound $2^{-cl}+\exp(-c\Log^{1/C_0}N)$  does not meet the polynomial-type requirement $2^{-cl}+N^{-c}$ for the method in \cite{KMPWW24}. \\
(2)  It suffices to prove Theorem \ref{Thm:uweyl1} under this assumption    $2^l\les\exp(\Log^{1/C_0}N)$, since for larger values of $l$, the support condition  becomes stronger,  and the conclusion (\ref{uweyl11}) is essentially unchanged. Moreover, we may further assume that $l$ is sufficiently large, since otherwise the claim follows from Minkowski's inequality,  (\ref{UP1}) and H\"{o}lder's inequality. \\
(3) As described in Subsection \ref{subsectionbasicnotation}, we will absorb the constants of most inequalities (except for a few) into the notation $``\les"$, since the dependencies will be clear.
\end{remark}
\subsection{Inverse theorems}\label{subsection inverse}
We adopt the notation from  \cite{KMT22} and   
 use the inner product 
$\langle f,g\rangle:=\sum_{x\in\Z}f(x)g(x)$
on ${\Schwartz }(\Z)$.   (In the following multilinear analysis, inserting a complex conjugation into the inner product does not provide any  advantage.)  Let $k\in\Z_+$,  $h,f_1,\ldots,f_k\in {\Schwartz }(\Z)$,   $N\ge 1$,  and $j\in [k]$. Let  $w:\Z_+\to \mathbb {C}$ be  a weight function. 
Denote  $\tilde{A}_{N,w}^{\mathcal P}:=\tilde{A}_{N,w;\Z}^{\mathcal P}$ and ${A}_{N,w}^{\mathcal P}:={A}_{N,w;\Z}^{\mathcal P}$, and 
 observe the identity 
\begin{equation}\label{dual1}
 \big\langle \tilde{A}_{N,w}^{\mathcal P}(f_1,\ldots,f_k), h \big\rangle =\big \langle \tilde{A}_{N,w}^{\mathcal P,*j}(f_1,\ldots,f_{j-1},h,f_{j+1},\ldots,f_k), f_j\big\rangle
 \end{equation}
 (with obvious modifications when $j\in \{1,k\}$), 
where   
the dual operator $\tilde{A}_{N,w}^{{\mathcal P}*j}$ is given by 
\beq\label{dfen:dualj}
\tilde{A}_{N,w}^{\mathcal P,*j}(g_1,\ldots,g_k)(x):=\frac{1}{\lfloor N\rfloor}\sum_{n\in J_N} w(n)
\prod_{i\in [k]}  g_i\big(x-{\ind {i\neq j}}P_i(n)+P_j(n)\big) 
\eeq
with $J_N=[N]\setminus[N/2]$. 
A similar identity also holds for ${A}_{N,w}^{\mathcal P}$ and its dual operator 
${A}_{N,w}^{{\mathcal P},*j}$ for each $j\in [k]$.
Additionally, we shall often abbreviate 
$$\mathfrak{A}_{N,w}:=\mathfrak{A}_{N,w}^{\mathcal P},\quad \mathfrak{A}_{N,w}^{*j}:=\mathfrak{A}_{N,w}^{\mathcal P,*j} \quad {\rm whenever}\quad  \mathfrak{A}\in \{A,\tilde A\}.$$ 
We will use  the notation $O(1)$ to denote a constant
independent of $l$, $\delta$ and $N$, which  may change from line to line.
\subsubsection{Unweighted inverse theorem}
\label{subsubsection unweightedinverse}
We first state an important inverse theorem for the unweighted averages. Similar theorems can be found in \cite{KMT22,KMPWW24,KMPWfield}. 
\begin{thm}\label{inverse1}
 Let  $N \geq 1$,  $0<\delta\le 1$,  $k\in\Z_+$,    and let $\mathcal P$ be a polynomial mapping satisfying (\ref{eq:42}) and  (\ref{eq:41}).
Let $N_*$ be a quantity with $N_* \sim N^{d_k}$, and fix $j\in [k]$. Let $h,f_1,\ldots,f_k\in {\Schwartz}(\Z)$  be 
1-bounded functions supported on $[\pm N_*]$, obeying  the lower bound 
\begin{equation}\label{assL}
\big|\big\langle \tilde{A}_{N,1}(f_1,\ldots, f_k), h \big\rangle\big| \geq \delta N^{d_k}.
\end{equation}
Then there exists a function $F_j\in \ell^2(\Z)$ satisfying 
\beq\label{aa21}
\|F_j\|_{\ell^\infty(\Z)}\les 1,\quad \|F_j\|_{\ell^1(\Z)}\les N^{d_k},
\eeq
 and the property that  $\F_\Z F_j$ is  supported in the $O(\delta^{-O(1)}/N^{d_j})$-neighborhood of some  $\A_j\in \Q/\Z$ of height $O(\delta^{-O(1)})$, such that 
$$|\langle f_j,F_j\rangle|\gtrsim \delta^{O(1)} N^{d_k}.$$
The same conclusion holds for the operator ${A}_{N,1}$ in place of  $\tilde{A}_{N,1}$ in (\ref{assL}).
\end{thm}

\begin{proof}[Proof of Theorem \ref{inverse1}]
Peluse's inverse theorem (see \cite[Theorem 3.3]{P2}) will play a key
role in our proof. Note that we can simply take $F_j=\mathfrak{A}_{N,1}^{*j}(f_1,\ldots,f_{j-1},h,f_{j+1},\ldots,f_k)$ with $\mathfrak{A}=A$ or $\tilde A$  whenever 
$N\les \delta^{-O(1)}$. Thus,   it suffices 
to consider   
\beq\label{sufflarge}
N\ge \mathcal{C} \delta^{-\mathcal{C}}
\eeq
with 
sufficiently large constant $\mathcal{C}$.
Following  
 the arguments yielding  \cite[Theorem 6.3 and Proposition 6.6]{KMT22},  
we can infer that 
\beq\label{j=1case1}
{\rm Theorem \ \ref{inverse1}\   for\  the\  case}\  j=1\  {\rm holds}.
\eeq
We next  prove Theorem \ref{inverse1} for the case   $j=2$. 
It suffices to prove the desired result for $\tilde{A}_{N,1}$,  as the case for  ${A}_{N,1}$
 is analogous. The remainder of the proof closely follows \cite[Theorem 4.30]{KMPWW24}, so we omit the details.
\end{proof}

\subsubsection{Weighted inverse theorem}
\label{winver2}
This subsection presents the weighted extension of Theorem \ref{inverse1}. An important distinction from the unweighted case is the necessity to impose a lower bound on the parameter $\delta$   in the weighted setting.
\begin{thm}\label{inverse2} 
 Let  $N, k,\mathcal{P}$ and $N_*$ be given as in Theorem \ref{inverse1}.  Fix $j\in [k]$ and assume 
 \beq\label{lcondition1}
\exp(-c'\Log^{1/C_0} N) \le \delta \leq 1
\eeq
for some  sufficiently small  constant $c'=c'(k,d_k)>0$.\footnote{Since $k\le d_k$, we can also write $c'=c'(d_k)$.} Let $h,f_1,\ldots,f_k\in {\Schwartz}(\Z)$  be
1-bounded functions supported on $[\pm N_*]$, obeying  the lower bound 
\begin{equation}\label{ass1}
\big|\big\langle \tilde{A}_{N,\La_N}(f_1,\ldots, f_k), h \big\rangle\big| \geq \delta N^{d_k}.
\end{equation}
Then the conclusions of Theorem \ref{inverse1} hold.
 \end{thm}
    As in the proof of Theorem \ref{inverse1}, it suffices to prove Theorem \ref{inverse2} under the assumption that (\ref{sufflarge}) holds. We need  the following  proposition, which crucially depends on  \eqref{lcondition1}. 
 \begin{prop}\label{Prop:Reduction linear}
  Under  the hypotheses and notation of Theorem \ref{inverse2}, with the additional condition (\ref{sufflarge}),  there exists  a positive 
  $ \om_\circ=O(\delta^{-O(1)})$ such that  
  \beq\label{aim3.1} 
|\langle \tilde{A}_{N,\La_{\HB,\om_\circ}}(f_1,\ldots,f_k), h\rangle|\gtrsim \delta N^{d_k},
\eeq
 where  the weight $\La_{\HB,\om_\circ}$ is a Heath-Brown approximant  given by (\ref{HBE1})
with $\om=\om_\circ$.
 \end{prop}
 \begin{proof}[Proof of Theorem \ref{inverse2} (accepting Proposition \ref{Prop:Reduction linear})]
Using the definition of $\La_{\HB,\om_\circ}$ in (\ref{HBE1}),  we deduce from (\ref{aim3.1}) that  there exist positive integers   
 $q\le \lfloor \om_\circ\rfloor$ and $r\in [q]^\times$ such that 
 \beq\label{eqn:type1}
|\langle \tilde{A}_{N,e(r\cdot/q)}(f_1,\ldots,f_k), h\rangle|\gtrsim \delta^{O(1)} N^{d_k}.
\eeq
We begin by splitting the  proof into two cases, depending on whether the polynomial $P_1$ in $\mathcal {P}$
  is linear. We first prove the linear case, then deal with 
   the nonlinear case. 

   \smallskip \paragraph{\bf Case~1}
For the linear case, we assume   $P_1(n)=an+b$ for some  $ a,b\in\Z$ with $a\neq 0$. Write
$$e(\frac{rn}{q})=e\big(-\frac{r}{aq}(x-an-b)\big)\,e\big(\frac{r}{aq}(x-b)\big),\quad 
\tilde{f}_1:=e(-\frac{r}{aq}\cdot)\,f_1,\quad \tilde{h}:= e(\frac{r}{aq}(\cdot-b))\,h.
$$ 
Then (\ref{eqn:type1}) gives
$$|\langle \tilde{A}_N(\tilde{f}_1,f_2,\ldots,f_k), \tilde{h} \rangle|\gtrsim \delta^{O(1)} N^{d_k}.$$ 
Note that 
$\tilde{f}_1,\tilde h\in {\Schwartz}(\Z)$  are also
1-bounded functions supported on $[\pm N_*]$.
Hence, 
Theorem \ref{inverse1} gives us a function $G_j\in \ell^2(\Z)$ satisfying the desired $\ell^q$ bounds as in (\ref{aa21}),
and obeying that $\F_\Z G_j$ is supported in the $O(\delta^{-O(1)}/N^{d_j})$-neighborhood of some  $\A_j\in \Q/\Z$ of height $O(\delta^{-O(1)})$ such that 
$$|\langle {\ind {j=1}} \tilde{f}_1+{\ind {j\neq 1}}{f}_j,G_j\rangle|\gtrsim \delta^{O(1)} N^{d_k}.$$ 
From this, we see that
$$F_j:={\ind {j=1}}\, e(-\frac{r}{aq}\cdot)\, G_1+{\ind {j\neq1}}\,G_j$$
has the desired properties, concluding the proof of Theorem \ref{inverse2} for the case when $P_1$  is  linear.

\smallskip \paragraph{\bf Case~2}
We consider the  case   $\deg P_1\ge 2$. 
Here we use a new reduction argument.  
Using (\ref{eqn:type1}), Fubini's theorem,  and the change of variables $x\to x+n$, we have
\beq\label{eqn:type1 dge2}
  \big|\E_{n\in J_N} e(rn/q)~\big(\sum_{x\in\Z} f_1(x-\tilde{P}_1(n))\cdots f_k(x-\tilde{P}_k(n))~h(x+n)\big) \big|
  \gtrsim \delta^{O(1)} N^{d_k}
\eeq
with 
$\tilde{P}_i(n):={P}_i(n)-n$ for all $i\in [k]$.
By (\ref{eq:41}) and 
$\deg P_1\ge 2$,  we have  $\deg \tilde{P}_i=\deg P_i$  for each $i\in [k]$. 
Denote 
$$H(x):=e(rx /q) h(x),\qquad x\in\Z.$$
Applying  (\ref{eqn:type1 dge2}), $e(rn/q)=e(r(x+n)/q)e(-rx/q)$ and   Fubini's Theorem, we deduce     
$$ \|\E_{n\in J_N} f_1(x-\tilde{P}_1(n))\cdots f_k(x-\tilde{P}_k(n))~H(x+n) \|_{\ell^1(x\in\Z)}
  \gtrsim \delta^{O(1)} N^{d_k},$$ 
  which, combined with    the dual arguments and the support condition for $h$, yields that  there exists a 
  1-bounded function 
  $h_0\in \ell^\infty(\Z)$ supported on $[\pm N_{**}]$ with $N_{**}\sim N^{d_k}$ such that 
$$|\langle \mathcal{A}_N^{\tilde {\mathcal P}}(H,f_1,\ldots,f_k), h_0\rangle |\gtrsim \delta^{O(1)} N^{d_k},$$
where  ${\tilde {\mathcal P}}(n)=(-n,\tilde{P}_1(n),\ldots,\tilde{P}_k(n))$ is a map from 
$\Z \to \Z^{k+1}$ and 
 $\mathcal{A}_N^{\tilde {\mathcal P}}$ is given by 
$$\mathcal{A}_N^{\tilde {\mathcal P}}(H,f_1,\ldots,f_k)(x):=
\E_{n\in J_N}  {H(x+n)} ~f_1(x-\tilde{P}_1(n))\cdots f_k(x-\tilde{P}_k(n)).$$
Note that ${\tilde {\mathcal P}}$ is a family of polynomials with distinct degrees as well as the  polynomial  with the lowest degree   in ${\tilde {\mathcal P}}$ is linear; and 
$H,f_1,\ldots,f_k,h_0$ are all 1-bounded functions supported on $[\pm N_{**}]$.  As a result, 
we can obtain the  conclusions in  Theorem \ref{inverse2} 
for  the  case with $\deg P_1\ge 2$ 
by applying Theorem \ref{inverse1} with $(k,N_*)$ replaced by $(k+1,N_{**})$. 
 \end{proof}
 We will now prove Proposition \ref{Prop:Reduction linear} by generalizing the method used in \cite{KMTT24} (which addresses the special bilinear case) to the more general multilinear case.  For convenience, we adopt 
 $$\La_{{\rm Cr},\om}:=\La_{ \Cramer,\om}\qquad {\rm whenever}\ \om\ge1.$$
\begin{proof}[Proof of Proposition \ref{Prop:Reduction linear}]
Without loss of generality, we assume that $\delta$ is sufficiently small.
By applying  Lemma \ref{VNT} and \eqref{stable:cramel} 
 in order,  
 for any $z_1\sim z_2 $ with $1\le z_1,z_2\le \exp(\Log^{1/10}N)$, 
\beq\label{In:c}
\begin{aligned}
  N^{-d_k}|\langle \tilde{A}_{N,\La_{{\rm Cr},z_1}-\La_{{\rm Cr},z_2}}(f_1,\ldots,f_k), h\rangle|\les&\  
\langle {\LL}\  z_1\rangle  \big(N^{-1}+\|\La_{{\rm Cr},z_1}-\La_{{\rm Cr},z_2}\|_{u^{d_k+1}}\big)^{1/K}\\ 
\les&\  z_1^{-c_\circ}  
\end{aligned}
\eeq
for some positive  $K=K(d_k)$ and $c_\circ=c_\circ(d_k)$, where we have used 
  $|\La_{{\rm Cr},z_1}-\La_{{\rm Cr},z_2}|\les \langle {\LL}\  z_1\rangle$ derived from  (\ref{up1}). 
Then, for any $1\le \om\le \exp(\Log^{1/C_0} N)$,  we can deduce from (\ref{In:c}) that 
\beq\label{3.2}
\begin{aligned}
|\langle \tilde{A}_{N,\La_{N}-\La_{{\rm Cr},\om}}(f_1,\ldots,f_k), h\rangle|
\le&\  \sum_{n=\Log \om}^{\infty}|\langle \tilde{A}_{N,\La_{{\rm Cr},2^{n+1}}-\La_{{\rm Cr},2^n}}(f_1,\ldots,f_k), h\rangle|\\
&\ +O(\om^{-c_\circ}N^{d_k}
)\\
\les&\  \big(\om^{-c_\circ}+\sum_{n=\Log  \om}^{\infty} 2^{-c_\circ n}\big) N^{d_k}\les 
\om^{-c_\circ}N^{d_k}.
\end{aligned}
\eeq 
Thus, to prove (\ref{aim3.1}), it suffices   to show that    there exists $c_1>0$ such that 
\beq\label{3.24343}
\begin{aligned}
|\langle \tilde{A}_{N,\La_{\HB,\om}-\La_{{\rm Cr},\om}}(f_1,\ldots,f_k), h\rangle|
\les&\  
\om^{-c_1}N^{d_k}
\end{aligned}
\eeq 
for all $1\le \om\le \exp(\Log^{1/C_0} N)$.  Indeed, by \eqref{lcondition1}, we see that 
$$\om=\om_\circ:=\delta^{-2/\min\{c_\circ,c_1\}}$$
is an admissible choice for $\omega$; and applying   
 the triangle inequality along with  (\ref{3.24343}) and (\ref{3.2}),  we can finish the proof of    (\ref{aim3.1})   
 using this $\om_\circ$. 
 
To prove (\ref{3.24343}), we 
use Lemma \ref{VNT} again.
 However,  since this weight $\La_{\HB,\om}$ 	
  lacks  a suitable  upper bound like $\La_{{\rm Cr},\om}$, we cannot apply Lemma \ref{VNT} directly. To overcome this,  we 
introduce a modified  weight
  $\La_{\HB,\om}^\e$ defined  by 
$$\La_{\HB,\om}^\e(n):=\La_{\HB,\om}(n)~{\ind {|\La_{\HB,\om}(n)|\le \om^{c_\circ\e}}}
$$
with  sufficiently small $\e$  fixed later, where  $1\le \om\le \exp(\Log^{1/C_0} N)$ and  $c_\circ$ is  given by (\ref{In:c}).
So 
\beq\label{up2}
|\La_{\HB,\om}^\e(n)|\le \om^{c_\circ \e},
\eeq
 and  $\La_{\HB,\om}^\e$  vanishes on  $\{n\in\Z:\ |\La_{\HB,\om}(n)|> \om^{c_\circ \e}\}$, which, with (\ref{1moment}) ($k=k_\e:=2+\lfloor\e^{-1}\rfloor$), gives 
$$
\begin{aligned}
\langle\log\om\rangle^{2^{k_\e}+k_\e}\gtrsim \E_{n\in[N]}|\La_{\HB,\om}(n)|^{k_\e}
\gtrsim &\  \E_{n\in[N]}|D_{\HB,\om}(n)|^{k_\e} {\ind {|\La_{\HB,\om}(n)|>\om^{c_\circ \e}}}\\
\gtrsim&\  \om^{c_\circ}~ \E_{n\in[N]}|D_{\HB,\om}(n)|,
\end{aligned}
$$
where $D_{\HB,\om}:=\La_{\HB,\om}-\La_{\HB,\om}^\e$;
 consequently, we have 
\beq\label{CD1}
\E_{n\in[N]}|D_{\HB,\om}(n)|
\les \om^{-c_\circ} \langle\Log\om\rangle^{O_\e(1)}.
\eeq
Then, since (\ref{CD1}) and $1\le \om\le \exp(\Log^{1/C_0} N)$,   we immediately obtain a crude  estimate 
\beq\label{ai3.5}
|\langle \tilde{A}_{N,D_{\HB,\om}}(f_1,\ldots,f_k), h\rangle|
\les \E_{n\in[N]}|D_{\HB,\om}(n)|\les_\e 
\om^{-c_\circ} N.
\eeq
By applying (\ref{CD1})  and (\ref{stable:cramelHB}),
there exists a positive constant  $c_\circ'=c_\circ'(d_k)$ such that  
\beq\label{ai3.6}
\begin{aligned}
\|\La_{{\rm Cr},\om}-\La_{\HB,\om}^\e\|_{u^{d_k+1}[N]}
\le&\ \|\La_{{\rm Cr},\om}-\La_{\HB,\om}\|_{u^{d_k+1}[N]}+
\E_{n\in[N]}|D_{\HB,\om}(n)|\\
\les&\  \om^{-c_\circ'}+
\om^{-c_\circ} \langle\Log\om\rangle^{O_\e(1)}.
\end{aligned}
\eeq
In addition, by the upper bounds in (\ref{up1}) and (\ref{up2}), we have  
\beq\label{qa1}
|\La_{{\rm Cr},\om}-\La_{\HB,\om}^\e|\les \om^{c_\circ\e}+\Log \om
\les_\e \om^{c_\circ \e}.
\eeq
Combining   Lemma \ref{VNT} with  
(\ref{ai3.6})-(\ref{qa1}), we  obtain 
\beq\label{ai3.8}
\begin{aligned}
N^{-d_k}|\langle \tilde{A}_{N,\La_{{\rm Cr},\om}-\La_{\HB,\om}^\e}(f_1,\ldots,f_k), h\rangle|
\les&\  \om^{c_\circ \e} (N^{-1}+\|\La_{{\rm Cr},\om}-\La_{\HB,\om}^\e\|_{u^{d_k+1}[N]} )^{1/K}\\
\les&\   \om^{c_\circ \e} (\om^{-c_\circ'/K}+\langle\Log\om\rangle^{O_\e(1)} \om^{-c_\circ/K}).
\end{aligned}
\eeq

Set $\e=\min\{1/(2K),c_\circ'/(2c_\circ K)\}$. By  
  (\ref{ai3.8}) and (\ref{ai3.5}), we have 
\beq\label{aa2}
|\langle \tilde{A}_{N,\La_{{\rm Cr},\om}-\La_{\HB,\om}}(f_1,\ldots,f_k), h\rangle|
\les \om^{-\min\{c_\circ,c_\circ'\}/(4K)}N^{d_k},
\eeq
which yields (\ref{3.24343}) with $c_1=\min\{c_\circ,c_\circ'\}/(4K)$, and completes the proof of Proposition \ref{Prop:Reduction linear}.
\end{proof}

\subsection{Structure of dual functions}
With Theorem \ref{inverse2} in hand, 
we can obtain the structure of the  dual function  $\tilde{A}_{N,\La_N}^{*j}(f_1,\ldots,  f_k)$ for each $j\in [k]$, as defined in  (\ref{dfen:dualj}). 
\begin{thm}\label{struct}
 Let  $N, k,\mathcal{P},N_*$ and $\delta$ be given as in Theorem \ref{inverse2}. 
  Fix $j \in [k]$. Then for all
$1$-bounded functions $f_1,\ldots,f_k\in \Schwartz(\Z)$
supported on $[ \pm N_*]$, we have
\begin{equation}\label{decomp}
\tilde{A}_{N,\La_N}^{*j}(f_1,\ldots,  f_k)
= \sum_{\A_j \in \Q/\Z:~\Height(\A_j)\les \delta^{-O(1)}}F_{\A_j} + E_1+E_2,
\end{equation} 
where each $F_{\A_j} \in \ell^2(\Z)$  has Fourier transform supported in the $O(\delta^{-O(1)}/N^{d_j})$-neighborhood of $\A_j\in \Q/\Z$ and satisfies  the estimates 
\beq\label{eqn:FA}
\|F_{\A_j}||_{\ell^\infty(\Z)}\les \delta^{-O(1)}\quad \text{and} \quad\|F_{\A_j}||_{\ell^1(\Z)}\les \delta^{-O(1)}N^{d_k},
\eeq
and the error terms $E_1\in \ell^1(\Z)$ and $E_2\in \ell^2(\Z)$ satisfy the following estimates:
\beq\label{eqn:E1E2}
\|E_1\|_{\ell^1(\Z)}\le \delta N^{d_k}\quad \text{and} \quad\|E_2\|_{\ell^2(\Z)}\le \delta.
\eeq
\end{thm}
\begin{proof}[Proof of Theorem \ref{struct}]
Combining Theorem \ref{inverse2} with  the  Hahn-Banach theorem, we can obtain the desired 
 result 
 by following  the argument in  \cite[Corollary 6.10]{KMT22}.  
\end{proof}
By using the decomposition (\ref{decomp}) for each $j\in [k]$, 
we obtain two useful  $\ell^2$ bounds  for $\tilde{A}_{N,\La_{N}}^{*j}$ in the following two propositions. We extend  the approach in \cite{KMT22} and employ the Ionescu-Wainger multiplier result in Theorem \ref{thmIW}.
\begin{prop}\label{Prop:Strcut1}
Fix  $k\in \Z_+$ and $j\in[k]$. Let $N, N_*\ge 1$ with $N_*\sim N^{d_k}$, and $l\in\N$ with 
\beq\label{lowerbound}
c'\Log^{1/C_0} N\ge l
\eeq
where $c'$ is given by (\ref{lcondition1}).
Then there exists a positive  $c_0$,  depending only on $c'$,  such that 
\beq\label{eqn:structure 1}
\|(1-\Pi_{\le l,\le d_j(N,l)})\tilde{A}_{N,\La_{N}}^{*j}(f_1,\ldots,f_k)\|
_{\ell^2(\Z)}\les 2^{-c_0l}N^{d_k/2} \|f_1\|_{\ell^\infty(\Z)}
\cdots \|f_k\|_{\ell^\infty(\Z)}
\eeq
with $d_j(N,l):=-d_j\Log N+d_jl$,
whenever $f_1,\ldots,f_k\in \Schwartz(\Z)$ are supported on $[-N_*,N_*]$.
\end{prop}
\begin{proof}[Proof of Proposition \ref{Prop:Strcut1}]
We  proceed as in the proof of \cite[Proposition 6.15]{KMT22}. 
\end{proof}
\begin{prop}\label{pp2}
Let $k,j,N,N_*,l$ and $d_j(N,l)$ be given as in Proposition \ref{Prop:Strcut1}. Let $\bar{\iota}:[k]\to \{1,k\}$ denote the map $\bar{\iota}(m)=k$ if $m\in [k-1]$ and $\bar{\iota}(k)=1$. 
Then there exists a positive $c_0$, depending only on $c'$ given in (\ref{lcondition1}),  such that  
\beq\label{eqn:structure 2}
\|(1-\Pi_{\le l,\le d_j(N,l)})\tilde{A}_{N,\La_{N}}^{*j}(f_1,\ldots,f_k)\|
_{\ell^2(\Z)}\les 2^{-c_0l} \|f_{\bar{\iota}(j)}\|_{\ell^2(\Z)} \prod_{i\in [k]\setminus \{\bar{\iota}(j)\}}
\|f_i\|_{\ell^\infty(\Z)}.
\eeq
\end{prop}
 \begin{remark}\label{remark:inftyto l2}
 If $k\ge 3$, the definition $\bar{\iota}(k)=1$  can be replaced by  $\bar{\iota}(k)=m$ with $m\in [k]\setminus\{1,k\}$.     
 \end{remark}
 The proof of Proposition \ref{pp2} requires the following weighted  $L^{\rm p}$ improving bounds established in \cite[Lemma 5.1]{KMTT24} which in turn relies
on the estimate \eqref{stable:cramel}.

\begin{lemma}\label{lp-improvvv}  
Let $N\ge 1$,  let $Q$  be a polynomial with integer coefficients with  
 $\deg Q=d\geq 1$. Then there exists some  $q_0\in (1,2)$ such that  
\beq\label{eqn:weightedimpro}
\begin{aligned}
 \big\|\E_{n\in [N]}\big(\La(n)+\La_N(n)\big)|H(\cdot-Q(n))|\big\|_{\ell^2(\Z)}
\les_Q&\  N^{d(1/2-1/q_0)}\|H\|_{\ell^{q_0}(\Z)}\quad {\rm and}\\
 \big\|\E_{n\in [N]}\big(\La(n)+\La_N(n)\big)|H(\cdot-Q(n))|\big\|_{\ell^{q_0'}(\Z)}
\les_Q&\  N^{d(1/q_0'-1/q_0)}\|H\|_{\ell^{q_0}(\Z)},
\end{aligned}
\eeq
where  $q_0'=q_0/(q_0-1)$.
\end{lemma}
\begin{proof}[Proof of Proposition \ref{pp2}]
The linear case $k=1$ can be handled using a standard argument
(Plancherel's identity and basic exponent sum estimates), so we will
only consider the case $k\ge 2$. By modifying the arguments in the proofs of  \cite[Corollaries 6.22 and 6.24]{KMT22} and combining 
the weighted  $L^{\rm p}$ improving bounds in Lemma \ref{lp-improvvv},
we can achieve the desired result. 
\end{proof}
\subsection{Proof of Theorem \ref{Thm:uweyl1}}
\label{subsection: minor}
As stated  in $Remark$ \ref{Remark:weylmulti} after Theorem \ref{Thm:uweyl1}, 
we
can assume that 
$l\le c'\Log^{1/C_0}N$ with $c'$ given by \eqref{lowerbound},   and $l, N$ are  sufficiently large. 
Then it suffices to show  
$$
\big\| \tilde A^{\mathcal{P}}_{N, \La_N}(f_1,\ldots, f_k) \big\|_{\ell^q(\Z)} 
\les 2^{-cl}  \|f_1\|_{\ell^{q_1}(\Z)}\cdots \|f_k\|_{\ell^{q_k}(\Z)}.
$$ 
This immediately follows by H\"{o}lder's inequality with estimate (\ref{UP1}),
interpolation, duality and inequality \eqref{eqn:structure 2}; see \cite[p.1059]{KMT22}.
\section{Proof of Theorem \ref{t3shifted}: reduction to the major arcs}
\label{section:mainreduction}
In this section, we will reduce the proof of Theorem \ref{t3shifted} to establishing the major arc estimates. 
 \subsection{Reduction of  Theorem \ref{t3shifted} by the multilinear Weyl inequality}
 \label{subsection:by multilinear Weyl inequality}
 We will use the minor arc estimates  in Theorem \ref{Thm:uweyl1} to 
 reduce the proof of Theorem \ref{t3shifted} to showing  the following. 

\begin{thm}\label{lemma:majorarcsest}
 Let $k\in\Z_+$ with $k\ge 2$,  let $\mathcal P$ be a polynomial mapping satisfying (\ref{eq:42}) and (\ref{eq:41}). Let  $1< q_1,\ldots, q_k<\infty$ such that ${1}/{q_1} + \cdots + {1}/{q_k}={1}/{q}\le 1$. Then,  for any $r>2$,    
  \beq\label{eqn:majorestimate1}
\begin{aligned}
\|\big(\tilde{A}_{N,\La_{N}}^{\mathcal{P}}(\Pi_{\le l_{(N)},\le -d_1(\Log N-l_{(N)})} &f_1,\ldots,\Pi_{\le l_{(N)},\le -d_k(\Log N-l_{(N)})}f_k)\big)_{N\in \D}\|_{\ell^q(\Z;\V^r)}\\
 \les_{{\mathcal P},\la,r,q_1,\ldots,q_k}& \  \|f_1\|_{\ell^{q_{1}}(\Z)}\cdots \|f_k\|_{\ell^{q_{k}}(\Z)}
 \end{aligned}
 \eeq
  for any   $f_j\in \ell^{q_{j}}(\Z)$ with  $j\in [k]$,
 where   the scale 
 $ l_{(N)}$  is defined by 
  \beq\label{defnlog}
  l_{(N)}:=\lfloor \tilde{c} \Log^{1/C_0} N \rfloor
  \eeq
 for some  sufficiently small $ \tilde{c}=\tilde{c}(k)>0$.
\end{thm}

Set $p_0\in 2\Z_+$ large enough such that the parameters $q,\tilde{c},q_{1},\ldots,q_{k}$ in Theorem \ref{lemma:majorarcsest} satisfy 
\beq\label{assumptionnew}
q,\tilde{c},q_{1},\ldots,q_{k} \in [p_0',p_0].
\eeq
By Minkowski's inequality, (\ref{varsum}), (\ref{up1}) and H\"{o}lder's inequality, it is easy to establish
(\ref{shift}) in Theorem \ref{t3shifted}
with the lacunary set  $\D$ replaced by the subset $\D\cap \{N\les 1\} $. Thus,   by \eqref{simple}, 
it suffices to prove  (\ref{shift})  with the set 
 $\D$ replaced  by
 \beq\label{mathc0}
\D^B:=\D\cap \{N\ge \mathcal{C}_0\},
\eeq
where  $\mathcal{C}_0$ (depending on $\mathcal{P}$, $p_0$ and $C_0$ given as in Section \ref{section:MWIII})
is sufficiently large. 
In fact, the above  
 $ l_{(N)}$ is selected  according to the decay factor on the right-hand side of  (\ref{uweyl11}). 
Below we prove (\ref{shift}) under the assumption that  Theorem  \ref{lemma:majorarcsest} holds. All but a few inequalities will have their constants absorbed into the notation $``\les"$, as the dependencies will be clear.
 \begin{proof}[Proof of Theorem \ref{t3shifted} (accepting Theorem  \ref{lemma:majorarcsest})]
Since  $\D$ is lacunary,  by the triangle inequality  and \eqref{varsum},   it suffices to show the following two inequalities: first, 
\beq\label{eqn:appro1}
\|\tilde{A}_{N,\La-\La_{N}}^{\mathcal{P}}(f_1,\ldots,f_k)\|_{\ell^q(\Z)}
\les \langle \Log N \rangle^{-2}\|f_1\|_{\ell^{q_{1}}(\Z)}\cdots \|f_k\|_{\ell^{q_{k}}(\Z)};
\eeq
and second, 
\beq\label{eqn:major1}
\|\big(\tilde{A}_{N,\La_{N}}^{\mathcal{P}}(f_1,\ldots,f_k)\big)_{N\in \D^B}\|_{\ell^q(\Z;\V^r)}
\les  \|f_1\|_{\ell^{q_{1}}(\Z)}\cdots \|f_k\|_{\ell^{q_{k}}(\Z)}.
\eeq
Here the parameters  $q,\{q_{i}\}_{i\in [k]},r$ are given as in  Theorem  \ref{lemma:majorarcsest}. We will first prove  (\ref{eqn:appro1}) which follows the arguments in \cite{KMTT24}, extending them to our more general setting.

By 
 \eqref{up1},  (\ref{UP1}),  \eqref{defn:choicespeciao} and H\"{o}lder's inequality,
we have 
$$
\|\tilde{A}_{N,\La-\La_{N}}^{\mathcal{P}}(f_1,\ldots,f_k)\|_{\ell^q(\Z)}
\le \sum_{\om\in \{\La,\La_{N}\}}\|\tilde{A}_{N,\om}^{\mathcal{P}}(|f_1|,\ldots,|f_k|)\|_{\ell^q(\Z)}
\les  \prod_{i\in [k]}\|f_i\|_{\ell^{q_i}(\Z)},
$$
which implies,  by  interpolation,  that the proof of  (\ref{eqn:appro1}) is reduced to showing inequality  
\beq\label{3bound2}
\|\tilde{A}_{N,\La-\La_{N}}^{\mathcal{P}}(f_1,\ldots,f_k)\|_{\ell^1(\Z)}
\les_M \langle\Log N \rangle^{-M} \|f_1\|_{\ell^{2}(\Z)}\|f_k\|_{\ell^{2}(\Z)}
\prod_{i\in [k]\setminus\{1,k\}}\|f_i\|_{\ell^{\infty}(\Z)}
\eeq
 for any  $M\ge 1$.
To prove (\ref{3bound2}),   we observe that by duality, it suffices to show 
that for any $h,f_i\  (i\in [k]\setminus\{1,k\})$ with $\|h\|_{\ell^\infty(\Z)}=\|f_i\|_{\ell^\infty(\Z)}=1$, 
\beq\label{2bound}
\|\tilde{A}_{N,\La-\La_{N}}^{\mathcal{P},*1}(h,f_2,\ldots,f_k)\|_{\ell^2(\Z)}
\les_M \langle\Log N \rangle^{-M} \|f_k\|_{\ell^2(\Z)}.
\eeq
 Moreover,
 by the local nature of
 $\tilde{A}_{N,\La-\La_{N}}^{\mathcal{P},*1}$, it suffices  to prove \eqref{2bound} for 
the functions  $h,f_2,\ldots,f_k$  supported in $[-N_*,N_*]$ with $N_*\sim N^{d_k}$. 
From  \cite[Theorem 1.1]{MSTT23}\footnote{While \cite{MSTT23} only  treats  the case 
$C_0=10$, their methods extend to arbitrarily large $C_0$.} we conclude   that for any large $M_\circ\ge 1$, 
\beq\label{eqn:La-LaN}
\|\La-\La_{N}\|_{u^{d_k+1}[N]}\les_{M_\circ,C_0} \langle \Log N \rangle^{-M_\circ},
\eeq
where $\|\cdot\|_{u^{d_k+1}[N]}$ is defined  by \eqref{little}  with $s=d_{k}$. 
Again by duality, and using Lemma \ref{VNT} and (\ref{eqn:La-LaN}),   we see that for any $M_\circ\ge 1$,
\beq\label{2bound1}
\begin{aligned}
N^{-d_k }\|\tilde{A}_{N,\La-\La_{N}}^{\mathcal{P},*1}(h,f_2,\ldots,f_k)\|_{\ell^1(\Z)}
&\le\ \ {\rm Log} N \, (N^{-1}+\|\La-\La_{N}\|_{u^{d_k+1}[N]})^{1/K}  \|f_k\|_{\ell^\infty(\Z)}\\
&\les_{M_\circ}\langle \Log N \rangle^{-M_\circ} \|f_k\|_{\ell^\infty(\Z)}.
\end{aligned}
\eeq
By  
(\ref{eqn:weightedimpro})$_2$ with $(H,d,Q)=(f_k,d_k,P_k-P_1)$, there exists    $\tilde{q}\in (1,2)$ 
such that 
\beq\label{2bound2}
\begin{aligned}
\|\tilde{A}_{N,\La-\La_{N}}^{\mathcal{P},*1}(h,f_2,\ldots,f_k)\|_{\ell^{\tilde{q}'}(\Z)}
\le&\    \|\E_{n\in [N]}(\La(n)+\La_{N}(n))|f_k(\cdot+P_1(n)-P_k(n))|\|_{\ell^{\tilde{q}'}(\Z)}\\
\les&\  N^{d_k (1/\tilde{q}'-1/\tilde{q})}  \|f_k\|_{\ell^{\tilde{q}}(\Z)},\qquad {\rm where} \ \  \tilde{q}'={\tilde{q}}/({\tilde{q}-1}).
\end{aligned}
\eeq
Hence, (\ref{2bound}) follows by 
interpolating (\ref{2bound1}) with (\ref{2bound2}). This  completes the proof of (\ref{eqn:appro1}).

Next, we prove (\ref{eqn:major1}). 
For convenience, we denote
\beq\label{eqn:setpi222}
\Pi_{\le m,\le n}^1:=\Pi_{\le m,\le n}\qquad {\rm and}\qquad \Pi_{\le m,\le n}^2:=1-\Pi_{\le m,\le n} 
\eeq
  for  any $(m,n)\in\Z^2$. 
  Then, the proof of  (\ref{eqn:major1}) reduces to proving 
\beq\label{eqn:1minor arcsss}
\|\big(\tilde{A}_{N,\La_{N}}^{\mathcal{P}}(\Pi_{\le l_{(N)},\le -d_1(\Log N-l_{(N)})}^{\ka_1} f_1,\ldots,\Pi_{\le l_{(N)},\le -d_k(\Log N-l_{(N)})}^{\ka_k}f_k)\big)_{N\in \D^B}\|_{\ell^q(\Z;\V^r)}\les 1
\eeq
for any $\|f_i\|_{\ell^{q_{i}}(\Z)}=1$ and  any $(\ka_1,\ldots,\ka_k)\in \{1,2\}^k$, where $q,\{q_{i}\}_{i\in [k]},r$ and $l_{(N)}$ are given as in Theorem  \ref{lemma:majorarcsest}. 
By  (\ref{eqn:majorestimate1}), the bound (\ref{eqn:1minor arcsss}) follows whenever $\ka_j=1$ for all $j\in [k]$. For the remaining cases,  
we apply  Theorem \ref{Thm:uweyl1}.  Since  $\D$ is lacunary,   it suffices to show 
that for every $(\ka_1,\ldots,\ka_k)\neq (1,\ldots,1)$,
\beq\label{eqn:1minor arcs}
\|\tilde{A}_{N,\La_{N}}^{\mathcal{P}}(\Pi_{\le l_{(N)},\le -d_1(\Log N-l_{(N)})}^{\ka_1} f_1,\ldots,\Pi_{\le l_{(N)},\le -d_k(\Log N-l_{(N)})}^{\ka_k}f_k)\|_{\ell^q(\Z)}
\les \langle \Log N\rangle^{-2}.
\eeq
 Assume $\ka_{j_0}=2$ for some $j_0\in [k]$.  
From  the definitions   (\ref{eqn:setpi222}) and \eqref{eq:132}, we observe  that the Fourier transform of  $ \Pi_{\le l_{(N)},\le -d_{j_0}(\Log N-l_{(N)})}^{2}f_{j_0}$ vanishes on $\mathfrak{M}_{\le l_{(N)}-1, \le -d_{j_0}(\Log N-l_{(N)}+1)}$.  Thus, since  $N\in \D^B$, 
by applying Theorem \ref{Thm:uweyl1} (with $j=j_0$),  Theorem \ref{thmIW}  and \eqref{defnlog}, we see that the left-hand side of (\ref{eqn:1minor arcs}) is at most 
$$
 \begin{aligned}
\les 2^{-cl_{(N)}} 2^{k {\bf C}_{p_0}(2^{l_{(N)}})}
\les 2^{-\frac{c}{2}l_{(N)}}\les \langle \Log N\rangle^{-2} ,
 \end{aligned}
 $$
 where  ${\bf C}_{p_0}(\cdot)$ is given in \eqref{constantspecia}. This
 shows (\ref{eqn:1minor arcs}), concluding 
 the proof of (\ref{eqn:major1}).
 \end{proof}
 \subsection{Reduction of Theorem  \ref{lemma:majorarcsest} via dyadic decompositions}
 It suffices to prove \eqref{eqn:majorestimate1} with the set $\D$ replaced by $\D^B$.
In this section, we will employ both  an 
``arithmetic" dyadic decomposition  and  a ``continuous" dyadic decomposition as in
\cite[p.1047]{KMT22}. 

Recalling the definition \eqref{eq:132}, 
we define 
$$\Pi_{m,\le n}=\Pi_{\le m,\le n}-\Pi_{\le m-1,\le n}$$
with the convention $\Pi_{\le -1,\le n}=0$. The decomposition
$$\Pi_{\le m,\le n}=\sum_{0\le m'\le m}\Pi_{m',\le n}$$
gives us the following
arithmetic decomposition of the averages in  (\ref{eqn:majorestimate1}): 
\beq\label{eqn.arithdecom}
\sum_{l_1\in \N_{\le l_{(N)}}}\dots\sum_{l_k\in \N_{\le l_{(N)}}}
\tilde{A}_{N,\La_{N}}^{\mathcal{P}}\left(\Pi_{l_1,\le -d_1(\Log N-l_{(N)})}f_1,\ldots,\Pi_{l_k,\le -d_k(\Log N-l_{(N)})}f_k\right).
\eeq
For each $i \in [k]$, we decompose the bump function $\eta_{\leq -d_i (\Log N - l_{(N)})}$ (defined in \eqref{eta-resc}) into non-oscillatory and highly oscillatory components as follows:
$$
\eta_{\leq -d_i (\Log N- l_{(N)})}=\eta_{N}^{0}+\sum_{s_i \in [l_{(N)}]} \eta_{N}^{ s_i}
$$
where 
\begin{align}
\label{eqn:support}
\eta_{N}^{ s_i}:=
\begin{cases}
\eta_{\leq -d_i(\Log N  - s_i)}-\eta_{\leq -d_i(\Log N  - s_i+1)} & \text{ if } s_i\in \Z_+,\\
\eta_{\le-d_i\Log N} & \text{ if } s_i=0.
\end{cases}
\end{align}
We denote
\begin{align}
\label{eq:31}
\tilde{\Pi}_{l_i,s_i}^{N}:=T_{\Z}^{\mathfrak{Q}_{l_i}}[\eta_{N}^{ s_i}]
\end{align}
where (recall the definition \eqref{eq:129})
 \beq\label{dande}
\mathfrak{Q}_{\le n}:=\mathcal{R}_{\le 2^n}\qquad {\rm and} \qquad\mathfrak{Q}_{ n}:=\mathcal{R}_{\le 2^n}\setminus \mathcal{R}_{\le 2^{n-1}},\quad n\in\N;
\eeq
 here we set   $\mathcal{R}_{\le 2^{-1}}=\emptyset$.

Combining  the support property (\ref{eqn:support}) and the definitions \eqref{eq:100}, \eqref{eq:106}, (\ref{eq:31}), we rewrite (\ref{eqn.arithdecom}) as
\begin{gather*}
\sum_{l_1, s_1\in\N_{\le l_{(N)}}} \dots \sum_{l_k, s_k\in\N_{\le l_{(N)}}}
\tilde A_{N,\La_{N}}^{\mathcal{P}}\left(\tilde{\Pi}_{l_1,s_1}^{N}f_1,\ldots, \tilde{\Pi}_{l_k,s_k}^{N}f_k\right).
\end{gather*}
To prove Theorem  \ref{lemma:majorarcsest}, it suffices to show that there exists  $c_0>0$ such that
\beq\label{t31}
 \|\big(\tilde{A}_{N,\La_{N}}^{\mathcal{P}}(\tilde{\Pi}_{l_1,s_1}^{N}f_1,\ldots,\tilde{\Pi}_{l_k,s_k}^{N}f_k)
\big)_{N\in\D_{l,s}}\|_{\ell^q(\Z;\V^r)}\\
\les\  2^{-c_0(l+s)}\|f_1\|_{L^{q_{1}}(\Z)}\cdots \|f_k\|_{L^{q_{k}}(\Z)}  
\eeq
 with $q,\{q_{i}\}_{i\in [k]},r$ given as in  Theorem  \ref{lemma:majorarcsest}, and
where  the set $\D_{l,s}$ is given by
\beq\label{setset34}
\begin{aligned}
\D_{l,s}=&\ \{N\in\D^B:\  l_{(N)}\ge \max\{l,s\}\}\qquad\ \  {\rm with}\\ 
l=&\ \max\{l_1,\ldots,l_k\}\qquad \qquad{\rm and}\qquad \quad s=\max\{s_1,\ldots,s_k\}.
\end{aligned}
\eeq
\subsection{Major arcs approximation and a reduction of   \eqref{t31}}
\label{subsect:majorarcsappro}
We begin by establishing some necessary preliminaries for the major arcs approximation. 
For any $N\ge 1$, define the exponential sum and its continuous counterpart  by
\beq\label{exponentialsum}
\begin{aligned}
    m_{N,\La_N}(\xi):=&\ \frac{1}{\lfloor N\rfloor}\sum_{n\in J_N}\La_N(n)~ e(\xi\cdot \mathcal{P}(n))\qquad (\xi\in \T^k)
      \quad {\rm and}\\
    \tilde{m}_{N,\R}(\zeta):=&\ \int_{1/2}^1e(\zeta\cdot \mathcal{P}(Nt))dt\qquad  \ \ \quad\qquad\ (\zeta\in \R^k),
\end{aligned}
\eeq
where $J_N=[N]\setminus[N/2]$. 
Let $\mathfrak{Q}_{l_1,\ldots,l_k}$ denote   the  set of rational fractions 
 \beq\label{eqn:setset1}
 \mathfrak{Q}_{l_1,\ldots,l_k}:=\{{a}/{q}\in (\Q/\Z)^k:\ (a,q)=1,\ \ {a}/{q}\in \mathfrak{Q}_{l_1}\times\cdots\times \mathfrak{Q}_{l_k}\}
 \eeq
 where each $\mathfrak{Q}_{l_i}$ is defined as  (\ref{dande}). For $a/q\in \mathfrak{Q}_{l_1,\ldots,l_k}$, define 
  the exponential sum  
 \beq\label{eqn:exponentialsum}
 G^\times(\frac{a}{q})=\E_{n\in [q]^\times} e(\frac{a}{q} \cdot \mathcal{P}(n))=\frac{1}{\varphi(q)}\sum_{n\in [q]^\times}e(\frac{a}{q} \cdot \mathcal{P}(n)).
 \eeq
 Since $l=\max\{l_1,\ldots,l_k\}$, we can see   
 \beq\label{eqn:kgeqset}
 2^l\le q\le 2^{kl} 
 \eeq
 whenever $a/q\in \mathfrak{Q}_{l_1,\ldots,l_k}$.  Importantly,  changing  $(a,q)\to (Ka,Kq)$ for any $K\in\Z_+$ does not  affect  the expression in  (\ref{eqn:exponentialsum}).
\begin{prop}\label{pn1}
Let $\mathcal P : \R\to\R^k$ be a polynomial mapping  whose components are polynomials with integer coefficients and different degrees.
Let $N\in \D^B$, $M_1,\ldots,M_k\in \R_+$ and
$l_1,\ldots, l_k\in\N$ with $l:=\max\{l_1,\ldots,l_k\}\le \tilde c\Log^{1/C_0} N$ with $\tilde c$ given by \eqref{defnlog}. 
For each  $\xi\in\T^k$, $\theta=(\theta_1,\ldots,\theta_k)\in \mathfrak{Q}_{l_1,\ldots,l_k}$ and  $|\xi_i-\theta_i|\le M_i^{-1}$ for all $i\in [k]$, we have
\beq\label{eqn:approx1}
\left|m_{N,\La_N}(\xi)-G^\times(\theta)\tilde{m}_{N,\R}(\xi-\theta)\right|
\les_{\mathcal{P}}  2^{kl} \exp(-c \Log^{4/5} N)  \big( 1+\max\{N^{d_i}/M_i:{i\in [k]}\}\big)
\eeq
with $c$ given as in \eqref{approx;resides}.
\end{prop}
\begin{proof}[Proof of Proposition \ref{pn1}]
Denote $\theta=a/q$.  For $\tilde c$ sufficiently small, 
    $q$ obeys (\ref{eqn:kgeqset}) and 
 $2^l\le q\le \exp(\Log^{1/C_0}N)$. 
First we split the sum $m_{N,\La_N}(\xi)$ into residue classes mod $q$ so that
$$
m_{N,\La_N}(\xi) = \sum_{b\in [q]} \frac{1}{\lfloor N\rfloor}\sum_{n\in J_N}
e \big(\xi\cdot \mathcal{P}(n)\big) \La_N(n)  {\ind {n\equiv b~ ({\rm mod}~ q)}}.
$$
Noting
$e(\theta\cdot \mathcal{P}(n)) = e(\theta \cdot \mathcal{P}(b))$  when $n\equiv b$ (mod $q$), we see that
$$
\begin{aligned}
m_{N,\La_N}(\xi)=\  \sum_{b\in [q]} e (\theta\cdot \mathcal{P}(b)) 
E_{N,q,b}^{\mathcal P}(\xi-\theta),
\end{aligned}
$$
where 
\beq\label{defn:aver1}
E_{N,q,b}^{\mathcal P}(\zeta):=\frac{1}{\lfloor N\rfloor}\sum_{n\in J_N}
e \big(\zeta\cdot \mathcal{P}(n)\big) \La_N(n)  {\ind {n\equiv b~ ({\rm mod}~ q)}},\qquad \zeta\in\R^k.
\eeq
Since 
$\La_N(n)=0$ unless $(b,q)=1$,
we have
$$
\begin{aligned}
m_{N,\La_N}(\xi)=&  \sum_{b\in [q]^\times} e (\theta\cdot \mathcal{P}(b)) 
E_{N,q,b}^{\mathcal P}(\xi-\theta).
\end{aligned}
$$
To prove (\ref{eqn:approx1}), by applying  \eqref{eqn:kgeqset}, it suffices to show that 
\beq\label{eqn:approx2}
\begin{aligned}
&\ \ \big|
E_{N,q,b}^{\mathcal P}(\zeta)
-\frac{{\ind {(b,q)=1}}}{\varphi(q)}\tilde{m}_{N,\R}(\zeta)\big|
\les\   \exp(-c \Log^{4/5} N)  \big( 1+\max\{M_i^{-1}N^{d_i}:{i\in [k]}\}\big)
\end{aligned}
\eeq
for each $\zeta=(\zeta_1,\ldots,\zeta_k)\in\R^k$ with $|\zeta_i|\le M_i^{-1}$.  

 By    \eqref{approx;resides}, for every  $(n,N)\in J_N\times \D^B$, we have
\beq\label{er1}
 \mathrm{a}_n:=\frac{S_n}{n}-\frac{{\ind {(b,q)=1}}}{\varphi(q)}\qquad {\rm satisfying }\qquad |\mathrm{a}_n|\les \exp(-c\Log^{4/5} N),
\eeq
where 
$S_n:=\sum_{j\in [n]}\La_N(j){\ind {j\equiv b~ ({\rm mod}~ q)}}$. Then
 we  rewrite the average  $E_{N,q,b}^{\mathcal P}$ in \eqref{defn:aver1} as
\beq\label{eqn:averEE}
\begin{aligned}
E_{N,q,b}^{\mathcal P}(\zeta)
=&\  \frac{1}{\lfloor N\rfloor}\sum_{n\in J_N}e(\zeta\cdot \mathcal{P}(n)) \big(\frac{{\ind {(b,q)=1}}}{\varphi(q)}+b_n\big)
\end{aligned}
\eeq
with $b_n=n\mathrm{a}_n -(n-1)\mathrm{a}_{n-1}$.
By the mean value theorem, 
$$
\begin{aligned}
\big|\frac{1}{\lfloor N\rfloor}\sum_{n\in J_N} e(\zeta\cdot \mathcal{P}(n)) -\tilde{m}_{N,\R}(\zeta)\big|
\le &\ \frac{1}{\lfloor N\rfloor}\sum_{n\in J_N} \big|e(\zeta\cdot \mathcal{P}(n))- \int_{I_n} e(\zeta\cdot \mathcal{P}(t))dt\big|
+O(\frac{1}{N})\\
\les&\  N^{-1}\big(1+\max\{M_i^{-1}N^{d_i}:{i\in [k]}\}\big),
\end{aligned}
$$
where $I_n$ denotes the interval $[n,n+1]$. As a consequence, we obtain 
\beq\label{eqn:averEd1}
\begin{aligned} 
\frac{{\ind {(b,q)=1}}}{\varphi(q)} \big|\frac{1}{\lfloor N\rfloor}\sum_{n\in J_N} e(\zeta\cdot \mathcal{P}(n))-\tilde{m}_{N,\R}(\zeta)\big|
\les&\  N^{-1}\big(1+\max\{M_i^{-1}N^{d_i}:{i\in [k]}\}\big).
\end{aligned}
\eeq
On the other hand, we apply summation by parts and the mean value theorem, along with the upper bound of ${\mathrm a}_n$ (as given in (\ref{er1})),  to derive  
$$
\begin{aligned}
\big|\frac{1}{\lfloor N\rfloor}\sum_{n\in J_N}e(\zeta\cdot \mathcal{P}(n)) b_n\big| 
\les&\   \exp(-c \Log^{4/5} N)  \left( 1+\max\{M_i^{-1}N^{d_i}:{i\in [k]}\}\right).
\end{aligned}
$$
This,  with (\ref{eqn:averEd1}) and (\ref{eqn:averEE}), yields the desired  (\ref{eqn:approx2}). 
\end{proof}
The following proposition uses Proposition \ref{pn1} to allow us to reduce
\eqref{t31} to a major arcs approximation.
Given functions $m : \T^k\to \mathbb{C}$ and $S  : (\Q/\Z)^k\to \mathbb{C} $, we will be working with multilinear operators of the form
\begin{align}
\label{eq:multioperatorB}
B^{l_1,\ldots, l_k}[S; m](f_1,\ldots, f_k)(x):=
\sum_{\theta\in\mathfrak{Q}_{l_1,\ldots, l_k}}S(\theta)\sum_{y=(y_1,\ldots,y_k)\in\Z^k}K_{\tau_\theta m}(y)\prod_{i \in [k]} f_i(x-y_i)
\end{align}
for $x\in\Z$,
where $\tau_\theta m(\xi):= m(\xi-\theta)$, and 
\begin{align}
\label{eq:27kernel}
K_m(y):= \int_{\T^k}m(\xi)e(-\xi\cdot y)d\xi.
\end{align}

Now we formulate our approximation result. 
\begin{prop}
\label{Prop:approxmajor1}
Let   $N\ge \mathcal{C}_0$ with $\mathcal{C}_0$ given by \eqref{mathc0},  and let $l_1, s_1,\ldots, l_k, s_k\in\N$ obeying  $l_{(N)}\ge \max\{l, s\}$. Then, for all $f_1\in \ell^{q_1}(\Z),\ldots, f_k\in \ell^{q_k}(\Z)$ with $q_1,\ldots, q_k\in (1,\infty)$ such that $\frac{1}{q_1} + \dots + \frac{1}{q_k}=\frac{1}{q}\le 1$, we have
\beq
\label{eq:approxmajor1}
\begin{aligned}
   \|\tilde A_{N,\La_N}\big(\tilde{\Pi}_{l_1,s_1}^{N}f_1,\ldots, \tilde{\Pi}_{l_k,s_k}^{N}f_k\big)
-&B^{l_1,\ldots, l_k}[G^\times; \tilde m_{N,\R} w_N^{s_1,\ldots, s_k}](f_1,\ldots, f_k)
\|_{\ell^q(\Z)}\\
\lesssim&\  \exp(-\frac{c}{2}\Log^{4/5}N)
 \|f_1\|_{\ell^{q_1}(\Z)}\cdots  \|f_k\|_{\ell^{q_k}(\Z)}, 
\end{aligned}
\eeq
with $c$ given as in \eqref{approx;resides}, where 
 $w_N^{s_1,\ldots, s_k}(\xi):=\prod_{i \in [k]} \eta_{N}^{ s_i}(\xi_i)$ with  $\eta_{N}^{ s_i}$ given by  (\ref{eqn:support}). 
\end{prop}
\begin{proof}[Proof of Proposition \ref{Prop:approxmajor1}]
With Proposition \ref{pn1} in hand, we can obtain the desired result using a process similar to that in the proof of \cite[Proposition 7.28]{KMPWW24}.
\end{proof}
By using Proposition \ref{Prop:approxmajor1}, we now reduce the proof of  (\ref{t31}) to showing 
\beq\label{t41}
  \|\big(B^{l_1,\ldots, l_k}[G^\times; \tilde m_{N,\R} w_N^{s_1,\ldots, s_k}](f_1,\ldots, f_k)\big)_{N\in\D_{l,s}}\|_{\ell^q(\Z;\V^r)}
\les\  2^{-c_0(l+s)}\prod_{i\in [k]}\|f_i\|_{\ell^{q_{i}}(\Z)}
\eeq
 with $q,\{q_{i}\}_{i\in [k]},r$ given as in  Theorem  \ref{lemma:majorarcsest}.
This will then complete the proof of Theorem \ref{t3shifted}.
\subsubsection{Model operators}
For each $i\in [k]$,  we 
denote 
\begin{align}
\label{eqn51}
\eta_{N, t}^{s_i}(\zeta_i):=e(\zeta_i P_i(Nt))\eta_{N}^{s_i}(\zeta_i),\quad \zeta=(\zeta_1,\ldots,\zeta_k)\in \R^k,
\end{align}
where  $\eta_{N}^{s_i}$ is defined by \eqref{eqn:support}. Thus  
 we can express   $\tilde m_{N,\R} w_N^{s_1,\ldots, s_k}$ as
\begin{align*}
\tilde m_{N,\R}(\zeta) w_N^{s_1,\ldots, s_k}(\zeta)=
\int_{1/2}^1 \big(\prod_{i\in[k]}\eta_{N, t}^{ s_i}(\zeta_i)\big)dt.
\end{align*}
We define a new quantity 
\begin{align}
\label{eq:111}
u:= 
\begin{cases}
100  k (s^{5/4}+1) & \text{ if } s \ge  C_*l,\\
100  k(l^{5/4}+1) & \text{ if } C_*l > s,
\end{cases}
\end{align}
where we fix $C_*$ later to be a sufficiently large constant.
\begin{remark}\label{whywedonotuseKMPWW24}
    The exponet  $5/4$,  which can be substituted with any value in  $(5/4,\infty)$,  is determined by the error in (\ref{approx;resides}) (or \eqref{eqn:error1} below).  This  is another reason  why we cannot employ the multiparamater norm interchanging inequality which  plays a crucial role in \cite{KMPWW24}.
\end{remark}
 We define the high-frequency case when $C_*l \le s$,
and the low-frequency case when $C_*l > s$. 
As we shall see later, the low-frequency case requires a more involved analysis.
In fact,  it is the low-frequency case that requires multi-frequency analysis techniques.
For the high-frequency case,  we will use the multilinear Weyl inequality in the continuous setting, as detailed in Appendix \ref{section:Appen1}.
To obtain the desired estimate for  the low-frequency case, we  divide the set $\D_{l,s}$
into two subsets $\I_{\le }$ and $\I_{> }$ given by 
\beq\label{notatII}
\I_{\le }:=\I_{l,s,\le }=\D_{l,s}\cap [1,2^{10 p_02^u}]\qquad {\rm and }\qquad \I_{> }:=\I_{l,s,> }=\D_{l,s}\cap (2^{10 p_02^u},\infty),
\eeq
and reduce matters to proving \eqref{t41} with the set $\D_{l,s}$  replaced by $\I_{\le }$ and $\I_{> }$, which  correspond to  the
low-frequency cases  at  
small  scale  and large 
 scale, respectively. We will use the multilinear Rademacher-Menshov inequality and the harmonic analysis of the adelic integers $\Ad_\Z$  to address the small-scale case and the large-scale case, respectively.

Using the notation in \eqref{eq:106} and  setting 
 \beq\label{LLK1}
F_{N,t}^{l_i,s_i}:=F_{N,t}^{l_i,s_i}(f_i):=T^{\mathfrak{Q}_{l_i}}_\Z[\eta_{N,t}^{s_i}]f_i,\qquad i\in [k],
\eeq
 we may rewrite   ${B}^{l_1,\ldots,l_k}[G^\times(\theta);\tilde{m}_{N,\R} w_N^{s_1,\ldots, s_k}](f_1,\ldots,f_k)$ as
\beq\label{eqn:61}
\int_{1/2}^1 B^{l_1,\ldots, l_k}[G^\times; \eta_{u}^*](F_{N, t}^{ l_1, s_1},\ldots, F_{N, t}^{l_k, s_k}) dt,
\eeq
where the function $\eta_u^*$ is defined by 
\beq\label{cutoff}
\eta_u^*(\zeta)=\eta_{\le -10 p_0  u}(\zeta_1)\cdots  \eta_{\le -10 p_0  u}(\zeta_k),\quad \zeta=(\zeta_1,\ldots,\zeta_k)\in \R^k,
\eeq
and $p_0$ is given in \eqref{assumptionnew}, and $u$ is given by \eqref{eq:111}. In fact,
\eqref{defnlog}, \eqref{mathc0} and \eqref{eq:111} implies
$\eta_u^*(\zeta)=1$ on the support $\prod_{i\in [k]}\eta_N^{s_i}(\zeta_i)$ and so \eqref{eqn:61} follows.

It now suffices to consider  
\beq\label{eqn:operB}
B^{l_1,\ldots, l_k}[G^\times; \eta_{u}^*](F_{N, t}^{ l_1, s_1},\ldots, F_{N, t}^{l_k, s_k})
\eeq
for any $t\in [1/2,1]$. Importantly, passing to \eqref{eqn:operB} allows us to move back to the averages
$\tilde{A}_{2^u,\La_{2^u}}$ but now the varying parameter $N$ is fixed at $2^u$. 

Repeating the arguments yielding (\ref{eq:approxmajor1}), 
we can deduce that the difference between \eqref{eqn:operB} and the multilinear  operator
\beq\label{eqn:operA}
\tilde{A}_{2^u,\La_{2^u}}\big(\Pi_{l_1,\le -10 p_0u}(F_{N, t}^{ l_1, s_1}),\ldots, \Pi_{l_k,\le -10 p_0u}(F_{N, t}^{ l_k, s_k})\big)
\eeq
satisfies the error bound 
\beq\label{eqn:error1}
\begin{aligned}
   \|\big((\ref{eqn:operB})-(\ref{eqn:operA})\big)_{N\in \D_{l,s}}\|_{\ell^q(\Z;\V^r)}
\les  \exp(-\frac{c}{2}u^{4/5})\prod_{i\in[k]}\|(F_{N,t}^{l_i,s_i})_{N\in \D_{l,s}}\|_{\ell^{q_i}(\V^r)}
\end{aligned}
\eeq
with $1<q_1,\ldots, q_k<\infty$ satisfying ${1}/{q_1} + \cdots + {1}/{q_k} ={1}/{q}\le 1$. 

\begin{prop}\label{Prop:approxmajor1-again}
With the notation as above, for each $i\in[k]$, we have 
\begin{align}
\label{eq:60}
\begin{split}
\big\|\big(F_{N, t}^{ l_i, s_i}\big)_{ N\in\D_{l, s}}\big\|_{\ell^{q_i}(\Z;{\bf V}^r)}
& \lesssim (s+ 1) (2^{{\bf C}_{q_i}(2^{l})}{\ind {l\ge 10}} +1)
 \|f_i\|_{\ell^{q_i}(\Z)}
\end{split}
\end{align}
where
${\bf C}_{q_i}(\cdot)$ is given in  \eqref{constantspecia}.
\end{prop} 

\begin{proof}[Proof of Proposition \ref{Prop:approxmajor1-again}]
We divide the proof into two cases: $s_i>0$ and $s_i=0$. For the case $s_i>0$,  
we can use (\ref{Iweq:376}) and \cite[Theorem B.1]{KMT22}  to establish the stronger shifted square functions estimate  from which  (\ref{eq:60}) follows via \eqref{varsum}. For the case $s_i=0$, 
 we decompose  $\eta_{N,t}^{s_i}=\eta_{N}^{s_i}+\mathcal{E}_{N,t}^i$.  The proof that (\ref{eq:60}) holds with   $\eta_{N,t}^{s_i}$ replaced by $\eta_{N}^{s_i}$ is obtained
by combining  
\cite[Theorem 3.27]{KMPWW24} and  L{\'e}pingle's inequality \cite{Lep,MSZ1}. To prove (\ref{eq:60}) with $\eta_{N,t}^{s_i}$ replaced by $\mathcal{E}_{N,t}^i$, we use
 Theorem \ref{thmIW} and standard Littlewood-Paley arguments (see,
for example, \cite{DF} or \cite{MSZ2}). 
This completes the proof of (\ref{eq:60}).
\end{proof}  
Proposition \ref{Prop:approxmajor1-again}, together with \eqref{eqn:error1}, gives
\beq\label{eqn:error2}
\begin{aligned}
   \|\big((\ref{eqn:operB})-(\ref{eqn:operA})\big)_{N\in \D_{l,s}}\|_{\ell^q(\Z;\V^r)}
\les 2^{-c(l+s)}\prod_{i\in[k]}\|f_i\|_{\ell^{q_i}(\Z)}
\end{aligned}
\eeq
 with $q,\{q_{i}\}_{i\in [k]},r$ given as in  Theorem  \ref{lemma:majorarcsest}. By Minkowski's inequality,  \eqref{UP1}, \eqref{LLK1},  (\ref{eq:60}) and Theorem \ref{thmIW}, we have the following preliminary bound for \eqref{eqn:operA}: 
\beq\label{eqn:091}
\|\big((\ref{eqn:operA})\big)_{N\in \D_{l,s}}\|_{\ell^q(\Z;\V^r)}\les
(s+ 1)^k (2^{{\bf C}_{q_1}(2^{l})+\cdots+{\bf C}_{q_k}(2^{l})}{\ind {l\ge 10}} +1) \|f_1\|_{\ell^{q_1}(\Z)}\cdots \|f_k\|_{\ell^{q_k}(\Z)}.
\eeq
for $f_1\in\ell^{q_1}(\Z),\ldots,f_k\in\ell^{q_k}(\Z)$ with $1<q_1,\ldots, q_k<\infty$ satisfying ${1}/{q_1} + \cdots + {1}/{q_k} ={1}/{q}\le 1$. 

We end this subsection by establishing another useful variational inequality, which will be used to achieve the case $q=1$ for  the low-frequency case at a large scale.
\begin{lemma}\label{lemma:addww1}
We have the inequality \eqref{eqn:091} with  $(q,q_1,\ldots,q_k)$ replaced by some
  $$(q_*,q_{*,1},\ldots,q_{*,k})\in (0,1)\times(1,\infty)\times\cdots\times (1,\infty),$$
  where
   $\frac{1}{q_{*,1}}+\cdots+\frac{1}{q_{*,k}}=\frac{1}{q_*}$.
\end{lemma}
\begin{proof}[Proof of Lemma \ref{lemma:addww1}]
   Using (\ref{eq:60}) and   a routine computation,
$$\|\big(\tilde{A}_{2^u,\La_{2^u}}(G_N^1,\ldots,G_N^k)\big)_{N\in \D_{l,s}}\|_{\V^r}
\le \tilde{A}_{2^u,\La_{2^u}}\big(\|(G_N^1)_{N\in \D_{l,s}}\|_{\V^r},\ldots,\|(G_N^k)_{N\in \D_{l,s}}\|_{\V^r}\big),
$$
we  reduce matters to showing the inequality
\beq\label{eq:endgoal1222}
\|\tilde A_{2^u,\La_{2^u}}(H_1,\ldots,H_k)\|_{\ell^{q_*}(\Z)}
\les \|H_1\|_{\ell^{q_{*,1}}(\Z)}\cdots \|H_k\|_{\ell^{q_{*,k}}(\Z)}
\eeq
with the implicit constant independent of $u$. 
By the local nature of $\tilde A_{2^u,\La_{2^u}}$, it suffices to show  (\ref{eq:endgoal1222}) with all $H_i$ supported in 
$[\pm C_\circ 2^{d_k u}]$ ($C_\circ\sim_\mathcal P 1$). By \eqref{eqn:weightedimpro}, we have estimate 
 $\|\tilde A_{2^u,\La_{2^u}}^{*1}\|_{\ell^\infty\times\cdots\times\ell^\infty\times \ell^{q_0}\to\ell^{2}}
\les 2^{d_k u({1}/{2}-{1}/{q_0})}$ for some $q_0\in (1,2)$, which, by duality, gives
\beq\label{enddcay}
\|\tilde A_{2^u,\La_{2^u}}(H_1,\ldots,H_k)\|_{\ell^{1}(\Z)}\les 2^{d_k u(\frac{1}{2}-\frac{1}{q_0})}\|H_1\|_{\ell^2(\Z)}\|H_k\|_{\ell^{q_0}(\Z)} \prod_{i\in [k]\setminus\{1,k\}}\|H_i\|_{\ell^{\infty}(\Z)}.
\eeq
Set $v\in (2/3,1)$ such that ${1}/{v}={1}/{2}+{1}/{q_0}.$
By the support condition for $H_i$ with $i\in[k]$,  the function $\tilde A_{2^u,\La_{2^u}}(H_1,\ldots,H_k)$ is supported in an interval of length $O(2^{d_k u})$. Then, by
H\"{o}lder's inequality and (\ref{enddcay}),  we have 
$$
\begin{aligned}
    \|\tilde A_{2^u,\La_{2^u}}(H_1,\ldots,H_k)\|_{\ell^{v}(\Z)}\les&\  2^{(1-v)d_k u/v}
\|\tilde A_{2^u,\La_{2^u}}(H_1,\ldots,H_k)\|_{\ell^{1}(\Z)}\\
\les&\  \|H_1\|_{\ell^2(\Z)}\|H_k\|_{\ell^{q_0}(\Z)} \prod_{i\in [k]\setminus\{1,k\}}\|H_i\|_{\ell^{\infty}(\Z)},
\end{aligned}$$
which, together with the bound $\|\tilde A_{2^u,\La_{2^u}}\|_{\ell^{q_1}\times \cdots \times \ell^{q_k}\to \ell^1}\les 1$ for all $\frac{1}{q_1}+\cdots+\frac{1}{q_k}=1$ ($q_i\in [1,\infty]^k$), shows the desired 
(\ref{eq:endgoal1222}) by interpolation.
\end{proof}
\subsubsection{Further reductions of the proof of (\ref{t41})}
\label{subsubsectionAAAA}
From \eqref{eqn:091} and \eqref{eqn:error1}, it suffices  to consider
\beq\label{largezero}
\max\{s,C_* l\}\ge 10.
\eeq
We shall reduce the proof of (\ref{t41}) to showing  the following  propositions.
\begin{prop}[High-frequency case]\label{prop:highfre}
Let  $k\in\Z_+$ with $k\ge 2$, $s\ge C_*l$ with \eqref{largezero},  and let $2 \leq q_1,\ldots, q_k<\infty$ such that ${1}/{q_1} + \cdots + {1}/{q_k}=1$. Then there exists  $c\in(0, 1)$ independent of $C_*$ such that, for all $f_1\in \ell^{q_1}(\Z^k),\ldots, f_k\in \ell^{q_k}(\Z^k)$, we have
   \beq\label{eqn:large frequnce}
   \big\|\sup_{N\in\D_{l,s}}\big|{B}^{l_1,\cdots,l_k}[G^\times;\tilde{m}_{N,\R} w_N^{s_1,\ldots, s_k}](f_1,\cdots,f_k)
\big|\big\|_{\ell^1(\Z)}
\les 2^{-cs}\prod_{i\in [k]}\|f_i\|_{\ell^{q_i}(\Z)}.
\eeq
\end{prop}
 \begin{prop}[Low-frequency case at a  small scale]\label{t4.1}
Let $k\in\Z_+$ with $k\ge 2$,  $C_* l>s$ with \eqref{largezero},   and let $1< q_1,\ldots, q_k<\infty$ such that ${1}/{q_1} + \cdots + {1}/{q_k}={1}/{q}\le 1$. Then there exists  $c_1\in(0, 1)$ such that for
any $r\ge 2$ and for any $f_1\in \ell^{q_1}(\Z^k),\ldots, f_k\in \ell^{q_k}(\Z^k)$, we have
\beq\label{smallscale11}
 \begin{aligned}
 &\ \Big\|\big({B}^{l_1,\cdots,l_k}[G^\times;\eta_u^*]\big(F_{N,t}^{l_1,s_1},\cdots,F_{N,t}^{l_k,s_k}\big)\big)_{N\in \I_{\le}}
 \Big\|_{\ell^q(\Z;\V^r)}
 \les\ 2^{-c_1 l}\prod_{i\in [k]}\|f_i\|_{\ell^{q_i}(\Z)}
 \end{aligned}
 \eeq
 for all  $t\in [1/2,1]$, where $\{F_{N,t}^{l_i,s_i}\}_{i\in [k]}$ and $\eta_u^*$ are given by \eqref{LLK1} and \eqref{cutoff}.
\end{prop}
\begin{prop}[Low-frequency case at  a large scale] \label{t4.2}
Let $k\in\Z_+$ with $k\ge 2$,  $C_* l>s$ with \eqref{largezero}.   
Assume that $kD_k <r<\infty$ with  large $D_k$ (say, $D_k=10 kd_k$). Then there exists  $c_2\in(0, 1)$ such that   for any   $f_1\in \ell^{kD_k}(\Z^k),\ldots, f_k\in \ell^{kD_k}(\Z^k)$, 
\beq\label{eqn:adeilc analysis}
 \begin{aligned}
\Big\|\big({B}^{l_1,\cdots,l_k}[G^\times,\eta_u^*]\big(F_{N,t}^{l_1,s_1},\cdots,F_{N,t}^{l_k,s_k}\big)\big)_{N\in \I_{>}}
 \Big\|_{\ell^{D_k}(\Z;\V^r)}
 \les\ 2^{-c_2 l}\prod_{i\in [k]}\|f_i\|_{\ell^{kD_k}(\Z)},
 \end{aligned}
 \eeq
  for all  $t\in [1/2,1]$, where  $\{F_{N,t}^{l_i,s_i}\}_{i\in [k]}$ and $\eta_u^*$ are given by \eqref{LLK1} and \eqref{cutoff}.
\end{prop}
\begin{proof}[Proof of (\ref{t41}) (accepting Propositions \ref{prop:highfre}-\ref{t4.2})]
By (\ref{eqn:error1}),   to obtain (\ref{t41}), it suffices to show 
\beq\label{eqn:88812}
\|\big((\ref{eqn:operA})\big)_{N\in \D_{l,s}}\|_{\ell^q(\Z;\V^r)}
\les 2^{-c_0(l+s)}
\prod_{i\in [k]}\|f_i\|_{\ell^{q_i}(\Z)} 
\eeq
for  all  $q_1,\ldots,q_k\in [p_0',p_0]$ with $i\in [k]$ such that ${1}/{q_1} + \cdots + {1}/{q_k}= {1}/{q}\le 1$.

   In the  high-frequency case  $s\ge C_*l$,
   by  combining  \eqref{eqn:large frequnce} and  \eqref{eqn:error1}, and setting  $C_*$  sufficiently large, we deduce  (\ref{eqn:88812}) with $r=\infty$, $q=1$ and $q_i$ given as in Proposition \ref{prop:highfre}. Together with 
    \eqref{eqn:091}, this  yields  the desired result via  interpolation. 
Next, we treat the  low-frequency case  $s< C_*l$. Combining
   inequality \eqref{smallscale11} with
    (\ref{eqn:error1}) 
    implies (\ref{eqn:88812}) for the low-frequency case at  a small scale.
   On the other hand, by combining \eqref{eqn:adeilc analysis} and \eqref{eqn:error1}, we deduce  (\ref{eqn:88812}) with $kD_k <r<\infty$,  $q=D_k$ and $q_i=kD_k$ for all $i\in [k]$.  Interpolating this with   \eqref{eqn:091} and the estimate in Lemma \ref{lemma:addww1}  gives  
     (\ref{eqn:88812}) for the low-frequency case  at a large  scale. 
\end{proof}

We conclude this section by providing the proof of Proposition \ref{prop:highfre}.
\begin{proof}[Proof of Proposition \ref{prop:highfre}]
Normalize $\| f_i \|_{\ell^{q_i}(\Z)} = 1$ for all
$i \in [k]$ and note that it suffices to prove
\beq\label{eqn:goal1}
\sup_{\theta\in \mathfrak{
Q}_{l_1,\ldots, l_k}}\big\|\sup_{N\in\D_{l, s}}\big|T_{\theta,N}(f_1,\ldots, f_k)\big|\big\|_{\ell^{1}(\Z)} 
\lesssim 2^{-c' s}
\eeq
for some $c'>0$ independent of $C_*$,
where the operator $T_{\theta,N}$ is given by 
$$
T_{\theta,N}(f_1,\ldots, f_k)(x):= \sum_{y\in\Z^k}K_{\tau_\theta (\tilde{m}_{N,\R} w_N^{s_1,\ldots,s_k})}(y)\prod_{i \in [k]} f_i(x-y_i),\quad x\in\Z.
$$
Indeed, using Minkowski's inequality,  the bound  $\#\mathfrak{Q}_{l_1,\ldots, l_k}\le 2^{2kl}$, and \eqref{eq:multioperatorB}, the bound \eqref{eqn:large frequnce} follows 
by setting $C_*$ large enough such that $C_* c'>100 k$. 

It remains to establish \eqref{eqn:goal1}.
For each  $\theta\in \mathfrak{
Q}_{l_1,\ldots, l_k}$, we have
\begin{align*}
T_{\theta,N}(f_1,\ldots, f_k)(x)=e((\theta_1+\cdots+\theta_k)x)~T_{0,N}({\mathrm M}_{\theta_1}f_1,\ldots, {\mathrm M}_{\theta_k}f_k)(x), \quad x\in\Z,
\end{align*}
where ${\mathrm M}_{\theta_i}$ is the modulation operator defined  by ${\mathrm M}_{\theta_i}f(x):= e(-\theta_i  x)f(x)$
for $i\in [k]$. 
In view of this identity,  to prove (\ref{eqn:goal1}), 
it suffices  to establish 
\beq\label{eqn:goal2}
\big\|\sup_{N\in \mathbb{D}_{l, s}}\left|T_{0,N}(g_1,\ldots, g_k)\right|\big\|_{\ell^{1}(\Z)} 
\lesssim 2^{-c' s}
\eeq
for any $g_i\in \ell^{q_i}(\Z)$ such that $\|g_i\|_{\ell^{q_i}(\Z)}=1$ with $i\in[k]$.
By a transference argument (see  \cite[Lemma 4.4]{Bour89} or \cite[Theorem 7.42]{KMPWW24}), the proof of   (\ref{eqn:goal2}) can be reduced to showing
\beq\label{eqn:goal3}
\big\|\sup_{N\in \mathbb{D}_{l, s}}|\tilde T_{0,N}(g_1,\ldots, g_k)|\big\|_{L^{1}(\R)} 
\lesssim 2^{-c' s}
\eeq
for any $g_i\in L^{q_i}(\R)$ such that $\|g_i\|_{L^{q_i}(\R)}=1$ for $i\in[k]$, where $\tilde T_{0,N}(g_1,\ldots, g_k)$ is defined by 
$$\tilde T_{0,N}(g_1,\ldots, g_k)(x):= \int_{\R^k} K_{\tilde{m}_{N,\R} w_N^{s_1,\ldots,s_k}}(y)\, g_1(x-y_1)\,\cdots\, g_k(x-y_k) \,dy,\quad x\in\R.
$$

 We next  prove (\ref{eqn:goal3}). 
 Write $\tilde T_{0,N}(g_1,\ldots, g_k)$ as 
\beq\label{eqn:001}
\int_{1/2}^1 (\eta_N^{s_1}(D)g_1)(x-P_1(Nt))\,\cdots \,(\eta_N^{s_k}(D)g_k)(x-P_k(Nt))dt,
\eeq
where 
$(\eta_N^{s_i}(D)g_i)(x):=\int_\R \eta_N^{s_i}(\xi)\,\F_\R g_i(\xi)\, e(\xi x)\,d\xi$.
 Without loss of generality, we may assume that  $s$ is sufficiently large otherwise \eqref{eqn:goal3} can be proved by Minkowski's inequality, H\"{o}lder's inequality and Young's inequality.  In addition, we can further assume
 \beq\label{skbig}
s_1\le s_2\le\ldots\le s_{k-1}\le s_k=s;
\eeq
the other cases can be handled similarly.  Set $s_0=0$  and  
let $n\in \N_{\le k}$ be the biggest  integer such that $s_n=0$. 
Obviously, $n\le k-1$ (otherwise $s=0$). Without loss of generality, we only consider 
the case  $1\le n\le k-1$, since the case $n=0$  can be treated in the same way. 
 For each $i\in [n]$ (so we have  $s_i=0$), by Taylor's expansion $e(\xi P_i(t))=\sum_{m_i\in \N}(-2\pi {\rm i} \xi P_i(t))^{m_i}/(m_i!)$,   we 
 write 
$$
\begin{aligned}
    (\eta_N^{s_i}(D)g_i)(x-P_i(t))=&\ \int_\R \eta_N^{0}(\xi)\,\F_\R g_i(\xi)\, e(-\xi x)
\,e(\xi P_i(t)) d\xi\\
=&\ \sum_{m_i\in\N}\frac{(-2\pi {\rm i} N^{-d_i} P_i(t))^{m_i}}{m_i!}\,
(\eta_N^{0,m_i}(D)g_i)(x)
\end{aligned}$$
where  $\eta_N^{0,m_i}(D)g_i(x)=
\int_\R  (\xi N^{d_i})^{m_i} \eta_N^{0}(\xi)\F_\R g_i(\xi) e(\xi x)  d\xi$
can be seen as  a  variant of $\eta_N^{0}(D)g_i$. Rewrite (\ref{eqn:001}) as
\beq\label{eqn:002}
\sum_{m_1\in \N}\cdots \sum_{m_n\in \N}
\frac{(-2\pi {\rm i})^{m_1+\cdots+m_n}}{m_1!\cdots m_n!}
\Big(\prod_{i\in [n]}\eta_N^{0,m_i}(D)g_i(x)\Big)
\mathcal{A}_{N,\R}(g_{n+1},\ldots,g_k)(x),
\eeq
where 
the operator $\mathcal{A}_{N,\R}$ is given by 
$$\mathcal{A}_{N,\R}(g_{n+1},\ldots,g_k)(x):=\int_{1/2}^1\left(\prod_{i=n+1}^k
 \left(\eta_N^{s_i}(D)g_i\right)(x-P_i(Nt))
\right) \mathcal{B}_N(t)dt
$$
with $\mathcal{B}_N(t):=\prod_{i\in [n]}\big(\frac{P_i(Nt)}{N^{d_i}}\big)^{m_i}$  obeying 
$|\mathcal{B}_N|+|\mathcal{B}_N'|\les \tilde{C}^{m_{[1,n]}}$  for some $\tilde{C}=\tilde{C}_{\mathcal{P}}>0$,
where $m_{[1,n]}:=m_1+\dots+m_n$.
Define $\frac{1}{q_*}:=\frac{1}{q_{n+1}}+\cdots+\frac{1}{q_{k}}$ and 
 $\frac{1}{q_{**}}:=\frac{1}{q_{1}}+\cdots+\frac{1}{q_{n}}$. This yields  $\frac{1}{q_{*}}+\frac{1}{q_{**}}=1$ and $q_*,q_{**}\in (1,\infty)$. By a routine computation, there exists   $C>0$ such that 
 \beq\label{eqn:weiyou}
 \big\|\sup_{N\in \mathbb{D}_{l, s}}\big|\prod_{i\in [n]}\eta_N^{0,m_i}(D)g_i\big|\big\|_{L^{q_{**}}(\R)}
 \les C^{m_{[1,n]}}   \prod_{i\in [n]}\|M_{HL}g_i\|_{L^{q_{i}}(\R)}
 \les C^{m_{[1,n]}} \prod_{i\in [n]}\|g_i\|_{L^{q_{i}}(\R)}
 \eeq
 with the implicit constant independent of $m_1,\ldots,m_n$, 
 where  $M_{HL}$ represents the continuous  Hardy-Littlewood maximal operator.
  By (\ref{eqn:weiyou}) and (\ref{eqn:002}), the bound (\ref{eqn:goal3}) follow from
\beq\label{eqn:goal4500}
\big\|\big(\sum_{N\in \mathbb{D}_{l, s}}\big|\mathcal{A}_{N,\R}(g_{n+1},\ldots,g_k)\big|^{q_*}\big)^{1/q_*}\big\|_{L^{q_*}(\R)} 
\lesssim \tilde{C}^{m_{[1,n]}} 2^{-c' s}\prod_{i=n+1}^k\|g_i\|_{L^{q_i}(\R)}.
\eeq
Recall that $s=s_k>0$. We now apply the multilinear Weyl inequality in the continuous setting, 
Theorem \ref{appen11} with $j=k$ (see $Remark$ \ref{remaAppes} following Theorem \ref{appen11}), to bound the left-hand side of \eqref{eqn:goal4500} by  
$$
\begin{aligned}
    \big(\sum_{N\in \mathbb{D}_{l, s}}\|\mathcal{A}_{N,\R}(g_{n+1},\ldots,g_k)\|_{L^{q_*}(\R)}^{q_*}\big)^{1/q_*}
\les&\ \tilde{C}^{m_{[1,n]}} 2^{-c' s} \prod_{i=n+1}^k\|\big(\eta_N^{s_i}(D)g_i\big)_{N\in \mathbb{D}_{l, s}}\|_{L^{q_i}(\R;\ell^{q_i})}\\
\les&\ \tilde{C}^{m_{[1,n]}} 2^{-c' s} \prod_{i=n+1}^k\|\big(\eta_N^{s_i}(D)g_i\big)_{N\in \mathbb{D}_{l, s}}\|_{L^{q_i}(\R;\ell^{2})},
\end{aligned}$$
where the condition $q_{n+1},\ldots,q_k\ge 2$ is used in the second inequality. This yields \eqref{eqn:goal4500} by repeatedly
applying  standard Littlewood-Paley arguments (since $s_i>0$ for all $i\in [k]\setminus [n]$). This concludes the proof of Proposition \ref{prop:highfre}.
\end{proof}
 The proofs of Propositions \ref{t4.1} and \ref{t4.2} are given in the following two sections. Without loss of generality, we assume that $l$ is large.
\section{Major
arc estimates for the low-frequency case at a small scale}
\label{section:small scale}
   In this section, we prove Proposition \ref{t4.1} by establishing a multilinear Rademacher-Menshov inequality.   A natural approach to derive this inequality would be to directly  generalize the arguments used for the bilinear version  in \cite{KMT22}. However, the techniques used in \cite{KMT22} grow substantially more intricate as the parameter $k$ increases.
This can be circumvented by employing an induction argument,  the desired multilinear Rademacher-Menshov inequality can be then deduced directly from the bilinear version established in \cite{KMT22}.

\begin{lemma}[Rademacher-Menshov for multilinear forms]
\label{mil} 
 Let   $k\in\Z_+$ with $k\ge 2$,  and $K \in \Z_+$. For each $i\in [k]$,  let $\{f_{i,N}\}_{N\in [K]}$ be the elements of some vector space $V_i$.  Let $0< q < \infty$, and let $B : V_1 \times \cdots \times V_k \to L^q(X)$ be a multilinear  map where $X$ is some measure space.  Then
\begin{equation}\label{msmall}
\begin{split}
\big\| \big(B(f_{1,N},&\ldots,f_{k,N})\big)_{N \in [K]} \big\|_{L^q(X;\V^2)} 
\lesssim_q\  \langle \Log K \rangle^{k-2+k\max\{1,{1}/{q}\}}\\
&\ \ \times\sup \big\| B\big( \sum_{j \in [K]} \varepsilon_{j}^{(1)} (f_{1,j}-f_{1,j-1}),\ldots,\sum_{j\in [K]} \varepsilon_{j}^{(k)} (f_{k,j}-f_{k,j-1})\big) \big\|_{L^q(X)}
\end{split}
\end{equation}
with the convention $f_{i,0}=0$ for all $i\in [k]$, where the supremum is taken over  
$$\varepsilon_j^{(i)}\in  \{1,-1\}\qquad{\rm for\ all }\quad
(i,j)\in [k]\times [K].
$$
\end{lemma}
\begin{proof}[Proof of Lemma \ref{mil}]
We argue by induction on the parameter 
$k\ge 2$, the case $k=2$ 
 was established in \cite[Corollary 8.2]{KMT22}. 
For every $m\in \{3,4,5,\ldots\}$, we assume that \eqref{msmall} holds for   $k= m-1$,
and we will  demonstrate that \eqref{msmall} also holds for $k=m$.  

For each $(\tilde n,n)\in [K]^2$, we  denote 
 \beq\label{afexpree}
a_{\tilde n,n}:=B(f_{1,\tilde n},\ldots,f_{m-1,\tilde n},f_{m,n}).
\eeq
By  \cite[Lemma 8.1]{KMT22}, which gives the two-parameter Rademacher-Menshov inequality, we have  
\begin{align}
    &
\| (a_{N,N})_{N \in [K]} \|_{\V^2} \lesssim  \sum_{M_1,M_2 \in 2^{\N} \cap [K]} 
\big\| (\Delta a_{M_1 j_1,M_2j_2})_{(j_1,j_2) \in [K/M_1] \times [K/M_2]} \big\|_{\ell^2},
\label{eqnshizia}
\end{align} 
where $\Delta a_{M_1 j_1,M_2j_2}:= a_{M_1 j_1,M_2j_2}-a_{M_1j_1,M_2(j_2-1)}-a_{M_1(j_1-1),M_2j_2}+a_{M_1(j_1-1),M_2(j_2-1)}$ with  $a_{k_1,k_2}=0$ whenever $k_1k_2=0$. By  (\ref{afexpree}),  we expand the sum on the  RHS of (\ref{eqnshizia}) as 
\begin{align} 
\sum_{M_1,M_2 \in 2^\N \cap [K]}\Big( \sum_{j_1 \in [K/M_1]}&\sum_{j_2 \in [K/M_2]}
\big|
B(f_{1,M_1j_1},\ldots,f_{m-1,M_1j_1},\tilde{f}_{m,M_2j_2})\\
&\ -B(f_{1,M_1(j_1-1)},\ldots,f_{m-1,M_1(j_1-1)},\tilde{f}_{m,M_2j_2}) \big|^2\Big)^{1/2}
\label{eqn:differ11}
\end{align} 
with $\tilde{f}_{m,M_2j_2}:=f_{m,M_2j_2}-f_{m,M_2(j_2-1)}$. 
 Denote 
 \beq\label{reductionGGG1}
 G_{m,K,M_2}(t):=\sum_{j_2\in [K/M_2]}\varepsilon_{j_2}^{(m)}(t)\tilde{f}_{m,M_2j_2},
 \eeq
 where 
 $\{\varepsilon_{i}^{(m)}(t)\}_{i=0}^\infty$ is  the sequence of Rademacher functions (see, e.g.,  \cite[Appendix C]{GG14}) on $[0,1]$ with    
\beq\label{eqn:xinqing}
 \|\sum_{i=0}^\infty z_i \varepsilon_{i}^{(m)}(t)\|_{L^\rho_t([0,1])}\sim (\sum_{i=0}^\infty|z_i|^2)^{1/2}\quad {\rm for\ any}\ \rho\in(0,\infty).
 \eeq
Then, using (\ref{eqn:xinqing}) with $\rho=r\in(0, \min\{q,2\})$,  Minkowski's inequality and  (\ref{reductionGGG1}), 
 we obtain from (\ref{eqn:differ11}) that  
\beq\label{eqn:differ12}
\begin{aligned} 
&\ \| (a_{N,N})_{N \in [K]} \|_{\V^2}\\
\les&\  \langle \Log K\rangle \sum_{M_2 \in 2^\N \cap [K]}
\Big\| 
\Big(B\big(f_{1,N},\ldots,f_{m-1,N},G_{m,K,M_2}(t)\big)\Big)_{N\in [K]}\Big\|_{L^r_t([0,1];\V^2)}.
\end{aligned}
\eeq
Taking the $L^q(X)$ norm on the both sides of (\ref{eqn:differ12}), we deduce by  Minkowski's inequality (for $q\ge 1$) or quasi-triangle inequality (for $0<q<1$) that
\begin{equation}\label{msmall1}
\begin{split}
\| \big(a_{N,N}\big)_{N \in [K]}& \|_{L^q(X;\V^2)} 
\les_q \ \langle \Log K \rangle^{1+\max\{1,{1}/{q}\}} \\
&\ \times \sup_{M_2 \in 2^\N \cap [K]}\sup_{t\in [0,1]} \big\| \big(
B(f_{1,N},\ldots,f_{m-1,N},G_{m,K,M_2}(t))\big)_{N\in [K]}\big\|_{L^q(X;\V^2)}.
\end{split}
\end{equation}

Now, by the induction hypothesis that   (\ref{msmall}) holds for $ k= m-1$,
we conclude 
$$
\begin{aligned}
    &\ \big\| \big(
B(f_{1,N},\ldots,f_{m-1,N},G_{m,K,M_2}(t))\big)_{N\in [K]}\big\|_{L^q(X;\V^2)}
\les_q\  \langle \Log K \rangle^{m-3+(m-1)\max\{1,{1}/{q}\}}\\
&\ \times\sup \big\| B\big( \sum_{j \in [K]} \varepsilon_j^{(1)} (f_{1,j}-f_{1,j-1}),\ldots,\sum_{j\in [K]} \varepsilon_j^{(m-1)} (f_{m-1,j}-f_{m-1,j-1}),G_{m,K,M_2}(t)\big) \big\|_{L^q(X)}
\end{aligned}
$$
for all  $M_2 \in 2^\N \cap [K]$ and ${t\in [0,1]}$, where 
 the supremum is taken over $\varepsilon_j^{(i)}\in \{1,-1\}$ for all  $(i,j)\in [m-1]\times [K]$. This estimate, combined  with (\ref{msmall1}) and (\ref{reductionGGG1}), gives 
 (\ref{msmall}) for the case $k=m$. This completes the proof of Lemma \ref{mil}.
\end{proof}

\begin{proof}[Proof of Proposition \ref{t4.1}]
By (\ref{eqn:error1}), it suffices to show that for any $t\in [1/2,1]$,
\beq\label{eqn:reductionoperA}
\|\big((\ref{eqn:operA})\big)_{N\in \I_{\le}}\|_{\ell^q(\Z;\V^r)}\les 2^{-cl}\|f_1\|_{\ell^{q_1}(\Z)}\cdots \|f_k\|_{\ell^{q_k}(\Z)}
\eeq
with $r,q,q_1,\ldots,q_k$ given as in Proposition \ref{t4.1}. 

Normalize $\|f_i\|_{\ell^{q_i}(\Z)}=1$ for all $i\in [k]$.   
We can enumerate the elements of ${\I}_{\le}$ in order as $N_1 < \dots < N_K$, which by (\ref{notatII})  gives  $K = O( 2^u )$.  By Lemma \ref{mil},  we can  reduce matters to showing the inequality 
\beq\label{reducetion are1}
\|\tilde A_{2^u,\La_{2^u}}\big(\Pi_{l_1,\le -10 p_0u}(\tilde F_{t}^{l_1,s_1}),\ldots,\Pi_{l_k,\le -10 p_0u}(\tilde F_{t}^{l_k,s_k})\big)\|_{\ell^q(\Z)}
\les 2^{-cl}
\eeq
with the functions $\tilde F_{t}^{l_i,s_i}$ $(i\in [k])$ given by 
$$\tilde F_{t}^{l_i,s_i}:=T^{\mathfrak{Q}_{l_i}}_\Z[\eta_{*,t}^{s_i}]f_i,\quad {\rm where }\quad 
\eta_{*,t}^{s_i}=\sum_{j\in [K]}\varepsilon_j^{(i)}(\eta_{N_j,t}^{s_i}-\eta_{N_{j-1},t}^{s_i})
$$
with   
$\varepsilon_j^{(i)}\in \{1,-1\}$ 
and  $\eta_{N_{0},t}^{s_i}=0$. 
Since $\F_\Z \tilde F_{t}^{l_i,s_i}$ vanishes on 
${\mathfrak M}_{\le l_i-1,\le -d_i (u- l_i+1)}$ for each $i\in [k]$,  we deduce by Theorem \ref{Thm:uweyl1} and Theorem \ref{thmIW} that 
\beq\label{applymultiweyl}
\|\tilde A_{2^u,\La_{2^u}}\big(\Pi_{l_1,\le -10 p_0u}(\tilde F_{t}^{l_1,s_1}),\ldots,\Pi_{l_k,\le -10 p_0u}(\tilde F_{t}^{l_k,s_k})\big)\|_{\ell^q(\Z)}
\les 2^{-cl }\prod_{i\in [k]}\|\tilde F_{t}^{l_i,s_i}\|_{\ell^{q_i}(\Z)}.
\eeq
By the shifted 
 Calder\'on-Zygmund theory (see \cite[Theorem B.1]{KMT22}), we see that
 \beq\label{real-11}
\|T_\R[\eta_{*,t}^{s_i}]f_i\|_{L^{q_i}(\R)}\les \langle s \rangle^{O(1)} \|f_i\|_{L^{q_i}(\R)},\quad i\in [k].
 \eeq
 Then, applying (\ref{Iweq:376}) to each $\tilde F_{t}^{l_i,s_i}$,  in combination  with  \eqref{real-11}, gives   
 \beq\label{CEShifted}
\|\tilde F_{t}^{l_i,s_i}\|_{\ell^{q_i}(\Z)}
\les    \langle s \rangle^{O(1)} 2^{{\bf C}_{q_{i}}(2^{l})},\quad i\in [k].
 \eeq
Since $u\sim l^{5/4}$,  $C_*l>s$ and taking $l$ sufficiently large,  
 we finally  obtain (\ref{reducetion are1}) by combining  (\ref{CEShifted}) and (\ref{applymultiweyl}). This completes the proof of Proposition \ref{t4.1}.
\end{proof}
\section{Major
arc estimates for the low-frequency case at a large scale}
\label{section:large scale Adelic group}
In this section, we shall prove Proposition \ref{t4.2} by  employing  harmonic analysis of 
the adelic integers $\Ad_\Z$, and establishing  
 a crucial  arithmetic multilinear estimate (Theorem \ref{lend1}). We highlight that the arithmetic bilinear estimate in \cite{KMT22} applies only to small exponents, Theorem \ref{lend1} gives us a multilinear estimate valid for large exponents. 
\subsection{Reduction of Proposition \ref{t4.2}}
We first rewrite (\ref{eqn:operB})
by replacing the function  $\eta_u^*$ in \eqref{eqn:adeilc analysis} by
$$
\eta_u^{**}(\zeta)=\eta_{\le -p_0 2^u}(\zeta_1)\cdots  \eta_{\le -p_0 2^u}(\zeta_k),\quad \zeta=(\zeta_1,\ldots,\zeta_k)\in \R^k.
$$
This is possible since $N\in \I_{>}$ and it is important to do so in order to apply
the quantitative Shannon sampling theorem   (see Theorem \ref{Sampling} below).

We now introduce the adelic model functions  $f_{l_1,\ma},\ldots,f_{l_k,\ma}$ defined by 
$$f_{l_i,\ma}(x,y) :=\sum_{\A_i\in (\Q/\Z)_{l_i}}\int_\R \eta_{\le -p_0 2^u}(\xi_i)
\F_\Z f_i(\A_i+\xi_i)
e(-(\xi_i,\A_i)\cdot (x,y))d\xi_i,\qquad i\in [k],
$$
where $(x,y)\in \Ad_\Z=\R\times \hat{\Z}$, or equivalently on the Fourier side,
$$\F_{\ma_\Z} f_{l_i,\ma}(\xi_i,\A_i)={\ind {\Height(\A_i)=2^{l_i}}}~  \eta_{\le -p_0 2^u}(\xi_i)
~\F_\Z f_i(\A_i+\xi_i),\qquad i\in [k],
 $$
 where  $\xi_i\in \R$ and $\A_i\in  \Q/\Z$. For each $i\in [k]$, using (\ref{eq:111}), we have 
 $$
 \Pi_{l_i,\le -p_0 2^u}f_i=\Sample f_{l_i,\ma} \ \ \ {\rm and}
  \ \ \
 f_{l_i,\ma}=\Sample_{\R_{\le -p_0 2^u}\times (\Q/\Z)_{l_i}}^{-1}\Pi_{l_i,\le -p_0 2^u}f_i
$$
where the operators $\Sample$ and 
$\Sample^{-1}_\Omega$ 
are given by \eqref{interpolationop333} and 
 \eqref{interpolationoer}, respectively. Recall the notation $\R_{\le k}$ from \eqref{majorarcs22}.
 
Similarly, using  \eqref{eq:100}, \eqref{multioperator} 
and the notation in Subsection \ref{subbb:Terminor},  we can  write (\ref{eqn:operB}) for each $t\in [1/2,1]$ (with $\eta_u^{**}$ in place of $\eta_u^{*}$) as 
 \beq\label{defn:adeilc B}
 \Sample B_{\Ad_\Z}[1\otimes m_{l,\hat{\Z}^\times}](\mathfrak{F}^{l_1,s_1}_{N,t,\Ad},\cdots,\mathfrak{F}^{l_k,s_k}_{N,t,\Ad}),
\eeq
with  $\mathfrak{F}^{l_i,s_i}_{N,t,\Ad}$  and $m_{l,\hat{\Z}^\times}$ given by  
\begin{align}
    \mathfrak{F}^{l_i,s_i}_{N,t,\Ad}(x,y):=&\ \mathfrak{F}^{s_i}_{N,t}(f_{l_i,\Ad})(x,y)\ \ := T_{\Ad_\Z}[\eta_{N,t}^{s_i}\otimes 1]f_{l_i,\ma}(x,y)\quad {\rm and}\label{FF1}\\
   m_{l,\hat{\Z}^\times}(\theta):=&\ m_{l,\hat{\Z}^\times}(\theta_1,\ldots,\theta_k):=  G^\times(\theta) \prod_{i\in [k]} {\ind {(\Q/\Z)_{l_i}}}(\theta_i). \label{sumdiscrete1}
\end{align}
Since $N$ is only related to the continuous part $\R$ in the adelic group, the same argument 
for Proposition \ref{Prop:approxmajor1-again} shows that
\beq\label{eqn:adelic continuous}
\|\big(\mathfrak{F}^{l_i,s_i}_{N,t,\Ad}\big)_{N\in \I_{>}}\|_{L^q(\Ad_\Z;\V^r)}
\les \langle  s \rangle^{O(1)} 
\|f_{l_i,\ma}\|_{L^q(\Ad_\Z)}
\eeq
holds for any $r>2$ and $q\in(1,\infty)$.

Since  $N>2^{10 p_02^u}$,  the above functions $(f_{l_1,\ma},\ldots,f_{l_k,\ma})$  
and $(\mathfrak{F}^{l_1,s_1}_{N,t,\Ad},\ldots,\mathfrak{F}^{l_k,s_k}_{N,t,\Ad})$ defined   on $\ma_\Z$  have Fourier support in the region 
$$\big(\R_{\le -2^u}\times (\Q/\Z)_{l_1}\big)\times\cdots\times \big(\R_{\le-2^u}\times (\Q/\Z)_{l_k}\big)$$
which is contained in 
\beq\label{eqn:supporte1}
\big(\R_{\le -2^u}\times (\frac{1}{Q_{\le l}}\Z/\Z)\big)\times\cdots\times \big(\R_{\le-2^u}\times (\frac{1}{Q_{\le l}}\Z/\Z)\big)
\eeq
with $Q_{\le l}:={\rm lcm}(q\in\Z_+: q\in  [2^l])\le 2^{O(2^l)}$ ($l=\max\{l_1,\ldots,l_k\}$). Note by 
(\ref{eq:111}), $Q_{\le l}$ is smaller than $2^{2^{u/4}}$.  Thus,
in this large-scale situation,
we can apply Theorem \ref{Sampling} which together with (\ref{Iweq:376}), gives that  for each $q\in [p_0',p_0]$, 
\beq\label{inequk11}
\|f_{l_i,\ma}\|_{L^q(\ma_\Z)}\les 2^{{\bf C}_{p_0}(2^{l})}\|f_i\|_{\ell^q(\Z)}
\eeq
with ${\bf C}_{p_0}(\cdot)$ given by (\ref{constantspecia}). Furthermore, 
applying  Theorem \ref{Sampling} again (with the normed vector space $B=\V^r$), we have 
$$
\|\big((\ref{defn:adeilc B})\big)_{N\in \I_{>}}\|_{\ell^q(\Z;\V^r)}
\sim \|\big( B_{\Ad_\Z}[1\otimes m_{l,\hat{\Z}^\times}](\mathfrak{F}^{l_1,s_1}_{N,t,\Ad},\ldots,\mathfrak{F}^{l_k,s_k}_{N,t,\Ad})\big)_{N\in \I_{>}}\|_{\ell^q(\Ad_\Z;\V^r)}
$$
for $t\in [1/2,1]$,
which,  combined with 
 inequality  (\ref{inequk11}), shows that for  Proposition \ref{t4.2}, it suffices to find a $c>0$ such that  for any $t\in [1/2,1]$ and $r\in (kD_k,\infty)$, 
\beq\label{eqn:zuihou11}
 \begin{aligned}
 &\ \|\big( B_{\Ad_\Z}[1\otimes m_{l,\hat{\Z}^\times}](\mathfrak{F}^{l_1,s_1}_{N,t,\Ad},\ldots,\mathfrak{F}^{l_k,s_k}_{N,t,\Ad})\big)_{N\in \I_{>}}\|_{L^{D_k}(\Ad_\Z;\V^r)}
 \les\ 2^{-cl}\prod_{i\in [k]}\|f_{l_i,\ma}\|_{L^{kD_k}(\ma_\Z)}.
 \end{aligned}
 \eeq
\subsection{Proof of (\ref{eqn:zuihou11})}
Using the notation in \eqref{multioperator} and \eqref{sumdiscrete1}, we 
rewrite the multilinear operator $B_{\Ad_\Z}[1\otimes m_{l,\hat{\Z}^\times}]$ in \eqref{defn:adeilc B}  as
\beq\label{eqn:BAdfor1}
\sum_{\theta\in (\Q/\Z)^k}
\int_{\xi\in \R^k}
G^\times (\theta) \big(\prod_{i\in [k]} \F_{\ma_\Z}\mathfrak{F}^{l_i,s_i}_{N,t,\Ad} 
(\xi_i,\theta_i)\big) e\big(-(x,y)\cdot(\xi_1+\cdots+\xi_k,\theta_1+\cdots+\theta_k)\big)
d\xi.
\eeq
By the Fourier inversion formula in the continuous setting, (\ref{eqn:BAdfor1})  equals
\beq\label{eqn:BAdfor2}
\sum_{\theta_1\in  (\Q/\Z)}
\cdots \sum_{\theta_k\in  (\Q/\Z)}
G^\times (\theta) \big(\prod_{i\in [k]} \F_{\hat{\Z}}\mathfrak{F}^{l_i,s_i}_{N,t,\Ad} 
(x,\theta_i)\big) e\big(-y\cdot (\theta_1+\cdots+\theta_k)\big).
\eeq
Recall that $\hat{\Z}=\prod_{p\in \mathbb P}\Z_p$ denotes the compact group of profinite integers. Let  $\hat{\Z}^\times=\prod_{p\in \mathbb P}\Z_p^\times$ where  
 $\mathbb{Z}_p^\times$ denotes  the multiplicative group of $p$-adic units.
Define the multilinear averaging operator $A_{\hat{\Z}^\times}$ by
$A_{\hat{\Z}^\times}(g_1,\ldots,g_k)(x):= \int_{ \hat{\Z}^\times}~ 
g_1(x-P_1(y))\cdots g_k(x-P_k(y)) d\nu_{\hat{\Z}^\times}(y)$
 where $x\in \hat{\Z}$ and $\nu_{\hat{\Z}^\times}$ is Haar measure on $\hat{\Z}^\times$.
Equivalently, 
for any  $Q\in \Z_+$, and
any functions  $g_1,\ldots,g_k: \Z/Q\Z \to \mathbb{C}$  (which one can also view as functions on $\hat \Z$ in the obvious way), 
$$
A_{\hat{\Z}^\times}(g_1,\ldots,g_k)(x)= \E_{n \in (\Z/Q\Z)^\times}~ 
g_1(x-P_1(n))\cdots g_k(x-P_k(n)).
$$

As mentioned previously, changing  $(a,q)\to (Ka,Kq)$ for any $K\in\Z_+$ does not affect  the expression   (\ref{eqn:exponentialsum}) for $G^\times(\frac{a}{q})$.  
Thus,
by writing out $G^\times$ and using  the Fourier inversion formula in the arithmetic setting, we see that to bound (\ref{eqn:BAdfor2}), it suffices to estimate the function 
$A_{\hat{\Z}^\times}(\mathfrak{F}^{l_1,s_1}_{N,t,\Ad}(x,\cdot),\ldots,\mathfrak{F}^{l_k,s_k}_{N,t,\Ad}(x,\cdot))(y)$.
This allows us to  reduce the proof of  (\ref{eqn:zuihou11}) to proving  that there exists 
${R_\circ}\in (kD_k,r)$ such that
for any $t\in [1/2,1]$ and $x\in\R$,
\beq\label{eqn:zuihou1122}
 \begin{aligned}
 \|\big( A_{\hat{\Z}^\times}(\mathfrak{F}^{l_1,s_1,x}_{N,t,\Ad},\ldots,\mathfrak{F}^{l_k,s_k,x}_{N,t,\Ad})\big)_{N\in \I_{>}}\|_{L^{D_k}(\hat{\Z};\V^r)}
 \les\ 2^{-cl}\prod_{i\in [k]}\|\big(\mathfrak{F}^{l_i,s_i,x}_{N,t,\Ad}\big)_{N\in \I_{>}}\|_{L^{kD_k}(\hat\Z;\V^{R_\circ})}
 \end{aligned} 
 \eeq
 where    
$ \mathfrak{F}^{l_k,s_k,x}_{N,t,\Ad}(y):=\mathfrak{F}^{l_k,s_k}_{N,t,\Ad}(x,y)$ with $y\in \hat{\Z}$, and the implicit constant is independent of $x$ and $t$.
Indeed, (\ref{eqn:zuihou11}) follows from (\ref{eqn:zuihou1122}) by taking 
the $L^{D_k}(x\in \R) $ norm on  both sides of (\ref{eqn:zuihou1122}),  using H\"{o}lder's inequality  and (\ref{eqn:adelic continuous}) with $r=R_\circ$.

It remains to prove inequality (\ref{eqn:zuihou1122}). Before we give its proof, we need the following important multilinear estimate.
Recall that $D_k = 10 k d_k$ given in Proposition \ref{t4.2}.
\begin{thm}[Arithmetic multilinear estimate]\label{lend1}
Let $k\in\Z_+$ with $k\ge 2$,  ${\mathfrak l}\in\N$ and $D_k\le q_1,\ldots,q_k<\infty$.
Suppose 
$1\le \tilde r<\infty$ and 
 fix $j\in [k]$. 
There exist two positive constants  $c$ and $C$, depending on ${\mathcal P}, {\tilde r},q_1,\ldots,q_k$, such that 
\beq\label{eqn:arith1}
\|A_{\hat{\Z}^\times}(g_1,\ldots,g_k)\|_{L^{\tilde r}(\hat{\Z})}
\le C 2^{-c{\mathfrak l}}\|g_1\|_{L^{kq_1}(\hat{\Z})}\cdots \|g_k\|_{L^{kq_k}(\hat{\Z})}
\eeq
for any $g_1\in L^{kq_1}(\hat{\Z}),\ldots,g_k\in L^{kq_k}(\hat{\Z})$ with $\F_{\hat{\Z}} g_j$ vanishing  on $(\Q/\Z)_{\le {\mathfrak l}}$.
\end{thm}
\begin{remark}\label{rrr32:section arith}
In contrast to the arithmetic bilinear estimate in \cite{KMT22} for small exponents, Theorem \ref{lend1} establishes an arithmetic multilinear estimate for large exponents, enabling us not to impose any restrictions on  the polynomial map 
 $\mathcal{P}$ beyond those specified in  (\ref{eq:42}) and (\ref{eq:41}).
 Moreover, as in Theorem \ref{thm:padic} below, we do not specify all admissible $q_1,\ldots,q_k$
  explicitly in order to keep the presentation concise.
\end{remark}
\begin{proof}[Proof of (\ref{eqn:zuihou1122}) accepting Theorem \ref{lend1}]
Using the norm interchanging trick from  \cite[Lemma 9.5]{KMT22},  
for any $1 \leq R_1 < {r_1} \leq \infty$, we have the inequality 
\beq\label{interchange}
\| (f_k)_{k \in [K]} \|_{L^{r_1}(X; \V^{r_1})} \lesssim_{{r_1},R_1} \| (f_k)_{k \in [K]} \|_{\V^{R_1}([K]; L^{r_1}(X))},
\eeq
where  $X$ is any measure space.   
Let $R\in (kD_k,r)$. 
By H\"{o}lder's inequality (noting $D_k<r$) and (\ref{interchange}) with $(R_1,r_1)=(R,r)$,  we bound the left-hand side of  (\ref{eqn:zuihou1122}) by 
\beq\label{inter21}
\begin{aligned}
  &\  \|\big( A_{\hat{\Z}^\times}(\mathfrak{F}^{l_1,s_1,x}_{N,t,\Ad},\ldots,\mathfrak{F}^{l_k,s_k,x}_{N,t,\Ad})\big)_{N\in \I_{>}}\|_{L^r(\hat{\Z};\V^r)}\\
\les&\  \|\big( A_{\hat{\Z}^\times}(\mathfrak{F}^{l_1,s_1,x}_{N,t,\Ad},\ldots,\mathfrak{F}^{l_k,s_k,x}_{N,t,\Ad})\big)_{N\in \I_{>}}\|_{\V^R(\I_{>};L^r(\hat{\Z}))}. 
\end{aligned}
\eeq
By the definitions \eqref{var-seminorm} and \eqref{vardef},
to estimate the right-hand side of (\ref{inter21}), 
it  suffices to bound the following two expressions:   first, 
\beq\label{endF1}
\sup \big(\sum_{n\in [J-1]}\|A_{\hat{\Z}^\times}(\mathfrak{F}^{l_1,s_1,x}_{N_{n+1},t,\Ad},\ldots,\mathfrak{F}^{l_k,s_k,x}_{N_{n+1},t,\Ad})-
A_{\hat{\Z}^\times}(\mathfrak{F}^{l_1,s_1,x}_{N_n,t,\Ad},\ldots,\mathfrak{F}^{l_k,s_k,x}_{N_n,t,\Ad})\|_
{L^r(\hat{\Z})}^R\big)^{1/R},
\eeq
where  the supremum is applied over all  $J\in \Z_+$ and sequences $\{N_n\}_{n\in \N_{\le J}}\subset \I_\ge$; and  second, 
\beq\label{END22}
\| A_{\hat{\Z}^\times}(\mathfrak{F}^{l_1,s_1,x}_{N_\circ,t,\Ad},\ldots,\mathfrak{F}^{l_k,s_k,x}_{N_\circ,t,\Ad})\|_{L^r(\hat{\Z})}, 
\eeq
with  $N_\circ\in \I_>$. 
We note that the functions $\mathfrak{F}^{l_i,s_i,x}_{N,t,\Ad}$
appearing in \eqref{endF1} and \eqref{END22} have the property that
$\F_{\hat{\Z}}  \mathfrak{F}^{l_i,s_i,x}_{N,t,\Ad}$ is  supported  
on $(\Q/\Z)_{l_i}$ and in particular, vanishes on $(\Q/\Z)_{\le l_i-1}$ 
and so we may apply Theorem \ref{lend1} to bound both \eqref{endF1} and \eqref{END22}.

Indeed, the desired estimate  for  \eqref{END22} is a direct result of Theorem \ref{lend1} (with $\tilde r=r$ and  $q_i=D_k$ for $i\in [k]$). Also, by  Theorem  \ref{lend1} (again with $\tilde r=r$ and $q_i=D_k$ for $i\in [k]$) and  a   routine  computation, we see that (\ref{endF1}) is 
\begin{align}
 \ \les\ 2^{-cl }\sum_{j\in [k]} \Big\{ 
   \sup\big(\sum_{n\in [J-1]}\|
   \bar{\mathfrak{F}}^{l_j,s_j,x}_{N_n,t,\Ad}
   \|_{L^{kD_k}(\hat{\Z})}^R\big)^{1/R}\Big\}
   \prod_{i\in [k]\setminus\{j\}}  \|\sup_{N\in \I_{>}}|\mathfrak{F}^{l_i,s_i,x}_{N,t,\Ad}|\|_{L^{kD_k}(\hat{\Z})}\label{eqn:bound1}
\end{align}
with $\bar{\mathfrak{F}}^{l_j,s_j,x}_{N_n,t,\Ad}:=\mathfrak{F}^{l_j,s_j,x}_{N_{n+1},t,\Ad}-\mathfrak{F}^{l_j,s_j,x}_{N_n,t,\Ad}$. Recall  $R\in (kD_k,r)$.
 Applying Minkowski's inequality to  inequality (\ref{eqn:bound1}), and using 
 \eqref{vardef}, 
 we see that
(\ref{eqn:zuihou1122}) follows with $R_\circ=R$, finishing the proof of 
(\ref{eqn:zuihou1122}) and so completing the proof of Proposition \ref{t4.2}.
\end{proof}

It remains to prove Theorem   \ref{lend1}.   
Let 
$$A_{\Z_p}(f_1, \ldots, f_k)(x):=\int_{ \Z_p}f_1(x-P_1(y))\cdots f_k(x-P_k(y))d\nu_{\Z_p}(y)$$
be a multilinear polynomial averaging operator acting on functions defined over the $p$-adic integers $\Z_p$ where $\nu_{\mathbb{Z}_p}$ is Haar measure on $\mathbb{Z}_p$, normalized to be a probability measure.
See Section \ref{section:local field}  for further  details. 
Abusing notation slightly,  for all $k\in\Z_+$, 
we  denote  $A_{\Z_p}(f_1, \ldots, f_k)$ as 
\beq\label{abuse-not1}
A_{\Z_p}(f_1, \ldots, f_k)(x):=\E_{n\in \Z_p}f_1(x-P_1(n))\cdots f_k(x-P_k(n))
\eeq
with $\deg P_i=d_i$ for $i\in [k]$ obeying  \eqref{eq:41}.  
\begin{thm}[$p$-adic multilinear estimate]\label{thm:padic}
Let $k\in\Z_+$ with $k\ge 2$,    
$1\le q\le \infty$ and 
  $D_k\le q_1,\ldots,q_k<\infty$.
Then we have   $\|A_{{\Z_p}}\|_{L^{q_1}(\Z_p)\times\cdots \times L^{q_k}(\Z_p) 
\to L^{q}(\Z_p)}\les_{\mathcal P}  1$, where the implicit constant is independent of the prime  $p$.
\end{thm}
The proof of Theorem \ref{thm:padic} is deferred to Section \ref{section:local field}. In fact, we will prove a more general result (see Theorem \ref{thm:padiccccc}) 
for all primes $p$ with extended parameter ranges. Notably, the endpoint case $q=\infty$
will be crucial; see the proof of (\ref{eqn:arith111}) below. Moreover, the distinct degree condition \eqref{eq:41} is not needed for Theorem \ref{thm:padic}.
\subsection{Proof of Theorem  \ref{lend1}}
By interpolation and H\"{o}lder's inequality,  we reduce  the proof of Theorem  \ref{lend1} to establishing  two key    estimates. The first estimate states that  
 for each $1<\tilde{q}_1,\ldots,\tilde{q}_k<\infty$ with 
 ${1}/{\tilde{q}_1}+\cdots+{1}/{\tilde{q}_k}={1}/{\tilde{q}}\le1$, there exists $c>0$ such that
\beq\label{eqn:arith3321}
\|A_{\hat{\Z}^\times}(g_1,\ldots,g_k)\|_{L^{\tilde{q}}(\hat{\Z})}
\les  2^{-c{\mathfrak l}} \|g_1\|_{L^{\tilde{q}_{1}}(\hat{\Z})}\cdots \|g_k\|_{L^{\tilde{q}_{k}}(\hat{\Z})}
\eeq
for any $g_1\in L^{\tilde{q}_{1}}(\hat{\Z}),\ldots,g_k\in L^{\tilde{q}_{k}}(\hat{\Z})$ such that  $\F_{\hat{\Z}} g_j$ vanishing  on $(\Q/\Z)_{\le {\mathfrak l}}$, while   the second estimate states that for each $1\le r<\infty$ and $D_k\le q_1,\ldots,q_k<\infty$, 
\beq\label{eqn:arith111}
\|A_{\hat{\Z}^\times}(g_1,\ldots,g_k)\|_{L^{ r}(\hat{\Z})}
\les  \|g_1\|_{L^{kq_1}(\hat{\Z})}\cdots \|g_k\|_{L^{kq_k}(\hat{\Z})}
\eeq
for any $g_1\in L^{k {q}_{1}}(\hat{\Z}),\ldots,g_k\in L^{k {q}_{k}}(\hat{\Z})$.
Here, the frequency support condition $\F_{\hat{\Z}} g_j=0$  on $(\Q/\Z)_{\le \rm l}$ applies only to the first estimate. The value of the 
positive constant $c$ may vary in each appearance. 

\begin{proof}[Proof of (\ref{eqn:arith3321})]
The proof adapts the approach from \cite{KMT22} to fit our setting.  
By Minkowski's inequality and H\"{o}lder's inequality, we can establish  (\ref{eqn:arith3321}) without the factor $2^{-c{\mathfrak l}}$. 
Consequently, through interpolation,  it  suffices to prove  (\ref{eqn:arith3321}) for the special  case where $(q,q_1,\ldots,q_k)=(1,k,\ldots,k)$.
By a standard limiting argument, we may assume that  the functions $\{g_i\}_{i\in [k]}$  on 
$\hat{\Z}$ are periodic modulo $Q$, that is, they factor through $\Z/Q\Z$. The matter therefore reduces to establishing the bound
\beq\label{eqn:pinxi}
\|A_{(\Z/Q\Z)^\times}(g_1,\ldots,g_k)\|_{L^1(\Z/Q\Z)}
\les 2^{-c{\mathfrak l}} \|g_1\|_{L^{k}(\Z/Q\Z)}\cdots\|g_k\|_{L^{k}(\Z/Q\Z)}
\eeq
for some $c>0$ under the assumption that   $\F_{\Z/Q\Z} g_j$ vanishes on $(\Q/\Z)_{\le {\mathfrak l}} \cap (\frac{1}{Q}\Z/\Z)$.
We may normalize $\|g_i\|_{L^{k}(\Z/Q\Z)}=1$ for all $i\in [k]$. 
 
Let 
$N$ and $R$ be two large natural numbers, which 
should be  considered  
as 
  being much larger than ${\mathfrak l}$ and $ Q$.
 Define  $g_{1,R},\ldots,g_{k,R} \in \Schwartz(\Z)$ by 
$$g_{i,R}(n):=R^{-1/{k}}\phi(n/R) g_i(n\ {\rm mod}\ Q),\qquad i\in [k],$$
where $\phi\in \Schwartz(\R)$ is a non-negative, even function with
$\int_\R \phi^k(x)dx=1$ and with Fourier support in $[-1,1]$. One has $g_{i,R}\in \ell^{2}(\Z)\cap \ell^{k}(\Z)$
and $\F_\Z g_{j,R} $ is supported on 
$\{\xi + \A : \xi \in [-1/R,1/R] \ {\rm and} \ \F_{\Z/Q\Z} g_j (\A)\neq 0\}$.
Hence  $\F_\Z g_{j,R} $ vanishes on $\mathfrak{M}_{\le {\mathfrak l},\le -d_i\Log N+d_i{\mathfrak l}}$, which by Theorem \ref{Thm:uweyl1},  gives 
\beq\label{OOeqn1}
\|A_{N,\La_N}(g_{1,R},\ldots,g_{k,R})\|_{\ell^1(\Z)}\les 2^{-c{\mathfrak l}} \|g_{1,R}\|_{\ell^{k}(\Z)}\cdots \|g_{k,R}\|_{\ell^{k}(\Z)}.
\eeq
 Using Riemann sums of $\phi^k$, we see that
$$\lim_{R\to \infty}\|g_{i,R}\|_{\ell^{k}(\Z)}=\|g_i\|_{L^{k}(\Z/Q\Z)}=1,\qquad i\in [k].$$
This, combined with \eqref{OOeqn1}, leads to 
\beq\label{RHSlim}
\limsup_{N\to \infty}\limsup_{R\to \infty}\|A_{N,\La_N}(g_{1,R},\ldots,g_{k,R})\|_{\ell^1(\Z)}\les
2^{-c{\mathfrak l}}. 
\eeq
 By splitting into residue classes mod $Q$, we can write   $A_{N,\La_N}(g_{1,R},\cdots,g_{k,R})(x)$
as
\beq\label{OO890}
\begin{aligned}
\sum_{a\in [Q]}
 \big(\prod_{i\in [k]} g_i(x~ {\rm mod}~Q-P_i(a))\big)~ \E_{n\in [N]} \Big(\La_N(n) {\ind {n\equiv a ({\rm mod}~ Q) }} \big\{R^{-1} \prod_{i\in [k]} \phi
 (\frac{x-P_i(n)}{R})\big\}\Big).
 \end{aligned}
\eeq
By the mean value theorem (applied to $\phi(\cdot/R)$) and \eqref{approx;resides}, 
we rewrite (\ref{OO890})  as 
\beq\label{OPeqn}
\begin{aligned}
 R^{-1}\phi^k(x/R)\big\{ 
A_{(\Z/Q\Z)^\times}(g_1,\ldots,g_k)(x~ {\rm mod}~Q)
 &+ O(Q^{O(1)} \exp(-c\Log^{4/5}N)\big\}\\
 &+O(N^{O(1)}R^{-2}\langle x/R\rangle^{-2}).
\end{aligned}
\eeq
Taking the $\ell^1(x\in\Z)$ norm on both sides of (\ref{OPeqn}), splitting the sum into residue classes mod $Q$ and using Riemann sums of $\phi^k$, we obtain    
$$\limsup_{N\to \infty} \limsup_{R\to \infty}\|A_{N,\La_N}(g_{1,R},\ldots,g_{k,R})\|_{\ell^1(\Z)}
=\|A_{(\Z/Q\Z)^\times}(g_1,\cdots,g_k)\|_{L^1(\Z/Q\Z)},
$$
which, together  with (\ref{RHSlim}),  leads to  the desired (\ref{eqn:pinxi}).
This finishes the proof of (\ref{eqn:arith3321}).
\end{proof}
\begin{proof}[Proof of (\ref{eqn:arith111})]
  {\bf Claim}: for any $1\le  r<\infty$ and any $D_k\le q_1,\ldots,q_k<\infty$, 
\beq\label{As1}
\begin{cases}
\arraycolsep=1.5pt
\begin{array}{ccl}
 \|A_{\hat{\Z}^\times}\|_{L^{2q_1}(\hat{\Z})\times L^{2q_2}(\hat{\Z})\times L^{\infty}(\hat{\Z})\times L^{\infty}(\hat{\Z})\times\cdots \times L^{\infty}(\hat{\Z})
\to L^{r}(\hat{\Z})}&\les_{\mathcal{P}} &  1,\\ 
\|A_{\hat{\Z}^\times}\|_{L^{\infty}(\hat{\Z})\times L^{2q_2}(\hat{\Z})\times L^{2q_3}(\hat{\Z})\times L^{\infty}(\hat{\Z})\times\cdots \times L^{\infty}(\hat{\Z})
\to L^{r}(\hat{\Z})}&\les_{\mathcal{P}} &  1,\\ 
\|A_{\hat{\Z}^\times}\|_{L^{\infty}(\hat{\Z})\times L^{\infty}(\hat{\Z})\times L^{2q_3}(\hat{\Z})\times L^{2q_4}(\hat{\Z})\times\cdots \times L^{\infty}(\hat{\Z})
\to L^{r}(\hat{\Z})}&\les_{\mathcal{P}} &  1,\\
 \vdots & \vdots & \\
\|A_{\hat{\Z}^\times}\|_{L^{2q_1}(\hat{\Z})\times L^{\infty}(\hat{\Z})\times L^{\infty}(\hat{\Z})\times L^{\infty}(\hat{\Z})\times\cdots \times L^{2q_k}(\hat{\Z})
\to L^{r}(\hat{\Z})}&\les_{\mathcal{P}} &  1.
\end{array}
\end{cases}
\eeq
By interpolation, the claim establishes 
(\ref{eqn:arith111}). Thus it remains to prove  \eqref{As1}.
 By H\"{o}lder's inequality, we further  reduce  the matter to showing the  estimate for the bilinear average $A_{{\hat \Z^\times}}^{P_i,P_j}$;
 that is, 
$$\|A_{{\hat\Z^\times}}^{P_i,P_j}\|_{L^{2q_i}(\hat{\Z})\times L^{2q_j}(\hat{\Z})\to L^{r}(\hat{\Z})}\les_{\mathcal{P}} 1$$ 
for any $i\neq j$,  $1\le r<\infty$ and 
$D_k\le q_i,q_j<\infty$. 

Without loss of generality, we only show the details for  the case $(i,j)=(1,2)$ since other cases can be treated similarly.  
By approximating $\hat {\Z}$ (and their unit group ${\hat \Z}^\times$) by the product of finitely many of the $p$-adic groups $\Z_p$ (and $\Z_p^\times$ respectively),  and following the argument in \cite[p.1084]{KMT22},  it suffices to  establish two key $p$-adic estimates for the averaging operator $A_{{\Z_p^\times}}^{\tilde{\mathcal{P}}}(g_1,g_2):= A_{{\Z_p^\times}}^{P_1,P_2}(g_1,g_2)$: first, 
\beq\label{and1}
\|A_{{\Z_p^\times}}^{\tilde{\mathcal{P}}}\|_{L^{{2q_1}}({\Z_p})  \times L^{{2q_2}}({\Z_p})  \to L^{r}({\Z_p})}\les_{P_1,P_2}  {\frac{p}{p-1}}
\eeq
for all primes $p$;   
and second, the improved uniform bound
\beq\label{and2222}
\|A_{{\Z_p^\times}}^{\tilde{\mathcal{P}}}\|_{L^{{2q_1}}({\Z_p})  \times L^{{2q_2}}({\Z_p})  \to L^{r}({\Z_p})}\le 1
\eeq
for all sufficiently large primes   $p$. 

Note that ${2q_1},{2q_2}\ge 2D_k$
since $q_1,q_2\ge D_k$.
By using  Theorem \ref{thm:padic} (in the bilinear case)
and the basic inequality
$|A_{{\Z_p^\times}}^{\tilde{\mathcal{P}}}(g_1,g_2)|\le \frac{p}{p-1}A_{{\Z_p}}^{\tilde{\mathcal{P}}}(|g_1|,|g_2|)$, 
we conclude 
\beq\label{AAAand1}
\|A_{{\Z_p^\times}}^{\tilde{\mathcal{P}}}\|_{L^{{2q_1}}({\Z_p})  \times L^{{2q_2}}({\Z_p}) \to L^{r}({\Z_p})}\les_{P_1,P_2}  {\frac{p}{p-1}},
\eeq
which gives
\eqref{and1} immediately. On the other hand, since  $|A_{{\Z_p^\times}}^{\tilde{\mathcal{P}}}(g_1,g_2)|\le A_{{\Z_p^\times}}^{\tilde{\mathcal{P}}}(|g_1|,|g_2|)$,
to prove \eqref{and2222}, it suffices to show that 
for any $2<r<\infty$, $D_k\le q_1,q_2<\infty$ and for sufficiently large primes $p$,
\beq\label{and2}
\|A_{{\Z_p^\times}}^{\tilde{\mathcal{P}}}(g_1,g_2)\|_{L^r(\Z_p)}
\le \|g_1\|_{L^{{2q_1}}({\Z_p})} \|g_2\|_{L^{{2q_2}}({\Z_p})}
\eeq
holds for  any non-negative $g_1\in L^{{2q_1}}({\Z_p}),g_2\in L^{{2q_2}}({\Z_p})$.

For  $i=1,2$,  we normalize   $\|g_i\|_{L^{2q_i}(\Z_p)}=1$, 
set $a_i:=\E_{n\in\Z_p}g_{i}$, and $g_{i,0}:=g_i-a_i$; a simple computation gives that   $\E_{n\in\Z_p}g_{i,0}(n)=0$ and $|a_i|\le\|g_i\|_{L^{2}(\Z_p)}\le 
1$ (by H\"{o}lder's inequality). Note that the only frequencies in ${\mathbb Z}_p^{*}$ of height at most $p$ is the zero frequency and since $g_{i,0}$ has mean zero, we can argue as in  (\ref{eqn:arith3321}) to conclude that there
exists a $\gamma_0>0$ such that  
\beq\label{eqn:padicdeca1}
\begin{aligned}
 \|A_{{\Z_p^\times}}^{\tilde{\mathcal{P}}}(g_{1,0},a_2)\|_{L^{2}(\Z_p)}
 + \|A_{{\Z_p^\times}}^{\tilde{\mathcal{P}}}(a_{1},g_{2,0})\|_{L^{2}(\Z_p)} \les&\  p^{-\gamma_0} (\|g_{1,0}\|_{L^{2}(\Z_p)}+\|g_{2,0}\|_{L^{2}(\Z_p)})\quad {\rm and }\\
  \|A_{{\Z_p^\times}}^{\tilde{\mathcal{P}}}(g_{1,0},g_{2,0})\|_{L^{1}(\Z_p)} \les&\  p^{-\gamma_0} \|g_{1,0}\|_{L^{2}(\Z_p)}\|g_{2,0}\|_{L^{2}(\Z_p)}.
  \end{aligned}
\eeq
Next, we will use \eqref{eqn:padicdeca1} to obtain \eqref{and2}. 
Denote 
$$B_1:=\|g_{1,0}\|_{L^{2}(\Z_p)}^2 \qquad {\rm and }\qquad B_2:=\|g_{2,0}\|_{L^{2}(\Z_p)}^2.$$
 Then   
$|a_i|^2=\|g_{i}\|_{L^{2}(\Z_p)}^2-B_i\le 1-B_i$ for $i=1,2$.
Note that  
the function $x \mapsto |x|^{r}$   is continuously twice differentiable. 
By Taylor's expansion, we have
$$
\begin{aligned}
|A_{{\Z_p^\times}}^{\tilde{\mathcal{P}}}(g_1,g_2)|^r
=&\ |a_1 a_2|^r +r |a_1a_2|^{r-1}L_k(g_1,g_2)\\
&\ +O_r(1)\big| a_1 a_2+\theta L_k(g_1,g_2)\big|^{r-2} \big(L_k(g_1,g_2)\big)^2,
\end{aligned}$$
for some $\theta\in (0,1)$, 
where
\beq\label{eqn:enddecom941}
L_k(g_1,g_2):=A_{{\Z_p^\times}}^{\tilde{\mathcal{P}}}(g_1,g_2)
-a_1a_2=A_{{\Z_p^\times}}^{\tilde{\mathcal{P}}}(g_{1,0},a_2)+A_{{\Z_p^\times}}^{\tilde{\mathcal{P}}}(a_1,g_{2,0})
+A_{{\Z_p^\times}}^{\tilde{\mathcal{P}}}(g_{1,0},g_{2,0}).
\eeq
As a consequence, we can deduce  
\beq\label{eqn:May40}
\begin{aligned}
    \|A_{{\Z_p^\times}}^{\tilde{\mathcal{P}}}(g_1,g_2)\|_{L^r(\Z_p)}^r
\le&\  |a_1a_2|^2 + O_r\big(\big|\E_{n\in \Z_p} L_k(g_1,g_2)(n)\big|\big)\\
&\ +O_r\big(\|L_k(g_1,g_2)\|_{L^2(\Z_p)}^2+\|L_k(g_1,g_2)\|_{L^r(\Z_p)}^r\big).
\end{aligned}
\eeq
Using \eqref{AAAand1}, \eqref{eqn:enddecom941} and $0\le a_i\le \|g_i\|_{L^{2q_i}(\Z_p)}=1$ $(i=1,2)$, we have
$$
\begin{aligned}
    \|L_k(g_1,g_2)\|_{L^\infty(\Z_p)}\le&\  \|A_{{\Z_p^\times}}^{\tilde{\mathcal{P}}}(g_{1,0},a_2)\|_{L^\infty(\Z_p)}+\|A_{{\Z_p^\times}}^{\tilde{\mathcal{P}}}(a_1,g_{2,0})\|_{L^\infty(\Z_p)}
+\|A_{{\Z_p^\times}}^{\tilde{\mathcal{P}}}(g_{1,0},g_{2,0})\|_{L^\infty(\Z_p)}\\
\les&\ \|g_{1,0}\|_{L^{{2q_1}}({\Z_p})}\|a_2\|_{L^{{2q_2}}({\Z_p})}+ \|a_1\|_{L^{{2q_1}}({\Z_p})}\|g_{2,0}\|_{L^{{2q_2}}({\Z_p})}\\
&\ +\|g_{1,0}\|_{L^{{2q_1}}({\Z_p})} \|g_{2,0}\|_{L^{{2q_2}}({\Z_p})}\\
\les&\  \|g_{1}\|_{L^{{2q_1}}({\Z_p})}\|g_2\|_{L^{{2q_2}}({\Z_p})}
\les\ 1,
\end{aligned}$$
and then 
$$ \|L_k(g_1,g_2)\|_{L^r(\Z_p)}^r
\le \|L_k(g_1,g_2)\|_{L^2(\Z_p)}^2\|L_k(g_1,g_2)\|_{L^\infty(\Z_p)}^{r-2}
\les_r \|L_k(g_1,g_2)\|_{L^2(\Z_p)}^2,
$$
which, with \eqref{eqn:May40}, gives 
\beq\label{eqn:May444}
\begin{aligned}
    \|A_{{\Z_p^\times}}^{\tilde{\mathcal{P}}}(g_1,g_2)\|_{L^r(\Z_p)}^r
\le&\ |a_1a_2|^2 + O_r\big(\big|\E_{n\in \Z_p} L_k(g_1,g_2)(n)\big|+\|L_k(g_1,g_2)\|_{L^2(\Z_p)}^2\big).
\end{aligned}
\eeq
A simple computation shows 
\beq\label{eqn:may3}
|a_1 a_2|^2\le (1-B_1)(1-B_2)\le 1-\frac{B_1+B_2}{2}.
\eeq
Since 
  $\E_{n\in\Z_p}g_{i,0}(n)=0$,   
 we   deduce  from \eqref{eqn:enddecom941} and    (\ref{eqn:padicdeca1}) that 
 \beq\label{ABAB2}
 |\E_{n\in \Z_p} L_k(g_1,g_2)(n)|=\| A_{{\Z_p^\times}}^{\tilde{\mathcal{P}}}(g_{1,0},g_{2,0})\|_{L^1(\Z_p)}
 \les \ p^{-\gamma_0} B_1^{1/2}B_2^{1/2}.
 \eeq
 In addition, by applying \eqref{eqn:enddecom941}, Minkowski's inequality, 
  \eqref{eqn:padicdeca1},  H\"{o}lder's inequality and \eqref{AAAand1}, we can deduce that $\|L_k(g_1,g_2)\|_{L^2(\Z_p)}^2$ is  
\beq\label{eqn:may1}
\begin{aligned}
\les&\ \|A_{{\Z_p^\times}}^{\tilde{\mathcal{P}}}(g_{1,0},a_2)\|_{L^2(\Z_p)}^2+\|A_{{\Z_p^\times}}^{\tilde{\mathcal{P}}}(a_1,g_{2,0})\|_{L^2(\Z_p)}^2
+\|A_{{\Z_p^\times}}^{\tilde{\mathcal{P}}}(g_{1,0},g_{2,0})\|_{L^2(\Z_p)}^2\\
\les&\ p^{-\gamma_0}(\|g_{1,0}\|_{L^2(\Z_p)}^2+\|g_{2,0}\|_{L^2(\Z_p)}^2)
+\|A_{{\Z_p^\times}}^{\tilde{\mathcal{P}}}(g_{1,0},g_{2,0})\|_{L^1(\Z_p)}\|A_{{\Z_p^\times}}^{\mathcal P}(g_{1,0},g_{2,0})\|_{L^\infty(\Z_p)}\\
\les&\ p^{-\gamma_0}(\|g_{1,0}\|_{L^2(\Z_p)}^2+\|g_{2,0}\|_{L^2(\Z_p)}^2
+\|g_{1,0}\|_{L^2(\Z_p)}\|g_{2,0}\|_{L^2(\Z_p)})\\
\les&\ p^{-\gamma_0} (B_1+B_2+B_1^{1/2}B_2^{1/2}).
\end{aligned}
\eeq
 Finally, 
substituting (\ref{eqn:may3})-\eqref{eqn:may1} into (\ref{eqn:May444}) gives
us a constant $C_r>0$  such that 
\beq\label{dne21}
\|A_{{\Z_p^\times}}^{\tilde{\mathcal{P}}}(g_1,g_2)\|^{r}_{L^{r}(\Z_p)}
\le 1-\frac{B_1+B_2}{2}+C_r p^{-{\gamma}_0}(B_1+B_2).
\eeq
Set the prime $p$  large enough such that $C_r p^{-{\gamma}_0}\le 1/4$. Then \eqref{and2} follows from \eqref{dne21}.
\end{proof}

\section{Multilinear polynomial averages over  {$p$}-adic fields}
\label{section:local field}
In this section, we establish an $L^{\rm p}$-improving inequality for polynomial averaging operators defined over $p$-adic fields ${\mathbb Q}_p$ 
for all  primes $p$. The field ${\mathbb Q}_p$ is the prototypical example of a {\it nonarchimedean} local field which is
the field of fractions of the ring of $p$-adic integers ${\mathbb Z}_p$ encountered above.
Let $|\cdot|_p$ (or simply $|\cdot|$) denote the usual $p$-adic absolute value so the associated unit ball $B_1(0) = \{x \in {\mathbb Q}_p : |x| \le 1\} = {\mathbb Z}_p$ is the ring of $p$-adic integers. Finally, let $\nu$ denote normalised Haar measure on ${\mathbb Q}_p$ so that $\nu(B_1(0)) = 1$. 
See \cite{N} for a comprehensive treatment of $p$-adic fields and local fields in general.

Let $m\in\Z_+$, and   consider a collection of polynomials 
$\mathcal{P}:=\{P_1,\ldots,P_m\}\subset  {{\mathbb Q}_p}[{\rm y}]$
 satisfying $1\le d_1\le d_2\le \cdots\le d_m$, where 
$d_i:=\deg P_i$ for $i\in [m]$. Define  the multilinear averaging operator 
$A:=A_{{\mathbb Q}_p}$ as 
$$
A(f_1, \ldots, f_m)(x):= \int_{\Z_p} f_1(x - P_1(y)) \cdots f_m(x - P_m(y)) \, d\nu(y),\quad x\in {\mathbb Q}_p.
$$
If $d_m = 1$ and $1/p_1 + \cdots + 1/p_m = 1$, H\"older's inequality implies the $L^{\mathrm{p}}$-improving bound
$$
\|A(f_1, \ldots, f_m)\|_{L^{\infty}({\mathbb Q}_p)} \ \lesssim_{\mathcal P} \ \|f_1\|_{L^{p_1}({\mathbb Q}_p)} \cdots \|f_m\|_{L^{p_m}({\mathbb Q}_p)}.
$$
By interpolation, 
this, combined with the trivial bound $L^{p_1}({{\mathbb Q}_p}) \times \cdots \times L^{p_m}({{\mathbb Q}_p}) \to L^1 ({{\mathbb Q}_p})$, establishes the $L^{p_1}({{\mathbb Q}_p}) \times \cdots \times L^{p_m}({{\mathbb Q}_p}) \to L^q ({{\mathbb Q}_p})$ bounds for all $1 \le q \le \infty$.
Hence,  we focus on the case  $d_m\ge 2$ in what follows.

We are interested in  establishing  
$L^{p_1}({{\mathbb Q}_p}) \times \cdots \times L^{p_m}({{\mathbb Q}_p}) \to L^q ({{\mathbb Q}_p})$ (or simply $L^{p_1} \times \cdots \times L^{p_m} \to L^q $) bounds for the operator $A$; that is,  
\begin{equation}\label{bound}
\|A(f_1, \ldots, f_m)\|_{L^q({\mathbb Q}_p)} \lesssim_{\mathcal P} \|f_1\|_{L^{p_1}({\mathbb Q}_p)} \cdots \|f_k\|_{L^{p_m}({\mathbb Q}_p)}.
\end{equation}
Let us introduce  three parameters $p_*,K$ and $q_{*,K}$ given by
$$
\frac{1}{p_{*}}:= \frac{1}{p_1} + \cdots + \frac{1}{p_m},\qquad
K := \frac{d_1}{p_1}  + \cdots  +  \frac{d_m}{p_m}
$$
and 
$$
q_{*,K}:=
\begin{cases}
\infty & \text{ if } K\le 1,\\
d_m/(K-1)& \text{ if } K>1.
\end{cases}
$$

The following is the main result of this section.
\begin{thm}[$L^{\rm p}({\mathbb Q}_p)$-improving estimate]\label{thm:padiccccc}
Let $m\in\Z_+$ with $m\ge 2$. Suppose $d_m\ge 2$,  $1< p_1,\ldots,p_{m}<\infty$, $1\le p_*<\infty$
and 
\begin{equation}\label{condition}
 \frac{d_1}{p_1} + \cdots + \frac{d_{m-1}}{p_{m-1}}< 1.
\end{equation}
Then \eqref{bound} holds 
for all $p_*\le q < q_{*,K}$ with an implicit constant that does not depend on the prime $p$.
Furthermore, if $K<1$ or if $K>1$ and $p_{*} > 1$, then 
\eqref{bound} holds for $p_*\le q \le q_{*,K}$.
\end{thm}
We give some remarks about Theorem~\ref{thm:padiccccc}.
\begin{itemize}
\item
 The following example demonstrates    the relevance of $K$ and the upper bound   of $q$ when $K>1$. Recall that $\Z_p$ is the unit ball in $\Q_p$.
 Fix a sufficiently small $\delta>0$, let $P_j(y) = y^{d_j}$,  and define the sets 
$E_j:=\{x \in \Z_p: |x| \le 2 \delta^{d_j}\}$
for $j\in [m]$. Then for any $x \in \Z_p$ with $|x| \le \delta^{d_m}$,
$$
A({\ind {E_1}}, \ldots, {\ind {E_m}})(x) \ge \int_{|y|\le \delta}
{\ind {E_1}}(x - y^{d_1}) \cdots {\ind {E_m}}(x- y^{d_m}) \, d\nu(y)  \sim  \delta.
$$
Consequently, if inequality  \eqref{bound} holds, we obtain 
\begin{equation*}
\delta^{1 + d_m/q} \ \lesssim_{\mathcal P} \ \delta^{K}, \ \ {\rm which\ \  implies} \ \ q \ \le \
\frac{d_{m}}{K - 1}= q_{*,K}\ (K>1).
\end{equation*}
This  
shows that the upper bound  $ q \le q_{*,K}$ in Theorem \ref{thm:padiccccc} is sharp.

\smallskip

\item 
 Theorem \ref{thm:padiccccc} recovers the $L^{\rm p}$-improving bound from \cite[Section 10]{KMT22}, proving that  for  the special  case $m=2$, $d_1 = 1$ and $d_2 \ge 2$,  inequality  \eqref{bound} holds 
 when 
$p_1 = p_2 = 2$ and $1 \le q < 2 d_2/(d_2 - 1) = q_{*,K}$ with $K=(d_2+1)/2$.

\item 
Since the linear averages ($m=1$) can be handled similarly (and indeed 
simpler), we focus on the multilinear case. Furthermore, we can omit considering endpoint cases where $p_k=\infty$ for some $k\in [m]$ due to
the inequality  $|A\,(f_1, \ldots, f_m)|\le \|f_j\|_{L^\infty({\mathbb Q}_p)}\,A\,(|f_1|, \ldots, 
|f_{j-1}|,|f_{j+1}|,\ldots,|f_m|)$ for any $j\in [m]$.  
\end{itemize}
\begin{proof}[Proof of Theorem \ref{thm:padic} (accepting Theorem \ref{thm:padiccccc})]
Applying Theorem \ref{thm:padiccccc} 
with  
$(m,p_1,\ldots,p_m,q)=(k,q_1,\ldots,q_k,q)$ (so we have $K<1$),
we can deduce  
$$
\|A({\rm f}_1,\ldots,{\rm f}_k)\|_{L^q(\Q_p)}
\les \|{\rm f}_1\|_{L^{q_1}(\Q_p)}\cdots \|{\rm f}_k\|_{L^{q_k}(\Q_p)}
$$
 for any $q_1,\ldots,q_k\in [D_k,\infty)$  and 
$q\in [q_*,\infty]$ with $1/q_*=1/q_1+\cdots+1/q_k$. Setting ${\rm f}_i= {\ind {\Z_p}}\, f_i$ for $i \in [k]$ and applying the stronger triangle inequality $|x+y| \le \max(|x|,|y|)$, we immediately deduce the conclusion of Theorem \ref{thm:padic} for the range $q \in [q_*, \infty]$. The result for the remaining range $q \in [1, q_*)$ then follows by using H\"{o}lder's inequality.
\end{proof}

It remains to prove Theorem \ref{thm:padiccccc}. We begin by gathering some preliminary results which will be useful in our arguments. 

\subsection*{Change of variables}
We will need a change of variables formula for polynomial maps. 
\begin{prop}\label{CoV-formula}
Let 
$\phi \in {\mathbb Q}_p[\rm x]$ be a polynomial which is injective
on an open set $U \subseteq{\mathbb Q}_p$. Then
\begin{equation}\label{CoV}
\int_U [f\circ \phi] \, |\phi'| \, d\nu \ = \ \int_{\phi(U)} f \, d\nu.
\end{equation}
\end{prop}
For a proof, see \cite{S}.

\subsection*{Hensel's lemma} We will need to use the classical Hensel's lemma which is known to hold
in any non-archimedean local field. See \cite{N} for a proof.

\begin{lemma}\label{classical-hensel}
Let $f \in {\mathbb Z}_p[\rm x]$ be a polynomial. Suppose
there is a $z_0 \in {\mathbb Z}_p$ such that $|f(z_0)| < 1$
and  $|f'(z_0)| = 1$.
Then
there exists a unique $z \in {\mathbb Z}_p$ such that $f(z) = 0$ and $|z-z_0|<1$.
\end{lemma}
\noindent We single out the uniqueness part of Hensel's lemma which will be useful in several places.

{\it Let $P \in {\mathbb Z}_p[\rm x]$ and
suppose that $z,w \in {\mathbb Z}_p$ satisfies the following three properties:
\begin{equation}\label{hyp}
(a) \ P(z) = P(w), \ \ \ (b) \ |z-w| < 1 \ \ \ {\rm and} \ \ (c) \ |P'(w)| = 1.
\end{equation}
Then $z = w$.}

\subsection*{Hardy-Littlewood-Sobolev fractional integration} In the final step of the proof of Theorem \ref{thm:padiccccc}, we will need
to apply the Hardy-Littlewood-Sobolev inequality which holds in 
any local field. See \cite{T} for a proof of the following proposition.

\begin{prop}\label{prop-3} For $0<\beta<1$, set
$$
I_{\beta}(f) \ = \ \int_{{\mathbb Z}_p} \frac{f(x - s)}{|s|^{\beta}} \, d\nu(s).
$$
Then
\begin{equation}\label{p>1}
\|I_{\beta}(f) \|_{L^q({{\mathbb Q}_p})} \ \lesssim_{p,q} \, \|f\|_{L^p({{\mathbb Q}_p})} \ \ {\rm whenever} \ 1<p<q<\infty \ {\rm satisfy} \ \frac{1}{p} - \frac{1}{q} = 1 - \beta.
\end{equation}
Moreover,  $I_{\beta}$ maps $L^1$ to $L^{1/\beta,\infty}$; that is, $I_{\beta}$ is weak-type $(1,1/\beta)$.
\end{prop}

\subsection{Two  propositions on $\Q_p[{\rm x}]$}
\label{section:twoprops}
The proof of Theorem \ref{thm:padiccccc} depends on the following two key propositions.

\begin{prop}\label{prop-1} 
Let $P \in {\Q_p}[{\rm x}]$ be a polynomial with degree $d = {\rm deg}(P) \ge 1$.  
Then there exist 
${\mathcal F}  \subset \Z_p$ (with $\# {\mathcal F} \le C_P$)  and pairwise disjoint open sets
$\{U_j\}_{j\in  [N_P]}$ such that
\begin{itemize}
\item[(i)] $\Z_p \setminus {\mathcal F} \  = \ \bigcup_{j\in [N_P]} U_j$; and
\vskip 7pt
\item[(ii)] $P$ is injective on each $U_j$.
\end{itemize}
The number $N_P$ of open sets in the above partition of $\Z_p \setminus {\mathcal F}$ is finite and depends only on $P$.
\end{prop}

The proof of Proposition \ref{prop-1} is given in Appendix \ref{Appendix:prop-1}.

\begin{remark}\label{Remark:field1}
(1) The proof of Proposition \ref{prop-1} for the real field $\R$ in place of $\Q_p$ is easier.  First,  partition  $\R$ into intervals whose endpoints comprise  the real zeros of $P'$. Then by the mean value theorem, $P$ is injective on each 
interval $I$. However, this proof strongly depends on the ordered structure of ${\mathbb R}$. We direct the reader to \cite{crw}, which gives a proof employing B\'ezout's theorem for the complex field ${\mathbb C}$. \\
(2) Proposition \ref{prop-1} extends to any local field ${\mathbb K}$ of characteristic zero. 
In the positive characteristic case, the proof (see Appendix \ref{Appendix:prop-1})  shows that the result holds if $p>d$ where $p$ is the characteristic of ${\mathbb K}$.
\end{remark}

Next, we give a precise structural description  for sublevel sets of polynomials. 

\begin{prop}\label{prop-2} 
Let
${\mathcal P} = \{P_1, \ldots, P_m\} \subset {\mathbb Q}_p[{\rm x}]$. For each $1\le j\le m$, write
$$
P_j'(y) \ = \ c_j \prod_{\ell} (y - \xi^j_{\ell})^{e_{\ell}^j}, \ \ \ {\rm where} \ \ \sum_{\ell} e_{\ell}^j \ = \ {\rm deg}\, P_j - 1 \ {\rm with}\   e_{\ell}^j\in\Z_+.
$$
Here  the roots $\{\xi^j_{\ell}\}$ lie in the algebraic closure ${\overline{{\mathbb Q}_p}}^{alg}$ of ${{\mathbb Q}_p}$,  and the absolute value
$|\cdot|$ on ${{\mathbb Q}_p}$ extends uniquely to ${\overline{{\mathbb Q}_p}}^{alg}$ which we continue to denote by $|\cdot|$. In what follows, we write 
$$
{\mathcal B}_r(\xi) := \{x \in \Z_p: |x - \xi| \le r\}  =  B^{alg}_r(\xi) \cap \Z_p  \ {\rm where} \  
B^{alg}_r(\xi)  =  \{z\in {\overline{{\mathbb Q}_p}}^{alg}: |z - \xi| \le r \},
$$ 
even if 
$\xi \in {\overline{{\mathbb Q}_p}}^{alg} \setminus {{\mathbb Q}_p}$. 
Then there exist $\delta = \delta_{\mathcal P}>0$
and $C = C_{\mathcal P}>0$ such that for each $1\le j \le m$,
the following statements  hold.
\begin{itemize}
\item[(a)] The balls $\{{\mathcal B}_{C\delta^{1/e_{\ell}^j}}(\xi^j_{\ell}) \}_{\ell}$ are pairwise disjoint and
$$
\bigl\{ y \in \Z_p : |P_j'(y)| \le \delta \bigr\} \ = \ \bigcup_{\ell} {\mathcal B}_{C\delta^{1/e_{\ell}^j}}(\xi^j_{\ell}).
$$
\vskip 7pt
\item[(b)] If $\xi^j_{\ell} \notin \Z_p$, then ${\mathcal B}_{C\delta^{1/e_{\ell}^j}}(\xi^j_{\ell}) = \emptyset$. 
\vskip 5pt
\item[(c)] If $j' \not= j$ and $\xi^{j'}_{\ell'} \not= \xi^j_{\ell}$, then
$$
{\mathcal B}_{C\delta^{1/e_{\ell}^j}}(\xi^j_{\ell}) \cap 
{\mathcal B}_{C\delta^{1/e_{\ell'}^{j'}}}(\xi^{j'}_{\ell'}) = \emptyset.
$$
\end{itemize}

\end{prop}

Proposition \ref{prop-2} follows from a result in \cite{KW}, which is inspired by a sublevel set bound by Phong-Stein-Sturm \cite{PSS} in the real field.

\subsection{Proof of Theorem \ref{thm:padiccccc}}
\label{subsection:proofthm7.1}
 The proof of Theorem \ref{thm:padiccccc} is divided into two parts, presented in Subsections \ref{subsubsection:prop7.1} and \ref{subsubsection2:prop7.1}. 

 \subsubsection{ Crucial decomposition using Propositions \ref{prop-1}
 and \ref{prop-2}}
\label{subsubsection:prop7.1}
We may assume that  the coefficients of  $P_1,\ldots,P_m$ lie  in $\Z_p$. In fact, we can find
a $c \in \Z_p$ with sufficiently small $|c|>0$ such that 
${\tilde{P}}_j :=  c P_j \in \Z_p[{\rm x}]$ for all $1\le j \le m$, 
and write  
$$
A(f_1, \ldots, f_m)(x)  =  \int_{\Z_p} f_1\big(c^{-1}( c x - {\tilde{P}}_1(y)\big) \cdots
f_m\big(c^{-1}( c x - {\tilde{P}}_m(y)\big) \, d\nu(y);
$$
then we can  write $A(f_1, \ldots, f_m)$ as
$$
A(f_1, \ldots, f_m)(x)  \ = \ {\tilde{A}}(g_1, \ldots, g_m)(c x)
$$
where $g_j(x) = f_j(c^{-1} x)$ and ${\tilde{A}}$ is the multilinear averaging operator associated to 
${\tilde{{\mathcal P}}} = \{{\tilde{P}}_1, \ldots, {\tilde{P}}_m\}$. Furthermore, 
by \eqref{CoV}, we have 
$$
\|A(f_1, \ldots, f_m)\|_{L^q}  =  |c|^{-1/q} \|{\tilde{A}}(g_1, \ldots, g_m)\|_{L^q}  \ \ {\rm and} \ \ 
\|g_j\|_{L^{p_j}}  = |c|^{1/p_j} \|f_j\|_{L^{p_j}},
$$
which yield that 
  the bound \eqref{bound} for ${\tilde{A}}$ in place of $A$ implies the bound \eqref{bound}. Therefore, without loss
of generality, we may assume $P_j \in \Z_p[{\rm x}]$ for all $1\le j \le m$.

Fix $1\le j \le m$. Proposition \ref{prop-1} gives  an $O_{\mathcal P}(1)$ collection  of pairwise disjoint
open sets $\{U_{j,k}\}_{k}$ covering $\Z_p$ up to a set of measure zero,
where  $P_j$ is injective on each $U_{j,k}$. Now invoke  part (a) in   Proposition \ref{prop-2} to partition each open set $U_{j,k}$ as
\begin{equation}\label{open}
U_{j,k} \ = \ \Bigl( U_{j,k} \cap \{|P'_j| > \delta\} \Bigr)  \ \cup \Big( \bigcup_{\ell} \ \bigl(U_{j,k} \cap 
{\mathcal B}_{C\delta^{1/e_{\ell}^j}}(\xi^j_{\ell}) \bigr)\Big),
\end{equation}
where $\{|P'_j| > \delta\}:=\{ y \in \Z_p : |P_j'(y)| > \delta \}$.  
Let us denote by $\{E^j_n \}_{1\le n \le M_j}$ the total collection of open sets in \eqref{open} 
such that on each $E^j_n$,
$P_j$ is injective and either $|P'_j(y)| > \delta$ or $|y - \xi^j_{\ell}| \le C \delta^{1/e_{\ell}^j}$ for some $\ell$. 
Note that  $M_j = O_{\mathcal P}(1)$ for all $1\le j \le m$. We decompose
$$
A(f_1, \ldots, f_m)(x) =  \sum_{{\underline{n}}\in {\mathcal M}} \int_{E^1_{n_1} \cap \cdots \cap E^m_{n_m}} f_1\big(x - P_1(y)\big)  \cdots f_m\big(x - P_m(y)\big) \, d\nu(y).
$$
where ${\mathcal M} = \{ {\underline{n}} = (n_1, \ldots, n_m): 1 \le n_j \le M_j \}$. Write the integral on the right-hand side  as
$A_{\underline{n}}(f_1, \ldots, f_m)(x)$ so that  
$$A(f_1, \ldots, f_m)(x) = \sum_{{\underline{n}}\in {\mathcal M}} A_{\underline{n}}(f_1,\ldots, f_m)(x).$$

Fix ${\underline{n}}=(n_1, \ldots, n_m)\in {\mathcal M}$, set
${\mathcal E}_{\underline{n}}:=  E^1_{n_1} \cap \cdots \cap E^m_{n_m}$ and let
$$
{\mathcal L}_{\underline n}:= \bigl\{ 1 \le j \le m :  \ E^j_{n_j} = U_{j,k} \cap \{ |P'_j| > \delta\} \ \ {\rm for \ some} \ k \bigr\}.
$$
For the case  ${\mathcal L}_{\underline n}=\{1,2,\ldots,m\}$,  the estimate $\|A\|_{L^{p_1}\times \cdots \times L^{p_m}\to L^{\infty}}\les 1$ follows  from  H\"{o}lder's inequality and \eqref{CoV}. Interpolating   this with the trivial bound $\|A\|_{L^{p_1}\times \cdots \times L^{p_m}\to L^{p_*}}\le 1$ gives  the claimed result. We thus focus subsequently on ${\mathcal L}_{\underline n}\subsetneqq \{1,2,\ldots,m\}$. 
For each $j\notin {\mathcal L}_{\underline n}$, 
$E^j_{n_j} = U_{j,k} \cap {\mathcal B}_{C\delta^{1/e^j_{\ell}}}(\xi^j_{\ell})$ for some $k$ and $\ell$. The conclusion from part (c) in 
Proposition \ref{prop-2} gives that,  if $\xi^j_ {\ell}\not= \xi^{j'}_{\ell'}$ for some $j, j' \notin {\mathcal L}_{\underline n}$, then $E^j_{n_j} \cap E^{j'}_{n_{j'}}
= \emptyset$ and so ${\mathcal E}_{\underline n} = \emptyset$. Therefore,  if ${\mathcal E}_{\underline n} \not= \emptyset$ and ${\mathcal L}_{\underline n}\subsetneqq \{1,2,\ldots,m\}$,  we write 
$$
{\mathcal E}_{\underline{n}}  =\,    {\mathcal U}_{\underline{n}} \ \cap \ \Bigl(\bigcap_{j\in {\mathcal L}_{\underline{n}}} \ \bigl\{ |P'_j| > \delta \bigr\}\Bigr) \
\cap \  {\mathcal B}_{C\delta^{\eta}}(\xi)
$$
for some $\eta>0$, 
where ${\mathcal U}_{\underline{n}}$ is an open set on which every $P_j$ is injective and $\xi$ is a common root of $P'_j$ for all $j\notin {\mathcal L}_{\underline n}$.  
Since the conclusion from part (b) in Proposition \ref{prop-2} implies  ${\mathcal B}_{C\delta^{\eta}}(\xi) = \emptyset$ whenever  $\xi \notin \Z_p$,     we can assume 
$\xi \in \Z_p$ in what follows.

Now,  making the change of variables $y \to y + \xi$, we have
$$
A_{\underline{n}}(f_1,\ldots, f_m)(x) \ = \ \int_{{\mathcal F}_{\underline{n}}}  g_1\big(x - Q_1(y)\big) \cdots g_m\big(x - Q_m(y)\big) \, 
d\nu(y)
$$ 
where ${\mathcal F}_{\underline{n}}:={\mathcal E}_{\underline{n}}-\xi$, and  for each $1\le j \le m$, $g_j(x)= f_j(x - P_j(\xi))$ and $Q_j(y) = P_j(y + \xi) - P_j(\xi) =y^{e_j} R_j(y)$ with polynomial $R_j$,  $1\le e_j \le d_j$ and 
$|Q'_j(y)| \sim_{{\mathcal P}} |y|^{e_j -1}$ on ${\mathcal B}_{C\delta^{\eta}}(0)$.  Furthermore,
$$
{\mathcal F}_{\underline{n}} \ = \ \bigl({\mathcal U}_{\underline{n}} - \xi \bigr) \ \cap \ \Bigl(\bigcap_{j\in {\mathcal L}_{\underline{n}}} \
\bigl\{ |Q'_j| > \delta \bigr\}\Bigr) \ \cap \  {\mathcal B}_{C\delta^{\eta}}(0).
$$
Note that  $\|f_j\|_{L^r} = \|g_j\|_{L^r}$ for all $1\le j\le m$,  and each $Q_j$ is injective on ${\mathcal U}_{\underline{n}} - \xi$.
\subsubsection{Main arguments}\label{subsubsection2:prop7.1}
Without loss of generality, we may assume that each $f_j$ (and therefore $g_j$) is nonnegative.  Our approach involves an iterative process of changes of variables: first setting $s_1 = Q_1(y)$,  and then $s_2 = Q_2(y)$,  
 and so forth, while systematically interlacing applications of 
H\"older's inequality. We use   \eqref{CoV} to conclude
$$
A_{\underline{n}}(f_1, \ldots, f_m)(x) = \int_{Q_1({\mathcal F}_{\underline{n}})} 
\frac{g_1( x- s) g_2(x - Q_2(\phi(s))) \cdots g_m(x - Q_m(\phi(s)))}{|Q_1'(\phi(s))|} \, d\nu(s)
$$
where $\phi(s) = Q_1^{-1}(s)$ is the inverse map of $Q_1$ on ${\mathcal F}_{\underline{n}}$. 
Note that $|s| = |Q_1(y)| \sim_{\mathcal P}  |y|^{e_1}$ and so $|y| = |\phi(s)| \sim |s|^{1/e_1}$.

\vskip.1in

Applying H\"older's inequality and then changing variables back again, recalling that $|Q'_1(y)| \sim_{{\mathcal P}}   |y|^{e_1 - 1}$ on ${\mathcal F}_{\underline{n}}$, we  have
$$
\begin{aligned}
A_{\underline{n}}(x)  \les_{{\mathcal P}}&  \   \|f_1\|_{L^{p_1}} \Bigl( \int_{Q_1({\mathcal F}_{\underline{n}})} 
\frac{g_2^{p'_1}(x - Q_2(\phi(s))) \cdots g_m^{p'_1}(x - Q_m(\phi(s)))}{|\phi(s)|^{(e_1 - 1)p'_1}} \, d\nu(s) \Bigr)^{1/p'_1}\\
\sim_{{\mathcal P}}& \  \|f_1\|_{L^{p_1}} \, \Bigl( \int_{{\mathcal F}_{\underline{n}}} 
\frac{g_2^{p'_1}(x - Q_2(y)) \cdots g_m^{p'_1}(x - Q_m(y))}{|y|^{(e_1 - 1)(p'_1 - 1)}} \, d\nu(y) \Bigr)^{1/p'_1},
\end{aligned}
$$
where we have used the abbreviation $A_{\underline{n}}(x) := A_{\underline{n}}(f_1, \ldots, f_m)(x)$. 
Next, changing variables $s = Q_2(y)$ so that $|s| = |Q_2(y)| \sim_{\mathcal P}  |y|^{e_2}$
(with  $y = \psi(s) = Q_2^{-1}(s)$ on ${\mathcal F}_{\underline{n}}$), we have
$$
\begin{aligned}
A_{\underline{n}}(x) \ \lesssim_{\mathcal P}& \ \|f_1\|_{L^{p_1}} \ \Bigl( \int_{Q_2({\mathcal F}_{\underline{n}})} 
\frac{g_2^{p'_1}(x - s) \cdots g_m^{p'_1}(x - Q_m(\psi(s)))}{|\psi(s)|^{(e_1 - 1)(p'_1-1)} |Q'_2(\psi(s))|} \, d\nu(s) \Bigr)^{1/p'_1} \\
\sim_{\mathcal P}&  \ \|f_1\|_{L^{p_1}} \ \Bigl( \int_{Q_2({\mathcal F}_{\underline{n}})} 
\frac{g_2^{p'_1}(x - s) g_3^{p'_1}(x - Q_3(\psi(s)))\cdots g_m^{p'_1}(x - Q_m(\psi(s)))}{|s|^{\alpha_2}} \, d\nu(s) \Bigr)^{1/p'_1},
\end{aligned}
$$
where $\alpha_2 $ is given by 
$$
\alpha_2  = \frac{e_2-1}{e_2} + \frac{e_1 - 1}{e_2}\, (p'_1-1)  =  \frac{e_2-e_1}{e_2} + 
\frac{e_1 - 1}{e_2}\, p'_1.
$$

Now,  iteratively apply H\"older's inequality and the change of variables $s = Q_j(y)$ ($j=3,\ldots,m$) to arrive at
\begin{equation}\label{iteration-final}
A(f_1, \ldots, f_m)(x)  \lesssim_{\mathcal P} \ \Big(\prod_{i\in [m-1]}\big\|f_i\big\|_{L^{p_i}} \Big)\,
\Bigl( \int_{Q_m({\mathcal F}_{\underline{n}})} \frac{g_m^{u_1' \cdots u_{m-1}'}(x-s)}{|s|^{\alpha_m}} \, d\nu(s) 
\Bigr)^{\frac{1}{u_1' \cdots u_{m-1}'}}, 
\end{equation}
where $u_1 = p_1$,   $u_j = p_j/(u_1' \cdots u_{j-1}')$ for $j=2,\ldots,m$, and 
$$
\alpha_m := \frac{e_m-e_{m-1}}{e_m} + \frac{e_{m-1} - e_{m-2}}{e_{m}}\, u_{m-1}' +  \cdots  +   \frac{e_1 -1}{e_m}\,
u_1' \cdots u_{m-1}',\quad m\ge2.
$$
Here  $u_j'$ is the H\"older conjugate of $u_j$ for each $j\ge 1$. 
We claim  
\begin{equation}\label{u-formula}
\frac{1}{u'_1 \cdots u'_j} =  1-\big( \frac{1}{p_1} + \frac{1}{p_2} + \cdots + \frac{1}{p_j}\big).
\end{equation}
This is clearly true for $j=1$.  Since $u_j = p_j/(u_1' \cdots u_{j-1}')$ for $j\ge 2$,  we have 
$$
\frac{1}{u'_1 \cdots u'_j}  =  \frac{1}{u'_1 \cdots u'_{j-1}} \, \Bigl(1 - \frac{1}{u_j}\Bigr)  =  \frac{1}{u'_1 \cdots u'_{j-1}}  -   \frac{1}{p_j},
$$
which yields \eqref{u-formula}  by induction. Furthermore,  
since    $p_{*}\ge 1$ and $1<p_1,p_m<\infty$, 
 we obtain from  \eqref{u-formula} that 
$$
  1 - \frac{1}{p_1} - \cdots - \frac{1}{p_{j-1}}  > \frac{1}{p_j} \qquad 
  {\rm and }\qquad u_j >1
$$
for all $1\le j \le m-1$.  This justifies the application of H\"older's inequality at each iterative step,
leading to \eqref{iteration-final}.
Using  \eqref{u-formula}, a straightforward computation establishes the following lemma.

\begin{lemma}\label{e-formula}
For each $m\ge 2$, 
 we have
$$
\alpha_m  = 1  +  \frac{u'_1 \cdots u'_{m-1}}{e_m} \, \Bigl( \,\sum_{j=1}^{m-1} \frac{e_j}{p_j} - 1 \, \Bigr).
$$
\end{lemma}

By our hypothesis \eqref{condition} and $1\le e_j \le d_j$ for all $1\le j \le m$,  we have 
$$
\Bigl( \, \sum_{j=1}^{m-1} \frac{e_j}{p_j} - 1 \, \Bigr) \ \le \ \Bigl( \, \sum_{j=1}^{m-1} \frac{d_j}{p_j}  - 1 \, \Bigr) \ < \ 0,
$$
which, with Lemma \ref{e-formula}, shows
$$
\alpha_m = 1 + \frac{u'_1 \cdots u'_{m-1}}{e_m} \, \Bigl( \, \sum_{j=1}^{m-1} \frac{e_j}{p_j} - 1 \, \Bigr) 
\ \le  1 + \frac{u'_1 \cdots u'_{m-1}}{d_m} \Bigl( \, \sum_{j=1}^{m-1} \frac{d_j}{p_j}  -  1 \, \Bigr)   =:  \beta_m.
$$
Now,  using  $Q_m({\mathcal F}_{\underline{n}}) \subseteq \Z_p$ 
and $|s|^{-\alpha_m} \le |s|^{-\beta_m}$ ($|s|\le 1$ when $s \in \Z_p$), we obtain from 
\eqref{iteration-final} that 
\begin{equation}\label{iteration-final-again}
A(f_1, \ldots, f_m)(x)  \les_{\mathcal P} \, \Big(\prod_{i\in [m-1]}\big\|f_i\big\|_{L^{p_i}} \Big)\,
\Bigl( \int_{\Z_p} \frac{g_m^{u_1' \cdots u_{m-1}'}(x-s)}{|s|^{\beta_m}} \, d\nu(s) 
\Bigr)^{\frac{1}{u_1' \cdots u_{m-1}'}}.  
\end{equation}
We claim $\beta_m\in (0,1)$. In fact, by applying Lemma \ref{e-formula} with the $e_j$'s replaced by the $d_j$'s, 
\begin{equation}\label{d-formula}
\begin{aligned}
 \beta_m \ =& \ 1  +  \frac{u'_1 \cdots u'_{m-1}}{d_m} \, \Bigl( \, \sum_{j=1}^{m-1} \frac{d_j}{p_j} - 1 \, \Bigr)\\
 =&\ \frac{d_m - d_{m-1}}{d_m} + \frac{d_{m-1} - d_{m-2}}{d_{m}}\, u_{m-1}' +  \cdots  +   \frac{d_1 -1}{d_m}\,
u_1' \cdots u_{m-1}',
\end{aligned}
\end{equation}
which gives  $\beta_m > 0$ (since otherwise $1 = d_1 = d_2 = \cdots = d_m$). Moreover, from  \eqref{d-formula} and our hypothesis \eqref{condition}  we see that  
 $\beta_m < 1$.

 Hence,  we have
\begin{equation}\label{q-est}
\|A(f_1, \ldots, f_m)\|_{L^q}  \les_{\mathcal P} \, \|f_1\|_{L^{p_1}} \cdots \|f_{m-1}\|_{L^{p_{m-1}}} \,
{\|I_{\beta_m}(h_m)\|_{L^{{\tilde q}}}^{1/(u'_1 \cdots u'_{m-1})}},
\end{equation}
where $\beta_m\in (0,1)$,  ${\tilde q} = q/(u'_1 \cdots u'_{m-1})$, $h_m(x) = g_m^{u_1' \cdots u_{m-1}'}(x)$ and $I_{\beta_m}(f)$ is  the  $p$-adic  field  version  of  fractional  integration of order $\beta_m$ given by 
$$
I_{\beta_m}(f)(x)  =  \int_{\Z_p} \frac{f(x-s)}{|s|^{\beta_m}} \, d\nu(s).
$$

We now apply the $p$-adic field  Hardy-Littlewood-Sobolev inequality to the last factor in either \eqref{iteration-final-again} or \eqref{q-est}, with the Lebesgue exponents 
$${\tilde p} = p_m/(u'_1 \cdots u'_{m-1})\qquad {\rm  and}\qquad  
{\tilde q} = q/(u'_1 \cdots u'_{m-1}),$$
where $q$ is  as in \eqref{bound}. 
The subsequent argument is divided into four cases.

\medskip \paragraph{\bf Case  $p_*=1$}  
If $1/p_*=1/p_1 + \cdots + 1/p_m = 1$, then  $K>1$,   $u'_1 \cdots u'_{m-1} = p_m$ (so we have  
${\tilde p} = 1$) and $\beta_m =p_m(K-1)/d_m= p_m/q_{*,K}$. Hence, by interpolating the
endpoint (weak-type) bound  $\|I_{\beta_m}\|_{ L^1 \to L^{{{1}/{\beta_m}},\infty}}\les 1$ with the trivial bound
$\|I_{\beta_m}\|_{L^1 \to L^1}\les 1$, we conclude that $\|I_{\beta_m}\|_{L^1 \to L^{\tilde q}}\les 1 $ for all $1\le \tilde{q} < 1/\beta_m$. This gives the 
bound \eqref{bound} for $1=p_*\le q < q_{*,K}$.

\medskip \paragraph{\bf Case  $p_*>1$ and $K>1$}  
Suppose that $1/p_*=1/p_1 + \cdots + 1/p_m < 1$ (so we have $u'_1 \cdots u'_{m-1} < p_m$, and then   ${\tilde p}>1$) and $K = d_1/p_1 + \cdots + d_m/p_m > 1$. 
By our hypothesis \eqref{condition}, 
$$
d_m \Bigl(\frac{1}{q_{*,K}} - \frac{1}{p_m} \Bigr) \ = \ \sum_{j=1}^{m-1} \frac{d_j}{p_j} -  1  <  0,
$$
which shows $q_{*,K} > p_m$ and 
$
 1 < {\tilde p} <  {{ q}_{*,K}}/(u'_1 \cdots u'_{m-1})=:{\tilde q}_{*,K}.
$
In addition, a simple computation gives  
$$
\begin{aligned}
 \beta_m  -  1  =& \ \frac{u'_1 \cdots u'_{m-1}}{d_m} \, \Bigl( \, \sum_{j=1}^{m-1} \frac{d_j}{p_j} -  1 \, \Bigr)  =\ 
\frac{u'_1 \cdots u'_{m-1}}{d_m} \, d_m \, \bigl(\frac{1}{q_{*,K}} - \frac{1}{p_m} \bigr) 
=\  \frac{1}{{\tilde q}_{*,K}}  -  \frac{1}{\tilde p}.
\end{aligned}
$$
By using the inequality 
\eqref{p>1}, we have 
$\|A\|_{L^{p_1}\times \cdots \times L^{p_m}\to L^{q_{*,K}}}\les_{\mathcal P} 1$. Interpolating 
this with  the trivial bound $\|A\|_{L^{p_1}\times \cdots \times L^{p_m}\to L^{p_*}}\le 1$ gives  \eqref{bound}
for all $p_*\le q \le q_{*,K}$.

\medskip \paragraph{\bf Case  $K=1$} 
In this case, we have 
$$
 1  =  \frac{d_1}{p_1} +  \frac{d_2}{p_2} +  \cdots  +  \frac{d_m}{p_m}
$$
and so $p_*>1$.
From \eqref{d-formula}, we have
$$
1 - \beta_m \ = \frac{u'_1 \cdots u'_{m-1}}{d_m} \, \Bigl( 1 -  \sum_{j=1}^{m-1} \frac{d_j}{p_j} \Bigr) \ = \ 
\frac{u'_1 \cdots u'_{m-1}}{p_m} \ = \ \frac{1}{\tilde p}.
$$
Now fix any $q\in ( p_m,\infty)$ and define $0< \beta < 1$ by
$$
\beta  := 1 + \frac{1}{\tilde q} -  \frac{1}{\tilde p}  \ge  \beta_m
$$
where ${\tilde q} = q/(u'_1 \cdots u'_{m-1})$. Since $|s|^{\beta} \le |s|^{\beta_m}$ for all $s \in \Z_p$, we see from  \eqref{iteration-final-again}  that
$$A(f_1, \ldots, f_m)(x) \les_{\mathcal P} \, \Big(\prod_{i\in [m-1]}\|f_i\|_{L^{p_i}}\Big)\,
\Bigl( \int_{\Z_p} \frac{g_m^{u_1' \cdots u_{m-1}'}(x-s)}{|s|^{\beta}} \, d\nu(s) 
\Bigr)^{\frac{1}{u_1' \cdots u_{m-1}'}}, 
$$
which, 
with \eqref{p>1} in Proposition \ref{prop-3} and the trivial bound $\|A\|_{L^{p_1}\times \cdots \times L^{p_m}\to L^{p_*}}\le 1$, implies 
the bound  \eqref{bound}
for all $p_*\le q<q_{*,K}=\infty$.

\medskip \paragraph{\bf Case  $K<1$}
In this case, we also have $p_*>1$ and
\beq\label{smallest1}
 \frac{d_m}{p_m} \ < \ 1 - \sum_{j=1}^{m-1} \frac{d_j}{p_j}.
\eeq
From \eqref{d-formula}, we have
$$
1 - \beta_m \ = \frac{u'_1 \cdots u'_{m-1}}{d_m} \, \Bigl( 1 -  \sum_{j=1}^{m-1} \frac{d_j}{p_j} \Bigr).
$$
This, with \eqref{smallest1} and   ${\tilde p} = p_m/(u'_1 \cdots u'_{m-1})$,   shows 
\begin{equation}\label{integrability}
\frac{1}{\tilde p} < 1 - \beta_m; \qquad{\rm that\  \ is,}\  \quad 
\beta_m \ {\tilde p}' < 1.
\end{equation}
Thus, an application of Hölder's inequality in \eqref{iteration-final-again} (with conjugate exponents $\tilde p$ and $\tilde p'$) yields the estimate
\beq\label{EDN1Q}
\|A(f_1, \ldots, f_m)\|_{L^\infty} \ \les_{\mathcal P} \, \big(\prod_{i\in [m-1]}\|f_i\|_{L^{p_i}}\big)\,  \|f_m\|_{L^{p_m}}  \,
\Bigl(\int_{\Z_p} \frac{1}{|s|^{\beta_m {\tilde p}'}} d\nu(s) \Bigr)^{\frac{1}{u'_1 \cdots u'_{m-1}{\tilde p}'}}.
\eeq
Since $\|A\|_{L^{p_1}\times \cdots \times L^{p_m}\to L^{p_*}}\le 1$ and the  integral on the right-hand side of \eqref{EDN1Q} is finite by \eqref{integrability}, as a consequence, interpolation gives that   \eqref{bound} holds for all $p_*\le q \le q_{*,K}=\infty$.

This completes the proof of Theorem \ref{thm:padiccccc}.

\appendix
\section{Multilinear Weyl inequality in the continuous setting}
\label{section:Appen1}
In this section, we will show the multilinear Weyl inequality in the continuous setting.
For $k\in\Z_+$, 
we define a multilinear operator  in the continuous setting 
\beq\label{Appdix1}
\begin{aligned}
\tilde A_{N; \mathbb R}^{\mathcal P}(f_1,\ldots, f_k)(x):=&\ 
\frac{1}{N} \int_{{N}/{2}}^N f_1(x-P_1(t))\cdots f_k(x-P_k(t)) dt,\qquad x\in \R,
\end{aligned}
\eeq
where $\mathcal P:= (P_1,\ldots, P_k)$ be a polynomial mapping with real coefficients and distinct degrees such that   (\ref{eq:41}) holds. 
\begin{thm}\label{appen11}
Let  $N \geq 1$, $l\in \N$, $k\in\Z_+$,    and let $\mathcal P$ be given as in (\ref{Appdix1}).  Let
$1<q_1,\ldots, q_k<\infty$ be exponents such that
$\frac{1}{q_1}+\cdots+\frac{1}{q_k}=\frac{1}{q}\le 1$. 
There exist two positive constants   $c$ and $C$,  depending on
$\mathcal P, q_1, \dots, q_k$, such that the following holds. Let 
$f_i \in L^{q_i}(\R)$ for each $i\in [k]$, and fix $j\in [k]$.
If $f_j \in L^{q_j}(\R) \cap L^{2}(\R)$  and
 $\F_{\R} f_j$ vanishes on $[-N^{-d_j}2^{d_j l},N^{-d_j}2^{d_j l}]$,
then
\begin{align}
\label{appuweyl11777}
\| \tilde A^{\mathcal{P}}_{N; \R}(f_1,\ldots, f_k) \|_{L^q(\R)} 
\le C (2^{-cl}+N^{-c})\|f_1\|_{L^{q_1}(\R)}\cdots \|f_k\|_{L^{q_k}(\R)}.
\end{align}
\end{thm}
\begin{remark}\label{remaAppes}
The  interval $[N/2,N]$ in \eqref{Appdix1} 
 can be replaced by $[c_1N,c_2N]$ for any $c_2>c_1\ge 0$. Moreover,  similar conclusion also holds for  the general  operator $T_{N; \mathbb R}^{\mathcal P}$ given by 
$$
\begin{aligned}
 T_{N; \mathbb R}^{\mathcal P}(f_1,\ldots, f_k)(x):=&\ 
 \int_{{1}/{2}}^1 f_1(x-P_1(Nt))\cdots f_k(x-P_k(Nt)) \rho_N(t) dt,
\end{aligned}
$$
where the function $\rho_N$ satisfies the bound  $|\rho_N'(t)|+|\rho_N(t)|\les 1$ for all $1/2\le t\le 1$. 
\end{remark}
\begin{proof}
Note that (\ref{appuweyl11777}) without the decay   $2^{-cl}+N^{-c}$ follows from 
  Minkowski's inequality and  H\"{o}lder's inequality.  
By interpolation, it suffices to show \eqref{appuweyl11777} for the case $q=1$.
Since  the proof can be achieved by following  the arguments in \cite[Section 7]{KMPWfield}, we omit the details.
\end{proof}
\section{Sampling map and quantitative Shannon sampling theorem}
\label{appenmap}
In this section, we introduce several important maps utilized in the  proof of the major arcs estimates. For motivations and detailed relationships, we refer to \cite[Section 4]{KMT22}.
The inclusion homomorphism $\iota: \Z \to \Ad_\Z$ 
is defined by
$ \iota(x):= \left(x, ( (x\mod p^j)_{j \in \N} )_{p \in \mathbb{P}} \right)$
and the addition homomorphism $\pi: \R \times \Q/\Z \to \T$ is given  by
\beq\label{pimap}
 \pi(\theta,\alpha):= \alpha + \theta;
\eeq
these two maps are Fourier adjoint to each other in the sense that
\begin{equation}
\label{adjoint} 
\iota(x) \cdot \xi = x \cdot \pi(\xi)
\end{equation}
for all $x \in \Z$ and $\xi \in \R \times \Q/\Z$.
On the major arcs, the map $\pi$ is injective (we call domains $\Omega \subset \R \times \Q/\Z$ on which $\pi$ is injective {\it non-aliasing}). We use these homomorphisms to ``approximate'' $\Z$ by $\Ad_\Z$, which in principle decouples the discrete harmonic analysis of $\Z$ from the continuous harmonic analysis of $\R$ and the arithmetic harmonic analysis of $\hat \Z$.

The  inclusion homomorphism $\iota$ mentioned above leads to  a sampling map $\Sample:\Schwartz(\Ad_\Z) \to \Schwartz(\Z)$, which is defined as follows:
\beq\label{interpolationop333}
\Sample f(x):= f( \iota(x) )
\eeq
for $x \in \Z$ and $f \in \Schwartz(\Ad_\Z)$.  Dually, the addition homomorphism $\pi : \R \times \Q/\Z \to \T$ leads to a projection map $\mathcal{T} :  \Schwartz(\R \times \Q/\Z) \to \Schwartz(\T)$ given  by 
$$ \mathcal{T} F(\xi) := \sum_{(\theta, \alpha) \in \pi^{-1}(\xi)} F(\theta, \alpha)$$
for $\theta \in \R$, $\alpha \in \Q/\Z$, and $F \in \Schwartz(\R \times \Q/\Z)$. It is important to note that the definition of $\Schwartz(\R \times \Q/\Z)$ ensures that this sum contains at most countably many non-zero terms.  
We obtain  from \eqref{adjoint}  the identity
$\F_\Z^{-1} \circ \mathcal{T} = \Sample \circ \F_{\Ad_\Z}^{-1}$
or equivalently 
$ \F_\Z \circ \Sample = \mathcal{T} \circ \F_{\Ad_\Z}.$

For a set of frequencies $\Omega \subseteq \G^{*}$ in the dual group of a LCA group $\G$, we denote by
$\Schwartz(\G)^\Omega$ the subspace of $\Schwartz(\G)$ consisting of functions that are Fourier supported on $\Omega$.
For a non-aliasing compact set of adelic frequencies $\Omega \subset \R \times \Q/\Z$ (so that $\pi : \Omega \to \T$ is injective), the sampling map $\Sample$ in \eqref{interpolationop333}, when restricted to $\Schwartz(\Ad_\Z)^\Omega$, is an injective map and this allows us to define the {interpolation operator} 
\beq\label{interpolationoer}
\Sample_\Omega^{-1} \colon \Schwartz(\Z)^{\pi(\Omega)} \to \Schwartz(\Ad_\Z)^\Omega.
\eeq
This operator, together with the sampling operator $\Sample$, forms a unitary correspondence between $\ell^2(\Z)^{\pi(\Omega)}$ and $L^2(\Ad_\Z)^\Omega$.
For all these facts about $\Sample$ and ${\mathcal T}$, see \cite{KMT22}.

Finally,  we  state  the quantitative Shannon sampling theorem (see \cite[Theorem 4.18]{KMT22}).
\begin{thm}\label{Sampling} 
Let $0 < q \leq \infty$, and  $B$ be a finite-dimensional normed vector space. If $F \in \Schwartz(\Ad_\Z; B)$ has Fourier support in $[-\frac{c_0}{Q},\frac{c_0}{Q}] \times \frac{1}{Q}\Z/\Z$ for some $Q \in \Z_+$ and some $0 < c_0 < {1}/{2}$, then we have 
\begin{equation}\label{e:SAMP} \| \Sample F \|_{\ell^q(\Z;B)} \sim_{c_0,q} \|F\|_{L^q(\Ad_\Z; B)}
\end{equation}
where we extend the sampling operator $\Sample$ to vector-valued functions in the obvious fashion.
\end{thm}

\section{Proof of Proposition \ref{prop-1}}
\label{Appendix:prop-1}
In this section, we prove Proposition \ref{prop-1}. 
Before we turn to its proof, we establish two useful lemmas giving the bound for the $p$-adic value of  $n!$ and desired partition for ``small balls". We will use the notation $O(1)$ or $O_Q(1)$ to denote a large constant/number that only depends on some polynomial $Q$ but importantly, does not depend on the prime $p$. The value of $O(1)$ or $O_Q(1)$  may change from line to line.

\subsection{Two lemmas} 
For $p$ a prime and $k\in {\mathbb \Z}_+$, let $\nu_p(k)\ge 0$ denote the additive $p$-adic valuation; that is $p^{\nu_p(k)} | \, k$ but $p^{\nu_p(k)+1} \not| \ k$ or in other words, $|n| = p^{-\nu_p(n)}$. Our first lemma
is the following.

\begin{lemma}\label{p-factorial} Let $n\ge 2$. Then $\nu_p(n!) \le n-2$ whenever $p\ge 3$ or  $p=2$ and $n \notin 2^\N$. Furthermore,  if $p=2$ and $n\in 2^\N$, we have $\nu_2(n!) = n-1$.
\end{lemma}

\begin{proof}  We have
$$
\nu_p(n!) \ = \ \sum_{i=1}^{\infty} \bigl\lfloor \frac{n}{p^i} \bigr\rfloor \ < \ n \ \sum_{i=1}^{\infty} \frac{1}{p^i} \ = \ \frac{n}{p-1} \ \le \ n - 1,
$$
where 
the last inequality holds if and only if $p\ge 3$ and $n\ge 2$. 
Hence $\nu_p(n!) \le n-2$.

When $p=2$ and $2^N < n < 2^{N+1}$, we have $\lfloor \log_2 n \rfloor = N$ and $1< n/2^N < 2$, implying $1 = \lfloor n/2^N \rfloor < n/2^N$. Hence
$$
\nu_2(n!) =  \sum_{i=1}^{\infty} \bigl\lfloor \frac{n}{2^i} \bigr\rfloor  = 
\sum_{i=1}^{N} \bigl\lfloor \frac{n}{2^i} \bigr\rfloor \ <
\ n  \sum_{i=1}^{N} 2^{-i}  \ =  \ n( 1 - 2^{-N} )  \  \le  \ n-1,
$$
where  
the last inequality follows since $2^N \le n$. Therefore,
$\nu_2(n!) \le n-2$.

If $n = 2^k$ so that $\log_2 n = k$, we have
$$
\nu_2(n!) \ =  \ \sum_{i=1}^{\infty} \bigl\lfloor \frac{n}{2^i} \bigr\rfloor  \ = \
\sum_{i=1}^{k} 2^{k-i} \ = \ 2^k -1 \ = \ n-1.
$$
\end{proof}

Recall the notation ${\mathbb N}_{\le p-1} = \{0,1,\ldots, p-1\}$ from Subsection \ref{subsectionbasicnotation}, and 
keep in mind that the set
${\mathbb N}_{\le p-1} $ is a complete set of representatives
for the residue field  ${\mathbb F}_p:=\Z_p/p\Z_p$.
For any $a \in {\mathbb Z}_p$, we have $a = \sum_{j\ge 0} a_j p^j$ where each $a_j \in {\mathbb N}_{\le p-1}$.
If $b = \sum_{j\ge 0} b_j p^j$ is another element of ${\mathbb Z}_p$, we have $|a-b|<1$ if $a_0 = b_0$. Also we denote ${\overline{a}} = [a_0] \in {\mathbb F}_p$ as the residue field element given by $a_0$; that is, ${\overline{a}}=a_0+p\Z_p$.

The second lemma applies to general polynomials
$Q \in {\mathbb Z}_p [{\rm x}]$ and examines the injectivity property of $Q$ on a ball $B_{p^{-1}}(y_{*})$
where $y_{*} \in {\mathbb Z}_p$.
More precisely,  we decompose $B_{p^{-1}}(y_{*})\setminus \{y_{*}\}$ into $O_Q(1)$  open sets and 
reduce matters to proving injectivity on each of these sets. First,
we partition the residue field ${\mathbb F}_p = U_1 \cup \cdots \cup U_L$ with $L = O_Q(1)$ so that
\begin{equation}\label{1-1-again}
\forall  1\le \ell\le L, \ \ {\rm if} \ \bar{e}, \bar{f} \in U_{\ell} \ \ {\rm and} \ \ \bar{e}^j = \bar{f}^j \ 
{\rm for \ some} \ 1\le j \le {\rm deg}(Q), \ {\rm then} \ \bar{e} = \bar{f}.
\end{equation}
(In other words, the conclusion is  $e_0=f_0$ or $|e-f|<1$  when $\bar{e}=[e_0]$ and $\bar{f}=[f_0]$.) Indeed,
we can achieve this since for any $w \in {\mathbb F}_p$, we have $\# \{z \in {\mathbb F}_p : z^j = w\} \le j$.
We now write
$$
B_{p^{-1}}(y_{*})\setminus \{y_{*}\} = \bigcup_{j=1}^L {\mathcal U}_j \ \ 
{\rm where} \ \ \ {\mathcal U}_j  =  \bigl\{ y_{*} + p^{\ell} u: {\rm for \ some} \ \ell \ge 1, |u| =1 \ {\rm and} \ {\overline{u}} \in U_j \bigr\}.
$$
Next, for a large positive integer $\ell_{*} = O_Q(1)$ (to be determined later), we further decompose each ${\mathcal U}_j$ as 
$$
{\mathcal U}_j \ = \ {\mathcal U}_j^{\ge \ell_{*}} \ \cup \ \Bigl[ \bigcup_{\ell =1}^{\ell_{*}-1} {\mathcal U}_j^{\ell} \Bigr], \ \ {\rm where} \ \  {\mathcal U}_j^{\ell} = \bigl\{ y_{*} + p^{\ell} u:  |u| =1 \ {\rm and} \ {\overline{u}} \in U_j \bigr\},
$$
where 
$$
{\mathcal U}_j^{\ge \ell_{*}} = \bigl\{ y_{*} + p^{\ell} u: {\rm for \ some} \ \ell \ge \ell_{*}, |u| =1 \ {\rm and} \ {\overline{u}} \in U_j \bigr\} .
 $$
\begin{lemma}\label{poly-basic-decomp} With the set-up as above, suppose the polynomial $Q \in {\mathbb Z}_p[{\rm x}]$ satisfies $|Q^{(k)}/k!| \equiv 1$ on $B_{p^{-1}}(y_{*})$ for some $k\ge 1$. For each $1\le j \le L$,  the following 
holds. 
If $p\ge 3$, 
 $Q$ is injective on ${\mathcal U}_j^{\ge \ell_{*}}$; and if $p=2$, 
we can partition  ${\mathcal U}_j^{\ge \ell_{*}}$ into two open subsets such that 
$Q$ is injective on each open subset. Hence the injectivity of $Q$ on $B_{p^{-1}}(y_{*}) \setminus \{y_{*}\}$ is reduced to examining $Q$ on each ${\mathcal U}_j^{\ell}$ for $1\le j \le L$ and $1\le \ell < \ell_{*}$.
\end{lemma} 
\begin{proof}
For $y_{*} + p^{\ell} u \in B_{p^{-1}}(y_{*})\setminus \{y_{*}\}$ with $|u|=1$,
we expand
\begin{equation}\label{P-expansion}
Q(y_{*} + p^{\ell} u) - Q(y_{*}) \ = \ \sum_{j\in {\mathcal E}} p^{L_j + j\ell} c_j u^j \ + \ c_k p^{k\ell} u^k \ + \ p^{(k+1)\ell} g(u)
\end{equation}
where $g \in {\mathbb Z}_p[{\rm x}]$, ${\mathcal E} = \{1\le j \le k-1 : Q^{(j)}(y_{*}) \not= 0 \}$, and for each
$j\in {\mathcal E}$, $Q^{(j)}(y_{*})/j! = p^{L_j} c_j$ for some unit $|c_j| = 1$.
When ${\mathcal E} = \emptyset$, \eqref{P-expansion} becomes
$$
Q(y_{*} + p^{\ell} u) - Q(y_{*}) \ = \  c_k p^{k\ell} u^k \ + \ p^{(k+1)\ell} g(u).
$$
Let $j_{*} := \min {\mathcal E}$ if ${\mathcal E} \not= \emptyset$ and $j_{*} := k$
if ${\mathcal E} = \emptyset$.

First, suppose ${\mathcal E} \not= \emptyset$. 
Let $\ell_{*}> 1$ be large so that for every $\ell \ge \ell_{*}$, we have 
\[
L_{j_{*}} + j_{*}\ell < L_j + j \ell \quad \text{for all} \quad j\not= j_{*},
\]
and $L_{j_{*}} + j_{*}\ell < k\ell$. Note that we can take $\ell_{*} = O_Q(1)$. 
Consider two elements $w = y_{*} + p^{\ell} u, \ z = y_{*} + p^m v \in {\mathcal U}_j^{\ge \ell_{*}}$ with $\ell,  m\ge \ell_{*}$ such that $Q(w) = Q(z)$. Then
$$
|Q(y_{*} + p^{\ell} u) - Q(y_{*})| \ = \ p^{-(L_{j_{*}} + j_{*}\ell)} \ = \ 
p^{-(L_{j_{*}} + j_{*}m)} \ = \ |Q(y_{*} + p^{m} v) - Q(y_{*})|,
$$
implying $\ell = m$. Furthermore, $Q(w) = Q(z)$ implies $F(u) = F(v)$, where
$$
F(u) \ := \ \frac{Q(y_{*} + p^{\ell} u) - Q(y_{*})}{p^{L_{j_{*}} + j_{*}\ell}} \ = \
c_{j_{*}} u^{j_{*}} \ + \ p^N \, h(u)
$$
and $h \in {\mathbb Z}_p[{\rm x}]$. Here $N$ can be taken large if $\ell_{*} > 1$ is large;
how large $N$ is needed will be determined later.

Next, suppose ${\mathcal E} = \emptyset$. Again, let $\ell_{*} > 1$ be large and consider two elements $w = y_{*} + p^{\ell} u, \ z = y_{*} + p^m v \in {\mathcal U}_j^{\ge \ell_{*}}$ with $\ell, m\ge \ell_{*}$ such that $Q(w) = Q(z)$. Then
$$
|Q(y_{*} + p^{\ell} u) - Q(y_{*})| \ = \ p^{-k\ell} \ = \ 
p^{- k m} \ = \ |Q(y_{*} + p^{m} v) - Q(y_{*})|,
$$
implying $\ell = m$. Furthermore, $Q(w) = Q(z)$ implies $F(u) = F(v)$ where
$$
F(u) \ := \ \frac{Q(y_{*} + p^{\ell} u) - Q(y_{*})}{p^{k \ell}} \ = \
c_{k} u^{k} \ + \ p^{\ell} \, g(u) 
$$
and $\ell$ can be taken large if $\ell_{*} > 1$ is large since $\ell \ge \ell_{*}$. In the ${\mathcal E} = \emptyset$ case, $\ell = N$.

 Hence $|u^{j_{*}}-v^{j_{*}}|<1$ 
which implies $|u-v|<1$ since $w,z \in {\mathcal U}_j$ and \eqref{1-1-again}. If $|j_{*}| = 1$,
then $|F'(v)| = |j_{*} c_{j_{*}} v^{j_{*} -1} + p^N h'(v)| = 1$ (recall $|u| = |v| =1$) and so Hensel's lemma implies $u = v$.  On the other hand, if $|j_{*}| = p^{-M}$ for some $M\ge 1$ and $u = u_0 + p \omega$, then
$$
F(u_0 + p \, \omega) - F(u_0) \ = \ c_{j_{*}} \Bigl(j_{*} p u_0^{j_{*}-1} \omega \ + \ \sum_{m=2}^{j_{*}}\binom{j_{*}}{m} p^m u_0^{j_{*}-m} \omega^{m} \Bigr)  \ + \ p^N \, {\tilde{h}}(\omega)
$$
for some ${\tilde{h}} \in {\mathbb Z}_p[{\rm x}]$.  Note that for each $2\le m \le j_{*}$,
\begin{equation}\label{j*choose}
\Bigl| \binom{j_{*}}{m} p^m u_0^{j_{*}-m} \Bigr| \le p^{-m-M+\nu_p(m!)}
\ \le \ p^{-m-M+m-2} \ = \ p^{-M-2}.
\end{equation}
by Lemma \ref{p-factorial} whenever $p\ge 3$. Therefore,  setting $N$ large enough such that  $N > M+1$, we deduce  
$$
\frac{F(u_0 + p \, \omega) - F(u_0)}{p^{M+1}} \ = \ a \, \omega \ + \ p g(\omega) \ \ {\rm for \ some \ unit} \ a \  {\rm and} \ \ g \in {\mathbb Z}_p [{\rm x}].
$$
Since $u = u_0 + p \omega$ and $v = u_0 + p \eta$ with $\omega, \eta \in {\mathbb Z}_p$, another application of Hensel's lemma gives $u=v$,
showing that in all cases, $Q$ is injective on ${\mathcal U}_j^{\ge \ell_{*}}$ for each $1\le j \le L$. 

Finally when $p=2$, we return to \eqref{j*choose} and note that $\big|\binom{j_{*}}{m}\big| \le 2^{-M+m-2}$ 
whenever $m\ge 3$. Hence, if $N> M+1$,
$$
G(\omega) :=\frac{F(u_0 + 2 \, \omega) - F(u_0)}{2^{M+1}} \ = \ a \, \omega \ + \ b \, \omega^2 \ + \ 2 g(\omega) \ \ {\rm for \ some \ units} \ a, \, b 
$$
and $g \in {\mathbb Z}_2 [{\rm x}]$. To apply Hensel's lemma and conclude $u = v$ as before, we need to
decompose ${\mathbb Z}_2 = \bigcup_{i\in\{0,1\}} W_i $ with 
$W_i=\{ \omega=\sum_{j\ge 0}\omega_j2^j \in {\mathbb Z}_2 : \omega_0 = i \}$. In fact, 
 $G$ is injective on each $W_i$. This gives that for every $1\le j \le L$,  
 we can write ${\mathcal U}_j^{\ge \ell_{*}}=\bigcup_{i\in\{0,1\}}\ {\mathcal U}_{j,i}^{\ge \ell_{*}}$, where 
$$
{\mathcal U}_{j,i}^{\ge \ell_{*}} = \bigl\{ y_{*} + 2^{\ell} u: {\rm for \ some} \ \ell \ge \ell_{*}, u=u_0+2\omega,\ {\rm with}\   |u_0| =1,  \omega\in W_i,  \ {\rm and} \ {\overline{u}} \in U_j \bigr\},
$$
such that 
$Q$ is injective on each ${\mathcal U}_{j,i}^{\ge \ell_{*}}$.
\end{proof}

\subsection{Proof of Proposition \ref{prop-1}}\label{subsection:c.1}
To prove Proposition \ref{prop-1}, 
we may assume that the coefficients of $P$ lie in ${\mathbb Z}_p$. Furthermore, we may assume $P(0) = 0$ and that one of the coefficients of $P$ is a unit. In fact, without loss of generality, we can replace $P$ with $cP$ for some constant $c \in {\mathbb Q}_p$.
Hence if $P(x) = \sum_{j=1}^d {a}_j x^j$, then each coefficient lies in ${\mathbb Z}_p$
and at least one coefficient ${a}_j$ has the property $|{a}_j| = 1$.

By induction, we will establish the following hypothesis for every integer $k\ge 1$.

\begin{Hk}
For every polynomial $Q \in {\mathbb Z}_p[{\rm x}]$ and for every open set $\Omega \subseteq {\mathbb Z}_p$ such that $|Q^{(k)}/k!| \equiv 1$ on $\Omega$, there exist a finite set ${\mathcal F}\subset {\mathbb Z}_p$ with $\# {\mathcal F} = O_Q(1)$ and an $O_Q(1)$ collection $\{V_j\}$ of open sets such that 

\begin{itemize}
\item[(i)] $\Omega \setminus {\mathcal F} \  = \ \bigcup_{j \le O_Q(1)} V_j$; and
\vskip 7pt
\item[(ii)] $Q$ is injective on each $V_j$.
\end{itemize}
    
\end{Hk}

Suppose now {\bf Hypothesis k} holds for every $k\ge 1$. For our polynomial $P$ of degree $d$ from Proposition \ref{prop-1} with the above normalisation ($P \in {\mathbb Z}[{\rm x}]$, $P(0) = 0$ and at least one coefficient is a unit), we claim that there is a $1\le k \le d$ such that $|P^{(k)}/k!| \equiv 1$ on ${\mathbb Z}_p$\footnote{This means $|P^{(k)}(x)/k!|=1 $ for all $x\in \Z_p$.}. In fact, if $|a_d| =1$, then $|P^{(d)}/d!| \equiv |a_d| = 1$ on ${\mathbb Z}_p$. If $|a_d|<1$ and $|a_{d-1}| = 1$, then $|P^{(d-1)}(x)/(d-1)!| = |a_{d-1} + d a_d x| \equiv 1$ on ${\mathbb Z}_p$. Repeating this process several times,   we can apply {\bf Hypothesis k} with $\Omega =  {\mathbb Z}_p$ to prove Proposition \ref{prop-1}.

It remains to prove {\bf Hypothesis k} for every $k\ge 1$.

\subsubsection{Proof of  Hypothesis 1}
We start with {\bf Hypothesis 1}. Let $\Omega \subseteq {\mathbb Z}_p$ be an open set such that $|P'| \equiv 1$ on $\Omega$. We may assume $\Omega \not= \emptyset$.  Consider
$$
{\mathcal G} \ := \ \bigl\{ y_0 \in {\mathbb N}_{\le p-1}: \Omega \cap B_{p^{-1}}(y_0) \not= \emptyset \bigr\} \ \ {\rm and \ note \ that} \ \ 
\Omega \ = \ \bigcup_{y_0 \in {\mathcal G}} (\Omega \cap B_{p^{-1}}(y_0)).
$$
For any $y = \sum_{j\ge 0} y_j p^j \in {\mathbb Z}_p$,  recall that ${\overline{y}} =  [y_0] \in {\mathbb F}_p$ denotes the residue field element such that $|y-y_0|<1$. For our polynomial $P(x) = \sum_{j=1}^d a_j x^j$, let ${\overline{P}}(x) := \sum_{j=1}^d {\overline{a_j}} x^j \in {\mathbb F}_p[{\rm x}]$ and note that ${\overline{P}}$ is not the zero polynomial in ${\mathbb F}_p[{\rm x}]$ by
our normalisation of $P$ that at least one coefficient is a unit. Therefore $\# ({\overline{P}})^{-1}(x) \le d$ for all $x \in {\mathbb F}_p$ and so we can partition
${\mathcal G} \ = \ U_1^\circ \cup U_2^\circ \cup \cdots \cup U_L^\circ$ with $L\le d$ so that (similar to \eqref{1-1-again})
\begin{equation}\label{1-1}
\forall 1\le j\le L, \ \ {\rm if} \ y_0', z_0' \in U_j ^\circ\ \ {\rm and} \ \ {\overline{P}}([y_0']) = {\overline{P}}([z_0']), \ {\rm then} \ y_0' = z_0'.
\end{equation}
Finally,  let $V_j := \bigcup_{y_0 \in U_j^\circ} B_{p^{-1}}(y_0)$ so that
$$
\Omega \ = \ \bigcup_{j=1}^L (\Omega \cap V_j)
$$
decomposes $\Omega$ into $O_P(1)$ open sets. 
For $y, z \in \Omega \cap V_j$, we have $y = y_0 + p u$ and $z = z_0 + p v$
where $y_0,z_0 \in U_j^\circ$. Also if $P(y) = P(z)$, then ${\overline{P}}({\overline{y}}) =
{\overline{P}}({\overline{z}})$ 
and so $|y - z|<1$ (that is, $y_0 = z_0$) by \eqref{1-1}. Also since $|P'(z)| = |P'(z_0)| = 1$, an application of Hensel's lemma shows $y = z$,
proving that $P$ is injective on $\Omega \cap V_j$. This establishes {\bf Hypothesis 1}.

\subsubsection{The inductive step}\label{subsubsection:The inductive step}
We now suppose that $k\ge 2$ and that {\bf Hypothesis ${\bf k}^{\prime}$} holds for all $1\le k'\le k-1$.
Let $\Omega \subset {\mathbb Z}_p$ be an open set such that $|P^{(k)}/k!| \equiv 1$ on $\Omega$. This is only possible if $k\le d = {\rm deg}\,P$.
We start as before and decompose 
$$
\Omega \ = \ \bigcup_{y_0 \in {\mathcal G}} (\Omega \cap B_{p^{-1}}(y_0)).
$$
For each $y_0 \in {\mathcal G}$, we note that $|P^{(k)}/k!| \equiv 1$ on $B_{p^{-1}}(y_0)$. Consider
$$
{\mathcal G'} \ := \ \bigl\{y_0 \in {\mathcal G} : |P^{(k-1)}(y_0)/(k-1)!| < 1 \bigr\}
$$
and note that if $y_0 \in {\mathcal G} \setminus {\mathcal G'}$, then
$|P^{(k-1)}/(k-1)!| \equiv 1$ on $B_{p^{-1}}(y_0)$. Therefore $|P^{(k-1)}/(k-1)!| \equiv 1$ on the open set
$$
V \ := \ \bigcup_{y_0 \in {\mathcal G}\setminus {\mathcal G'}} (\Omega \cap B_{p^{-1}}(y_0)).
$$
We can use {\bf Hypothesis k-1} for $P$ and $V$, giving us the desired decomposition of $V$ on which $P$ is injective. Hence matters are reduced to decomposing 
$$
V' \ := \ \bigcup_{y_0 \in {\mathcal G'}} (\Omega \cap B_{p^{-1}}(y_0)). 
$$

If $|k|<1$, then $p | k$, implying $p\le d$ and so $\# {\mathbb F}_p = p = O(1)$. In this case $\# {\mathcal G'} \le \# {\mathbb F}_p = O(1)$. 
If $|k|=1$
and $y_0 \in {\mathcal G'}$, then $|P^{(k-1)}(y_0)/(k-1)!| < 1$ and we can apply Hensel's lemma to $Q = P^{(k-1)}/(k-1)!$, noting $|Q'(y_0)| = 1$ since $|k|=1$, which gives us a unique ${\underline{y}} \in B_{p^{-1}}(y_0)$ such that $P^{(k-1)}({\underline{y}}) = 0$. Since ${\rm deg}(P^{(k-1)}) \le d-1$, we see that $\# {\mathcal G'} \le d-1$. Hence, in either case (whether $|k|<1$ or $|k|=1$), we have $\# {\mathcal G'} = O(1)$ and so matters are reduced to decomposing the open ball $B_{p^{-1}}(y_0)$ for each $y_0 \in {\mathcal G'}$.

When $|k|=1$, we have $B_{p^{-1}}(y_0) = B_{p^{-1}}({\underline{y}})$. Let $y_{*}$ denote either $y_0$ or
${\underline{y}}$, depending if $|k|<1$ or $|k|=1$, respectively. We will decompose $B_{p^{-1}}(y_{*})\setminus \{y_{*}\}$. We now apply Lemma \ref{poly-basic-decomp} to $P$ and $B_{p^{-1}}(y_{*})$,
noting $|P^{(k)}/k!| \equiv 1$ on $B_{p^{-1}}(y_{*})$ since $y_0 \in {\mathcal G'}$. This reduces matters to showing that $P$ is injective
on each ${\mathcal U}_j^{\ell}$ for $1\le j \le L$ and $1\le \ell < \ell_{*}$. 
We will decompose ${\mathcal U}_j^{\ell} = \bigl\{ y_{*} + p^{\ell} u:  |u| =1 \ {\rm and} \ {\overline{u}} \in U_j \bigr\}$ with the desired injectivity property of $P$.

For any $y_{*} + p^{\ell} u \in {\mathcal U}_j^{\ell}$, we recall \eqref{P-expansion} with $Q$ replaced by $P$, that is, 
$$
P(y_{*} + p^{\ell} u) - P(y_{*}) \ = \ \sum_{j\in {\mathcal E}} p^{L_j + j\ell} c_j u^j \ + \ c_k p^{k\ell} u^k \ + \ p^{(k+1)\ell} g(u).
$$
Let $m_{\ell} = \min\{L_j + j\ell: j=k \ {\rm or}\  j\in {\mathcal E}\}$ with $L_k=0$
so that
$$
 F(u) \ := \ \frac{P(y_{*} + p^{\ell} u) - P(y_{*})}{p^{m_{\ell}}} \ = \ \sum_{j\in J} b_j u^j \ + \ p \, {\tilde{g}}(u) \ \ {\rm for \ some} \ {\tilde{g}} \in {\mathbb Z}_p[{\rm x}],
$$
where $J \subset \{1,2,\ldots, k\}$ is nonempty and $|b_j| =1$ for all $j\in J$. 

\subsubsection{The case $|k| = 1$}\label{subsubsection:|k|=1}

In this case, $P^{(k-1)}(y_{*}) = 0$ and so $J\subset \{1,2,\ldots, k-2, k\}$. If $k \in J$, then
$$
F^{(k-1)}(u)/(k-1)!  \ = \ k \, b_k \, u \ + \ p \, {\tilde{g}}^{(k-1)}(u)/(k-1)!,
$$ 
implying $|F^{(k-1)}(u)/(k-1)!| \equiv 1$ on the open set $U := \{u\in\mathbb Z_p: |u| =1\}$. On the other hand, if $k \notin J$, set $j_{*} = \max J \le k-2$ so that 
$$
\bigl|F^{(j_{*})}(u)/j_{*}! \bigr| \ = \ \bigl| b_{j_{*}} +  p \, {\tilde{g}}^{(j_{*})}(u)/j_{*}!\bigr| \ \equiv \ 1
$$
on $U$. In either case, we can employ {\bf Hypothsis k-1} or {\bf Hypothsis ${\bf j_{*}}$}
to $F$ and $U$ which gives us the desired decomposition of $U$ (and hence ${\mathcal U}_j^{\ell}$) on which $F$ (and hence $P$) is injective. This completes the proof of Proposition \ref{prop-1} when $|k|=1$.

\subsubsection{The case $|k| < 1$.}\label{subsubsection:|k|<1}
We now turn our attention to the case $|k|<1$ where $y_{*} = y_0$. In this case, we adjust the definition of $m_{\ell}$ and define $m_{\ell}' = \min\{L_j + j\ell:  j \in {\mathcal E}\}$.

\smallskip

{\bf Case I}: $m_{\ell}' \le k\ell - 1$ (this implies that ${\mathcal E} \not= \emptyset$). As above, we consider
$$
F(u) \ := \ \frac{P(y_{0} + p^{\ell} u) - P(y_{0})}{p^{m_{\ell}'}} \ = \ \sum_{j\in J} b_j u^j \ + \ p \, {\tilde{g}}(u) \ \ {\rm for \ some} \ {\tilde{g}} \in {\mathbb Z}_p [{\rm x}]
$$
but now $J \subset \{1,2,\ldots, k-1\}$ is nonempty and $|b_j| =1$ for all $j\in J$. Again set $j_{*} = \max J \le k-1$ and note that
$$
\bigl|F^{(j_{*})}(u)/j_{*}! \bigr| \ = \ \bigl| b_{j_{*}} \ + \ p \, {\tilde{g}}^{(j_{*})}(u)/j_{*}!\bigr| \ \equiv \ 1
$$
on $U$, allowing us to invoke  {\bf Hypothsis ${\bf j_{*}}$} for $F$ and $U$, giving us the desired decomposition for ${\mathcal U}_j^{\ell}$ as before.

\smallskip

{\bf Case II}: $m_{\ell}' \ge k\ell$ in which case $m_{\ell} = k\ell$ (recall that
$m_{\ell} = \min\{L_j + j\ell: j=k \ {\rm or}\  j\in {\mathcal E}\}$ with $L_k=0$). In this case, we consider 
$$
F(u) \ := \ \frac{P(y_{0} + p^{\ell} u) - P(y_{0})}{p^{k\ell}} \ = \ \sum_{j\in {\mathcal E}} 
p^{L_j + \ell j - k\ell} c_j u^j \ + \ c_k u^k \ + \ \sum_{j\ge k+1} d_j p^{\ell(j-k)} u^j
$$
for $|c_j|=|c_k|=1$ and $|d_j|\le 1$,  and observe $|F^{(k)}(u)/k!| \equiv 1$.
Suppose that $|k| = p^{-t}$ for some $t\ge 1$ so that $k = a p^{t}$ for some unit $a$. Consider
$$
s \ := \ \begin{cases} L_{k-1} + \ell (k-1) - k \ell, & \text{if $k-1 \in {\mathcal E}$}\\
\infty, & \text{if $k-1 \notin {\mathcal E}$}. \end{cases}
$$

\smallskip

{\bf Subcase IIa}: $0\le s \le t-1$ (this implies $k-1 \in {\mathcal E}$). Here we observe that
\begin{equation}\label{F-IIa}
\frac{1}{p^{s}} \frac{F^{(k-1)}(u)}{(k-1)!} \ = \ c_{k-1} \ + \ \frac{k c_k}{p^{s}} u \ + \ 
\sum_{j\ge k+1} d_j  p^{-s} {\binom{j}{k-1}} p^{\ell(j-k)} u^{j-k+1}
\end{equation}
and for $j\ge k+1$,
$$
\Bigl|p^{-s} {\binom{j}{k-1}} p^{\ell(j-k)}\Bigr| \ = \ \Bigl|\frac{k}{p^s} \frac{p^{\ell(j-k)}}{(j-k+1)!} \frac{j!}{k!}\Bigr| \ \le \ p^{-(t-s)} p^{-\ell(j-k)} p^{\nu_p((j-k+1)!)},
$$
where $\nu_p((j-k+1)!) \le j-k$ by Lemma \ref{p-factorial}. Since $\ell\ge 1$,
$$
\Bigl|p^{-s} {\binom{j}{k-1}} p^{\ell(j-k)}\Bigr| \ \le \ p^{-1} p^{-(j-k)} p^{j-k} \ = \ p^{-1}.
$$
Since $|k p^{-s}| \le p^{-1}$ as well, this shows that for all $u \in {\mathbb Z}_p$, 
\begin{equation}\label{F-s}
\Bigl|\frac{1}{p^s} \frac{F^{(k-1)}(u)}{(k-1)!} \Bigr| \ = \ 1 \ \ {\rm and} \ \
|F^{(j)}(u)/j!| \le p^{-(j-k)} \  \ \text{for all} \ \ j\ge k+1.
\end{equation}
Next, we will decompose $U = \{u\in\mathbb Z_p: |u| =1\} = \bigcup_{m\le O(1)} U_m$ for some $U_m$,  which in turn will give
a decomposition of ${\mathcal U}_j^{\ell} = \bigcup_{m\le O(1)} {\mathcal U}_{j,m}^{\ell}$
where ${\mathcal U}_{j,m}^{\ell} = \{y_0 + p^{\ell} u \in {\mathcal U}_j: u \in U_m\}$.

For $A = O(1)$ (to be determined) and for every $u_{**} = u_0 + u_1 p + \cdots + u_{A-1} p^{A-1}$ with $u_0 \not= 0$, we define 
$U_{u_{**}} := \{u = u_{**} + p^A v : v \in {\mathbb Z}_p\}$. Since we are in the $|k|<1$ case (and hence $\# {\mathbb F}_p = O(1)$), this gives a decomposition of $U$ into $O(1)$ open sets. For $u_{**} + p^A \omega \in U_{u_{**}}$, we write
$F(u_{**} + p^A \omega) - F(u_{**})  =$ 
$$ 
\sum_{j=1}^{k-2} \frac{F^{(j)}(u_{**})}{j!} p^{Aj} \omega^j  +  \frac{F^{(k-1)}(u_{**})}{(k-1)!} p^{A(k-1)} \omega^{k-1}  + \frac{F^{(k)}(u_{**})}{k!} p^{Ak} \omega^{k} +  \sum_{j\ge k+1} d_j p^{Aj} \omega^j 
$$
$$
= \sum_{j=1}^{k-2} c_j p^{L_j + Aj} \omega^j + a p^{A(k-1) + s} \omega^{k-1} + b p^{Ak} \omega^k + \sum_{j\ge k+1} d_j p^{Aj} \omega^j,
$$
where $a,b$ and the $c_j's$ are all units, and $|d_j|\le 1$ for all $j\ge k+1$. Here we used \eqref{F-s} and the fact that
$|F^{(k)}/k!| \equiv 1$.

Let $m := \min\big\{\min\{L_j + Aj: 1\le j \le k-2\}, \ A(k-1)+s\big\}$ and note that
$m \le A(k-1) + s \le Ak-1$ for any $A > t$. With such an $O(1)$ choice for $A$, we consider
$$
H(\omega) \ := \ \frac{F(u_{**} + p^A \omega) - F(u_{**})}{p^m} \ = \
\sum_{j\in J} b_j \omega^j \ + \ p \, g(\omega)
$$
for some $g \in {\mathbb Z}_p[{\rm x}]$. Here $J \subset \{1,2,\ldots, k-1\}$ is nonempty and $|b_j| = 1$ for every $j\in J$. Let $j_{*} = \max J \le k-1$ and note that
$|H^{(j_{*})}(\omega)/j_{*}!| = |b_{j_{*}}| \equiv 1$. Hence   we can apply {\bf Hypothesis ${\bf j}_{*}$} (to $H$ and $U_{u_{**}}$) to give us our desired decomposition on which $P$ is injective.

\smallskip

{\bf Subcase IIb}: $s \ge t$ (this includes the case ${\mathcal E} = \emptyset$).
We begin with a couple preliminary observations about the function $F$ defined at the outset of {\bf Case II}: first, 
$$
\frac{F^{(k-1)}(u)}{k!} \ = \ \frac{1}{k} \, \frac{F^{(k-1)}(u)}{(k-1)!} \ = \ 
\frac{{\tilde{a}}}{p^{t}} \frac{F^{(k-1)}(u)}{(k-1)!},
$$
where $|{\tilde a}| = 1$. Furthermore, we claim $F^{(k-1)}/k! \in {\mathbb Z}_p[{\rm x}]$. In fact, when
$s < \infty$,
\begin{equation}\label{F-poly}
\frac{F^{(k-1)}(u)}{k!} \ = \ c_{k-1} \frac{p^{s-t}}{a} \ + \ c_k \, u \ + \ 
\sum_{j\ge k+1} d_j  \frac{j!}{k!} \frac{p^{\ell(j-k)}}{(j-k+1)!}  u^{j-k+1},
\end{equation}
and when $s = \infty$,
$$
\frac{F^{(k-1)}(u)}{k!} \ = \  c_k \, u \ + \ 
\sum_{j\ge k+1} d_j  \frac{j!}{k!} \frac{p^{\ell(j-k)}}{(j-k+1)!}  u^{j-k+1}.
$$
Since $\nu_p((j-k+1)!) \le j-k$ by Lemma \ref{p-factorial}, we see that
$p^{\ell(j-k)}/(j-k+1)! \in {\mathbb Z}_p$ for every $j\ge k+1$ and so
$F^{(k-1)}/k! \in {\mathbb Z}_p[{\rm x}]$.

With these observations, we begin our analysis of $F$ as in {\bf Case IIa} by decomposing ${\mathcal U}_j^{\ell} = \ \bigcup_{m\le O(1)} {\mathcal U}_{j,m}^{\ell}$ as before with respect to some large $A = O(1)$. As before, we write the difference
$F(u_{**} + p^A \omega) - F(u_{**})  =$ 
$$ 
\sum_{j=1}^{k-2} \frac{F^{(j)}(u_{**})}{j!} p^{Aj} \omega^j  +  \frac{F^{(k-1)}(u_{**})}{(k-1)!} p^{A(k-1)} \omega^{k-1}  + \frac{F^{(k)}(u_{**})}{k!} p^{Ak} \omega^{k} +  \sum_{j\ge k+1} d_j p^{Aj} \omega^j .
$$
Although it is still true that $|F^{(k)}/k!| \equiv 1$ as observed at the outset of {\bf Case II}, it is no longer the case that the analogue of \eqref{F-s} holds. We therefore
split into two further subcases.

{\bf b1}: 
$|F^{(k-1)}(u_{**})/k!| = |{\tilde{a}} p^{-t} F^{(k-1)}(u_{**})/(k-1)!| = 1$. This implies that
$$
F(u_{**} + p^A \omega) - F(u_{**})  =  \sum_{j=1}^{k-2} c_j p^{L_j + Aj} \omega^j + {\tilde{a}} p^{A(k-1) + t} \omega^{k-1} + b p^{Ak} \omega^k + \sum_{j\ge k+1} d_j p^{Aj} \omega^j,
$$
where ${\tilde{a}}, b$ and the $c_j's$ are all units. We now proceed exactly as in {\bf Subcase IIa} with $s$ replaced by $t$ to obtain the desired decomposition of each ${\mathcal U}_{j}^{\ell}$ on which our polynomial $P$ is injective.

{\bf b2}: $|F^{(k-1)}(u_{**})/k!| = |{\tilde{a}} p^{-t} F^{(k-1)}(u_{**})/(k-1)!| < 1$.

Consider $Q(u) := F^{(k-1)}(u)/k!$ which was observed in \eqref{F-poly} to belong to ${\mathbb Z}_p[{\rm x}]$. In this case, $|Q(u_{**})| < 1$ whereas $|Q'| = |F^{(k)}/k!| \equiv 1$ and so by Hensel's theorem, there is a unique ${\underline{u}}_{**}$ such that
$F^{(k-1)}({\underline{u}}_{**}) = 0$ and $|{\underline{u}}_{**} -u_{**}|<1$. We will decompose $B_{p^{-1}}({\underline{u}}_{**})\setminus \{{\underline{u}}_{**}\}$ into $O(1)$ open sets on each of which $F$ is injective. This will then give a decomposition of each ${\mathcal U}_{j}^{\ell}$ on which $P$ is injective, finishing the proof of Proposition \ref{prop-1}.

Here we apply Lemma \ref{poly-basic-decomp} to $F$ and $B_{p^{-1}}({\underline{u}}_{**})$, reducing matters
to examining $F$ on each ${\mathcal U}_j^{\ell} = \{{\underline{u}}_{**} + p^{\ell} u : |u|=1 \ {\rm and} \ 
{\overline{u}} \in U_j\}$ for each $1\le j \le L$ and $1\le \ell < \ell_{*}$. With respect to ${\underline{u}}_{**} + p^{\ell} u \in {\mathcal U}_j^{\ell}$, we expand
$$
F({\underline{u}}_{**} + p^{\ell} u) - F({\underline{u}}_{**}) \ = \ \sum_{j\in {\mathcal E}} p^{L_j + j\ell} c_j u^j \ + \ c_k p^{k\ell} u^k \ + \ \sum_{j\ge k+1} d_j p^{j\ell} u^j,
$$
where ${\mathcal E} = \{1\le j \le k-1 : F^{(j)}({\underline{u}}_{**}) \not= 0 \}$, and for each $j\in {\mathcal E}$,   $F^{(j)}({\underline{u}}_{**})/j! = p^{L_j} c_j$ for 
some unit $|c_j| = 1$. The important gain we have achieved in this $|k|<1$ case is that now
$F^{(k-1)}({\underline{u}}_{**}) = 0$ and so ${\mathcal E}\subset \{1\le j \le k-2\}$.

Also from \eqref{F-s}, for every $j\ge k+1$, we have
$$
|d_j| \ = \ |F^{(j)}({\underline{u}}_{**})/j!| \ \le \ p^{-(j-k)}
$$
and so
$$
F({\underline{u}}_{**} + p^{\ell} u) - F({\underline{u}}_{**}) \ = \ \sum_{j\in {\mathcal E}} p^{L_j + j\ell} c_j u^j \ + \ c_k p^{k\ell} u^k \ + \ \sum_{j\ge k+1} e_j p^{j\ell + (j-k)} u^j
$$
for some $e_j \in {\mathbb Z}_p$.

As before, we define
$m_{\ell} = \min\{L_j + j\ell: j=k \ {\rm or}\  j\in {\mathcal E}\}$ with $L_k=0$ so that
$$
 H(u) \ := \ \frac{F({\underline{u}}_{**} + p^{\ell} u) - F({\underline{u}}_{**})}{p^{m_{\ell}}} \ = \ \sum_{j\in J} b_j u^j \ + \ p \, {\tilde{g}}(u) \ \ {\rm for \ some} \ {\tilde{g}} \in {\mathbb Z}_p[{\rm x}],
$$
where $J \subset \{1,2,\ldots, k-2, k\}$ is nonempty and $|b_j| =1$ for all $j\in J$. 
If $k \notin J$, set $j_{*} := \max J \le k-2$ and note that $|H^{(j_{*})}(u)/j_{*}!|  \equiv|b_{j^{*}}|= 1$. Hence, we can apply {\bf Hypothesis ${\bf j}_{*}$} to $H$ and $U = \{u\in\mathbb Z_p: |u|=1\}$ to give us our desired decomposition.

This leaves us with the case $k \in J$ so that $m_{\ell} = k \ell$. In this case, a more precise formulation of $H$ is
$$
H(u) = \frac{F({\underline{u}}_{**} + p^{\ell} u) - F({\underline{u}}_{**})}{p^{k\ell}} = 
 \sum_{j\in {\mathcal E}} p^{L_j + (j-k)\ell} c_j u^j  +  c_k  u^k  +  \sum_{j\ge k+1} e_j \, p^{(j-k)(\ell+1)} u^j,
$$
defined on $U = \{u\in\mathbb Z_p: |u|=1\}$. Importantly, we have
$$
\frac{H^{(k-1)}(u)}{(k-1)!} \ = \ c_k k u \ + \ \sum_{j\ge k+1} e_j 
{\binom{j}{k-1}} p^{(j-k)(\ell +1)} u^{j-k+1}
$$
where
$$
{{\binom{j}{k-1}} }  \ = \ \frac{j(j-1)\cdots k}{(j-k+1)!} \ = \ k \, a_{j,k} 
p^{-\nu_p((j-k+1)!)} \ \ {\rm where} \ a_{j,k} \in {\mathbb Z}_p.
$$
Lemma \ref{p-factorial} implies
$$
\Bigl|{{\binom{j}{k-1}} } p^{(j-k)(\ell+1)}\Bigr|  \ \le \ |k| p^{-(j-k)(\ell+1) + \nu_p((j-k+1)!)} \ \le \ |k| p^{-(j-k)\ell}
\ < \ |k|,
$$ 
since $\ell \ge 1$ and $j\ge k+1$,  so we have 
\begin{equation}\label{H-k}
|H^{(k-1)}(u)/(k-1)!| \equiv |k| \ \ {\rm  for \ all} \ \ |u| = 1. 
\end{equation}
 For each $u_{*} \in \N_{\le p-1}\setminus \{0\}$, 
 since $p=O(1)$, 
 it suffices to decompose $B_{p^{-1}}(u_{*})\setminus \{u_{*}\}$
into $O(1)$ open sets on which $H$ is injective.

We now apply Lemma \ref{poly-basic-decomp} to $H$ and $B_{p^{-1}}(u_{*})$, reducing matters 
to examining $H$ on each ${\mathcal U}_j^{\ell} = \{u_{*} + p^{\ell} u : |u|=1 \ {\rm and} \ 
{\overline{u}} \in U_j\}$ for each $1\le j \le L$ and $1\le \ell < \ell_{*}$.
We expand
$$
H(u_{*} + p^{\ell} v) - H(u_{*}) \ = \ \sum_{j\in {\mathcal E}} p^{L_j + j\ell} c_j v^j \ + \ c_k p^{k\ell} v^k \ + \ \sum_{j\ge k+1} d_j p^{j\ell} v^j,
$$
but now ${\mathcal E} := \{1\le j \le k-1 : H^{(j)}(u_{*}) \not= 0\}$. From \eqref{H-k} and $|k|=p^{-t}$, we
see that $k-1 \in {\mathcal E}$ and $L_{k-1} = t$. By using the same argument as in {\bf Case I} above, matters are further reduced to $m_{\ell} = k\ell$ and examining
$$
G(v) \ := \ \frac{H(u_{*} + p^{\ell} v) - H(u_{*})}{p^{k\ell}} \ = \ \sum_{j\in {\mathcal E}} 
p^{L_j + \ell (j - k)} c_j v^j \ + \ c_k v^k \ + \ \sum_{j\ge k+1} d_j p^{\ell(j-k)} v^j
$$
for $|v|=1$. We are now in the same position as the outset of {\bf Case II} with $F$ replaced by $G$ (equivalently, $P$ replaced by $H$) and $y_0$ replaced by $u_{*}$. However now, $s := L_{k-1} + \ell(k-1 -  k) = t - \ell \le t-1$ since $\ell\ge 1$. Thus  we are in {\bf Subcase IIa}, which gives us the desired decomposition for $\{v\in\mathbb Z_p:|v|=1\}$ and hence $B_{p^{-1}}(u_{*})\setminus \{u_{*}\}$ on which $H$ is injective. 

This completes the proof of Proposition \ref{prop-1}.

\end{document}